\documentclass[11pt,oneside]{amsart}
\usepackage{amsmath,amsfonts,amsgen,amstext,amsbsy,amsopn,amsthm, amssymb}
\usepackage{multirow, verbatim, cancel, mathabx}
\usepackage{bbold}
\usepackage{enumitem}
\usepackage[T1]{fontenc}
\usepackage[colorlinks=true,hyperindex=true]{hyperref}
\usepackage{tikz}
\usetikzlibrary{calc}
\usepackage{verbatim}
\usetikzlibrary{arrows,chains,matrix,positioning,scopes}
\usepackage{mathrsfs}
\usepackage{color,soul}
\usepackage{xcolor}
\usepackage{tabulary}
\usepackage{pgfplots}
\usepackage{mathrsfs}
\usepackage{bbm}
\usepackage{dsfont}
\usepackage{bm}
\usepackage{mathtools}
\usepackage{float}

\DeclarePairedDelimiter\floor{\lfloor}{\rfloor}
\usepackage{arydshln}
\usepackage{color,soul}
\usepackage{subcaption}
 \usepackage{soul,xcolor}
 
 \usepackage[font={small,it}]{caption}

\usepackage{float}
\usepackage{flafter}

\usepackage{array}
\newcolumntype{R}[1]{>{\raggedleft\arraybackslash}p{#1}}

\usepackage{tabu} 
\usepackage{booktabs}

\usepackage[mathscr]{euscript}

\DeclareSymbolFont{rsfs}{U}{rsfs}{m}{n}
\DeclareSymbolFontAlphabet{\mathscrsfs}{rsfs}

\makeatletter

\usepackage{tkz-fct} 
\usepackage{xcolor}

\usepackage{tkz-fct} 
\usepackage{pgfplots}
\usepackage{setspace}

\usetikzlibrary{matrix}
\usepackage{tikz}

\numberwithin{equation}{section}

\newcounter{smallarabics}

\newcounter{smallroman}
\newenvironment{romanenumerate}
{\begin{list}{{\normalfont\textrm{(\roman{smallroman})}}}
  {\usecounter{smallroman}\setlength{\itemindent}{0cm}
   \setlength{\leftmargin}{5ex}\setlength{\labelwidth}{4ex}
   \setlength{\topsep}{0.75\parsep}\setlength{\partopsep}{0ex}
   \setlength{\itemsep}{0ex}}}
{\end{list}}

\newcommand{\ben}{\begin{romanenumerate}}  
\newcommand{\een}{\end{romanenumerate}}  

\newtheorem{theoreme}{theorem }[section]
\newtheorem{theorem}[theoreme]{Theorem}
\newtheorem{proposition}[theoreme]{Proposition}
\newtheorem{Lemma}[theoreme]{Lemma}
\newtheorem{conjecture}[theoreme]{Conjecture}

\newtheorem{corollary}[theoreme]{Corollary}
\newtheorem{remark}{Remark}[section]
\newtheorem{example}[theoreme]{Example}

\newcolumntype{L}{>{\centering\arraybackslash}m{3cm}}

\newcommand\nn\nonumber
\renewcommand\leq\varleq
\renewcommand\geq\vargeq

 \newcommand{\R}{\mathbb{R}}
 \newcommand{\N}{\mathbb{N}}
\newcommand{\Z}{\mathbb{Z}} \newcommand{\C}{\mathbb{C}}
 \newcommand{\F}{\mathcal{E}}
\newcommand{\E}{\mathcal{E}}

\newcommand{\grad}{\nabla}

\renewcommand{\i}{\mathrm{i}}

\renewcommand{\F} {\mathcal{F}}

\renewcommand{\epsilon}{\varepsilon}

\newcommand\blue[1]{\textcolor{blue}{#1}}
\newcommand\red[1]{\textcolor{red}{#1}}

\usepackage{tikz-cd}

\usepackage{xcolor}
\usepackage{dcolumn}
\usepackage{tabu}

\newcolumntype{A}{D{.}{.}{2.3}}

\usepackage{tikz,tkz-tab}
\usepackage{pgfplots}
\pgfplotsset{compat=1.11}
\usetikzlibrary{positioning}
\makeatletter
      \def\@setcopyright{}
      \def\serieslogo@{}
      \makeatother
\begin{document}

\author{Marc-Adrien Mandich}
      \address{Independent researcher, Jersey City, 07305, NJ, USA}
	\email{marcadrien.mandich@gmail.com}
   

   \title[LAP for discrete Schr\"odinger operator]{Thresholds and more bands of A.C.\ spectrum for the Molchanov--Vainberg Schr\"odinger operator with a more general long range condition}

   \begin{abstract}
   The existence of absolutely continuous (a.c.) spectrum for the discrete Molchanov--Vainberg Schr\"odinger operator $D+V$ on $\ell^2(\Z^d)$, in dimensions $d\geq 2$, is further investigated for potentials $V$ satisfying the long range condition $n_i(V-\tau_i ^{\kappa}V)(n) = O(\ln^{-q}(|n|))$ for some $q>2$, $\kappa \in \N$ even, and all $1 \leq i \leq d$, as $|n| \to \infty$. $\tau_i ^{\kappa} V$ is the potential shifted by $\kappa$ units on the $i^{\text{th}}$ coordinate. In this article \textit{finite} linear combinations of conjugate operators are constructed. These lead to more bands of a.c.\ spectrum being found. However, the new bands of a.c.\ spectrum are justified mainly by graphical evidence because the coefficients of the linear combinations are obtained by numerical polynomial interpolation. At the same time, an infinitely countable set of thresholds is rigorously identified (these will be defined exactly in the article). We conjecture that the spectrum of $D+V$ in dimension 2 is void of singular continuous spectrum, and that consecutive thresholds constitute endpoints of a band of a.c.\ spectrum.
    
   \end{abstract}

%

 

\maketitle
\hypersetup{linkbordercolor=black}
\hypersetup{linkcolor=blue}
\hypersetup{citecolor=blue}
\hypersetup{urlcolor=blue}
\tableofcontents

\begin{figure}[htb]
  \centering
 \includegraphics[height=7cm,width=12cm]{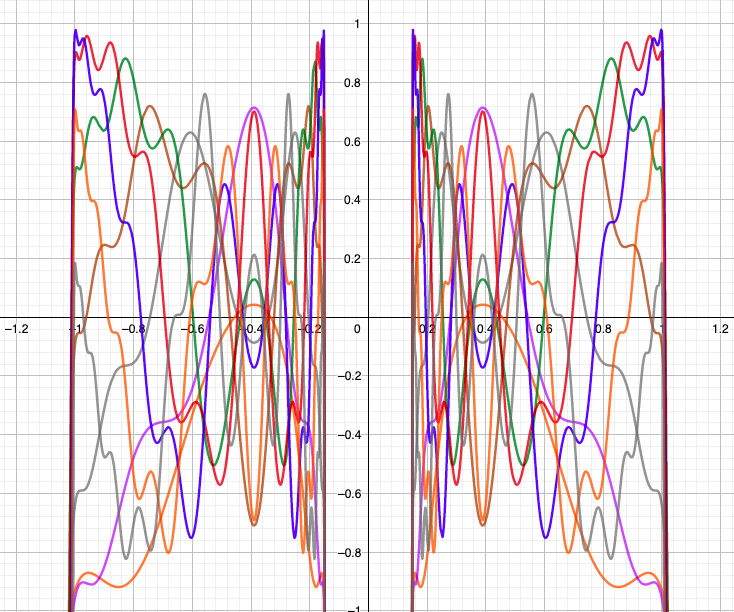}
\caption{This article is about constructing linear combinations of these polynomials, the $g_{j\kappa}^E(x)$}
\label{fig:cover}
\end{figure}

\section{Introduction}

This article is concerned with the study of spectral properties of the Molchanov--Vainberg Laplacian operator when perturbed by a class of long range perturbations. The Molchanov--Vainberg Laplacian, introduced in \cite{MV}, is a type of discrete Schr\"odinger operator on the lattice $\Z^d$ and can be used to model quantum phenomena in media with discrete postions such as crystals, or more general media by means of discretisation. This article is a sequel to \cite{GM2}, and parallels \cite{GM3} which covers the same topic but for the standard Laplacian, which we will denote by $\Delta$. The Molchanov--Vainberg Laplacian acts in the Hilbert space $\mathscr{H} := \ell^2(\Z^d)$ as follows:
\begin{equation}
\begin{aligned}
\label{def:std}
D[d] = D &:= \prod_{i=1}^d \Delta_i, \quad \text{where } \  \Delta_i := 2^{-1} (S_i + S_i ^*) \\
&= 2^{-d} \sum_{\nu \in \Gamma} S_1^{\nu_1} S_2^{\nu_2} \cdot \cdot \cdot S_d^{\nu_d}, \quad \text{where } \ \Gamma :=  \{-1,1\}^d, \ \nu = (\nu_1, \ldots, \nu_d).
\end{aligned}
\end{equation}
Here $S_i = S_i^1$ and $S^*_i = S^{-1}_i$ are the shifts to the right and left respectively on the $i^{th}$ coordinate. So $(S_i ^{\pm 1} u)(n) = u(n_1,\ldots,n_i \mp 1,\ldots,n_d)$ for $u \in \mathscr{H}$, $n=(n_1,\ldots, n_d) \in \Z^d$. Set $|n|^2 = n_1^2 +...+n_d^2$. Let $\sigma(\cdot)$ denote the spectrum of an operator. A Fourier transformation shows that the spectra of $\Delta_i$ and $D$ are purely absolutely continuous (a.c.), $\sigma(\Delta_i) \equiv [-1,1]$ and $\sigma(D) = [-1,1]$. 

Let $V$ model a discrete electric potential and act pointwise, i.e.\ $(Vu)(n) = V(n) u(n)$, for $u \in \mathscr{H}$. We always assume $V$ is real-valued and goes to zero at infinity. Thus the essential spectrum of $D+V$ equals $\sigma(D)$. Let $\N$ and $\N^*$ be the positive integers, including and excluding zero respectively. Let $\N_o$ and $\N_e$ be the \textit{odd} and \textit{even} positive integers respectively. Fix $\kappa \in \N^*$. The potential shifted by $\pm \kappa$ units is defined by 
\begin{equation*}
(\tau_i ^{\pm \kappa} V) u(n) := V(n_1,\ldots, n_i \mp \kappa,\ldots n_d)u(n), \quad \forall  1 \leq i  \leq d.
\end{equation*} 
As in \cite{GM2} and \cite{GM3}, we study potentials $V$ satisfying a non-radial condition of the form 
\begin{equation}
\label{generalLR condition}
n_i (V - \tau_i ^{\kappa} V)(n) = O(g(n)), \ \text{as} \ |n| \to \infty, \quad \forall  1 \leq i  \leq d,
\end{equation}
 where $g(n)$ is a radial function which goes to zero at infinity at an appropriate rate, e.g.\ $g(n) =  \ln ^{-q}(|n|+1)$, $q>2$. We refer to \cite{GM2} for some examples of Schr\"odinger operators that satisfy \eqref{generalLR condition}. 
 
The two discrete Laplacians, $\Delta$ and $D$, are identical in dimension 1, and so we will only consider $D$ for $d \geq 2$. Also, $\Delta$ and $D$ are isomorphic in dimension 2, and so an interesting aspect is to compare the dimension 2 results for $D+V$ with $V$ satisfying \eqref{generalLR condition} with $2\kappa$ with those for $\Delta+V$ in \cite{GM3} with $V$ satisfying \eqref{generalLR condition} with $\kappa$. The reason the comparison should not be done $\kappa$ to $\kappa$, but rather $2\kappa$ to $\kappa$, is detailed in the isomorphism in \cite[section V]{GM2}. In dimensions $\geq 3$ the two Laplacians are not isomorphic and so this article also presents results specific to $D$. It apppears that many, if not all, of the interesting observations or results for $\Delta+V$ mentioned in \cite{GM3} have corresponding ones for $D+V$. There is a lot of truth in saying that to write this article it is a simple matter of mindlessly replacing $a+b$ by $a \times b$, $a-b$ by $a \div b$, $a \times b$ by $a ^b$, and $a \div b$ by $a ^{1/b}$ in the key formulas in \cite{GM3}...but then it's not completely true! There are some notable and non-trivial differences that make the study of $D+V$ worthwhile and interesting. Because this article parallels \cite{GM3} which covers the exact same topic but for $\Delta$ we may be more brief in the conceptual exposition at times. 

Let $I \subset \R$ be a closed interval, let $I_{\pm}:= \{ z \in \C_{\pm}: \mathrm{Re} (z) \in I\}$, $\C_\pm:=\{z\in \C, \pm\mathrm{Im}(z)>0\}$. The limiting absorption principle (LAP) is a statement about the boundary values of the maps 
\begin{equation}
\label{LaP_generic}
I_{\pm} \ni z \mapsto (D+V-z)^{-1}. 
\end{equation}
One important implication of the LAP on an interval $I$ is the absence of s.c.\ spectrum for $D+V$ on that set. This article produces such type of results. In Mourre theory, which was extensively studied in \cite{Mo1}, \cite{Mo2} and \cite{ABG}, one of the strategies to obtain a LAP \eqref{LaP_generic} on an interval $I \subset \sigma(D+V)$ depends roughly on the ability to prove two key estimates. The first estimate is a \textit{strict Mourre estimate} for $D$ with respect to some self-adjoint conjugate operator $\mathbb{A}$ on this interval, that is to say, $\exists \gamma >0$ such that 
\begin{equation}
\label{mourreEstimate123}
1_I(D) [D, \i \mathbb{A}] _{\circ} 1_{I} (D) \geq \gamma 1_{I} (D),
\end{equation}
where $1_I(D)$ is the spectral projection of $D$ on $I$, and $[ \cdot , \i \mathbb{A} ] _{\circ}$ initially defined on the compactly supported sequences is the extension of the commutator between an operator and $\mathbb{A}$ to a bounded operator on $\mathscr{H}$ (this definition suffices for this article). The second estimate is one involving $V$, and according to a more recent perspective on the theory (see \cite{GM1}) is such as
\begin{equation}
\label{compactcomm}
 \ln^p(1+|n|) \cdot [V, \i \mathbb{A}]_{\circ}  \cdot \ln^p(1+|n|) \quad is \ a \ compact \ operator \ on \ \mathscr{H} \ for \ some \ p >1.
\end{equation}
Denote the position operators $(N_i u)(n) := n_i u(n)$. These are required to specify our choice of $\mathbb{A}$. As in \cite{GM3} we consider a (\textit{finite}) linear combination of conjugate operators of the form
\begin{equation}
\label{LINEAR_combinationA}
\mathbb{A} = \sum _{j \geq 1} \rho_{j\kappa} \cdot A_{j \kappa}, \quad \rho_{j\kappa} \in \R, \quad A_{j \kappa} := \sum_{1 \leq i \leq d} A_i (j,\kappa),
\end{equation}
where each $A_i (j,\kappa)$, initially defined on compactly supported sequences, is the closure in $\mathscr{H}$ of :
\small
\begin{equation}
\label{generatorDilations_0111}
A_{i}(j,\kappa) := \frac{1}{2\i} \bigg [  \frac{j \kappa}{2}(S_i ^{j \kappa} + S_i^{-j \kappa}) + (S_i ^{j \kappa} - S_i^{-j \kappa}) N_i \bigg ]  = \frac{1}{4\i} \bigg[  (S_i ^{j \kappa} - S_i^{-j \kappa})N_i + N_i (S_i^{j \kappa} - S_i^{-j \kappa}) \bigg].
\end{equation} 
\normalsize
Each $A_{j\kappa}$ is self-adjoint in $\mathscr{H}$ by an adaptation of the case $(j,\kappa) = (1,1)$, and so $\mathbb{A}$ is self-adjoint, at least whenever it is a finite sum. The reason choice \eqref{LINEAR_combinationA} is relevant is that 
$$[V, \i A_{j\kappa} ]_{\circ} = \sum_{1 \leq i \leq d}  (4\i)^{-1} \left( (V - \tau_i ^{j \kappa} V) S_i ^{j \kappa} - (V - \tau_i ^{-j\kappa}V)S_i ^{-j\kappa} \right) N_i + \text{hermitian conjugate},
$$
and so \eqref{generalLR condition} implies \eqref{compactcomm}, again, at least when $\mathbb{A}$ is a finite sum and $g(n) = \ln^{-q}(1+|n|)$, $q>2$. While the frequencies of the $A_{i}(j,\kappa)$ are in sync with the long range frequency decay of $V$, the coefficients $\rho_{j \kappa}$ need to be chosen so that \eqref{mourreEstimate123} holds and this is a challenge. We partition the spectrum of $D$ into two sets : $\boldsymbol{\mu}_{\kappa}(D)$ and $\boldsymbol{\Theta}_{\kappa}(D)$. $\boldsymbol{\mu}_{\kappa}(D)$ are energies $E \in \sigma(D)$ for which there is a self-adjoint linear combination (finite or infinite) of the form \eqref{LINEAR_combinationA}, an interval $I \ni E$ and $\gamma > 0$ such that the Mourre estimate \eqref{mourreEstimate123} holds. $\boldsymbol{\Theta}_{\kappa}(D)$ are energies $E \in \sigma(D)$ for which there is no self-adoint linear combination (finite or infinite) of the form \eqref{LINEAR_combinationA}, no interval $I \ni E$ and no $\gamma >0$ such that \eqref{mourreEstimate123} holds. By definition $\sigma(D)$ is a disjoint union of $\boldsymbol{\mu}_{\kappa}(D)$ and $\boldsymbol{\Theta}_{\kappa}(D)$. From Mourre theory $\boldsymbol{\mu}_{\kappa}(D)$ is an open set and so $\boldsymbol{\Theta}_{\kappa}(D)$ is closed. In this article, including title and abstract, energies in $\boldsymbol{\Theta}_{\kappa}(D)$ are referred to as \textit{thresholds}. \textit{This definition depends on the modeling assumption of $\mathbb{A}$}. Theorem \ref{lapy305} below highlights the usefulness of the sets $\boldsymbol{\mu}_{\kappa}(D)$. Let $\sigma_p(D+V)$ be the point spectrum of $D+V$. Let $\langle \mathbb{A} \rangle := \sqrt{1+\mathbb{A}^* \mathbb{A} }$. 
\begin{theorem}
\label{lapy305} 
Let $q > 2$, $\kappa \in \N^*$ be such that $\limsup \left( |V(n)|, | n_i (V - \tau_i^{\kappa} V)(n) | \right) = O \left(  \ln ^{-q} (|n|) \right )$, as $|n| \to \infty$ and $\forall 1 \leq i \leq d$.

Let $E \in \boldsymbol{\mu}_{\kappa}(D) \setminus \sigma_p(D+V)$. Let $\mathbb{A} = \sum_j \rho_{j \kappa} A_{j \kappa}$ be a finite sum such that \eqref{mourreEstimate123} holds in a neighorhood of $E$. Then there is an open interval $I$, $I \ni E$, such that
\begin{enumerate}
\item $\sigma_p(D+V) \cap I$ is at most finite (including multiplicity),
\item$\forall p >1/2$ the map $I_{\pm} \ni z \mapsto (D+V-z)^{-1} \in \mathscr{B}\left( \mathcal{K}, \mathcal{K}^* \right)$ extends to a uniformly bounded map on $I$, with $\mathcal{K} = L^2_{1/2,p}(\mathbb{A}) = \left \{ \psi \in \mathscr{H} : \| \langle \mathbb{A} \rangle ^{1/2}  \ln ^{p} (\langle \mathbb{A} \rangle) \psi \| < \infty \right \}$,
\item The singular continuous spectrum of $D+V$ is void in $I$.
\end{enumerate}
\end{theorem}
This theorem can be refined, see \cite{GM2} and references therein. Operator regularity is a necessary and important topic underlying our methodology. According to the standard literature, the regularity $D, V \in C^{1,1}(\mathbb{A})$, or adequate variations thereof, are required. It is clear that $D$, $V$ belong to $C^1(A_{j\kappa})$ for $1 \leq j \leq N < \infty$, and this implies $D$, $V \in C^1(\sum_{1 \leq j \leq N} \rho_{j\kappa} A_{j\kappa})$. Although the $C^1(\mathbb{A})$ compliance falls short of the required regularity, since $C^1(\mathbb{A}) \subset C^{1,1}(\mathbb{A})$, we do not expect regularity to be a problem in this article. \\

\noindent \textbf{Problem of article :} determine for as many energies $E \in \sigma(D)$ if $E \in \boldsymbol{\mu}_{\kappa}(D)$ or $E \in \boldsymbol{\Theta}_{\kappa}(D)$. \\

In this article we will only consider \textit{even} values of $\kappa \in \N^*$. For $\kappa$ odd it is an open problem to decide if $\boldsymbol{\mu}_{\kappa}(D)$ is empty or not -- note that this is harder to prove than \cite[Lemma IV.3]{GM2}, because there the linear combination $\mathbb{A} = $ \eqref{LINEAR_combinationA} consisted of just the first term.

As a result of \cite{GM2}, $\boldsymbol{\mu}_{\kappa=2}(D) \supset [-1,1] \setminus \{-1,0,1\}$ in any dimension $d$ and this was proved choosing $\mathbb{A} = A_1 = \sum_{i=1}^d A_i (1,2)$. Actually, equality holds and this is easy to prove (see Lemma \ref{lemSUMcosINTRO_00} below). So onwards we only consider $\kappa \geq 4$, $\kappa \in \N_e$. For these values of $\kappa$ incomplete results for $\boldsymbol{\mu}_{\kappa}(D)$ and $\boldsymbol{\Theta}_{\kappa}(D)$ were obtained in \cite{GM2} by using $\mathbb{A} = A_{\kappa} = \sum_{i=1}^d A_i (1,\kappa)$. Note that this corresponds to \eqref{LINEAR_combinationA} with $\rho_{j\kappa} = 1$ if $j=1$ and $\rho_{j \kappa} = 0$ if $j \geq 2$. Table \ref{tab:table1012sup} displays the intervals already determined (numerically) to belong to $\boldsymbol{\mu}_{\kappa}(D)$ for $2 \leq \kappa \leq 8$ (cf. \cite[Tables XV and XVII] {GM2}). In this article we continue to determine $\boldsymbol{\mu}_{\kappa}(D)$ and $\boldsymbol{\Theta}_{\kappa}(D)$ for $d \geq 2$, and $\kappa \geq 4$, even.

The overall high level strategy is the same as in \cite{GM3}. It is :
\begin{enumerate}
\item Fix a dimension $d$ and a $\kappa \geq 4$, $\kappa \in \N_e$. (we really only treat $d=2$ ; $d=3$ very briefly). 
\item Determine as many threshold energies in $\boldsymbol{\Theta}_{\kappa}(D)$ as possible. To do this we use the same idea as in \cite{GM3}. This involves solving systems of equations, which yield threshold energies $\E= \prod_{1 \leq l \leq d} x_l$ and their decomposition into coordinate-wise energies $\vec{x} = (x_1, ..., x_d)$. These are key in the next step.  
\item Pick 2 consecutive threshold energies $\E_{i_1}$ and $\E_{i_2}$ determined in the previous step (consecutive means that there aren't any other thresholds between $\E_{i_1}$ and $\E_{i_2}$), and try to construct a conjugate operator $\mathbb{A}$ of the form \eqref{LINEAR_combinationA} such that for \textit{every} $E \in (\E_{i_1}, \E_{i_2})$, there is an interval $I \ni E$ and $\gamma>0$ such that \eqref{mourreEstimate123} holds. \textit{The $\mathbb{A}$ is the same for every $E \in (\E_{i_1}, \E_{i_2})$}. To determine the coefficients $\rho_{j \kappa}$, we perform polynomial interpolation. 
\end{enumerate}  
Rigorous proofs for the existence of the thresholds in step (2) are produced. As for step (3), the polynomial interpolation is implemented numerically. Thus the coefficients $\rho_{j \kappa}$ are found numerically, and these are then used to plot a functional representation of \eqref{mourreEstimate123}. This is the evidence we share to see that strict positivity is in fact obtained.

We now describe the idea to get thresholds. Let $U_{\kappa}$ be the Chebyshev polynomials of the second kind of order $\kappa$. As 
$[D, \i A_{j\kappa}]_{\circ} = \sum_{i=1} ^d D \Delta_i ^{-1} (1-\Delta_i^2) U_{j\kappa - 1} (\Delta_i)$ and the $\Delta_i$ are self-adjoint commuting operators we may apply functional calculus. To this commutator associate the polynomial 
\begin{equation}
\label{def:gE_intro}
[-1,1]^d \ni \vec{x} \mapsto g_{j\kappa} (\vec{x}) :=  \sum_{i=1} ^d \left( \prod_{1 \leq l \leq d} x_l \right) x_i ^{-1} (1-x_i^2) U_{j\kappa - 1} (x_i) \in \R, \quad \vec{x} = (x_1, ..., x_d).
\end{equation}
If the linear combination of conjugate operators is $\mathbb{A} = \sum_{j \geq1} \rho_{j\kappa} \cdot A_{j\kappa}$, set
\begin{equation}
\label{def:GGE_intro}
[-1,1]^d \ni \vec{x} \mapsto G_{\kappa} (\vec{x}) :=  \sum_{j \geq 1} \rho_{j\kappa} \cdot g_{j\kappa} (\vec{x}) \in \R.
\end{equation}
$G_{\kappa}$ is a functional representation of $[D, \i \mathbb{A}]_{\circ}$. Consider the constant energy $E \in \sigma(\Delta)$ surface 
\begin{equation}
\label{constE_MV}
S_E := \left \{ \vec{x} \in [-1,1]^d : E =  \prod_{1 \leq l \leq d} x_l \right \}.
\end{equation} 
By functional calculus and continuity of the function $G_{\kappa}$ we have $E \in \boldsymbol{\mu}_{\kappa} (D)$ iff $ \left.G_{\kappa} \right|_{S_E} >0$. \\

\noindent \textit{Definition of $\boldsymbol{\Theta}_{0,\kappa}(D)$}. $E \in \boldsymbol{\Theta}_{0,\kappa}(D)$ iff $\exists \ \vec{x} := (x_1, ..., x_d) \in S_E$ such that $g_{j \kappa} (\vec{x}) = 0$, $\forall j \in \N^*$. 

\noindent If $\vec{x}$ is such a solution, then for any choice of coefficients $\rho_{ j \kappa} \in \R$, \eqref{def:GGE_intro} $= G_{\kappa} (\vec{x})  =0$.

\noindent \textit{Definition of $\boldsymbol{\Theta}_{m,\kappa}(D)$, $m \in \N^*$}. $E \in \boldsymbol{\Theta}_{m,\kappa}(D)$ iff there are $(\vec{x}_q)_{q=0} ^{m} := (x_{q,1}, ...,x_{q,d})_{q=0}^m \subset S_E$, and $(\omega_q)_{q=0}^{m-1} \subset \R$, $\omega_q \leq 0$ (crucial), $\forall \ 0 \leq q \leq m-1$, such that 
\begin{equation}
\label{special relationship_intro}
g_{j \kappa} (\vec{x}_m) = \sum_{q=0} ^{m-1} \omega_q \cdot g_{j \kappa} (\vec{x}_q), \quad \forall j \in \N^*.
\end{equation} 
If the $\vec{x}_q$ are such a solution, then for any choice of coefficients $\rho_{ j \kappa} \in \R$,
\begin{equation}
\label{special G relationship_intro}
G_{\kappa} (\vec{x}_m) = \sum _{j\geq 1} \rho_{j \kappa} \cdot g_{j  \kappa} (\vec{x}_m) = \sum _{q=0} ^{m-1} \omega_q \sum _{j\geq 1} \rho_{j \kappa} \cdot g_{j \kappa} (\vec{x}_q) = \sum _{q=0} ^{m-1} \omega_q \cdot G_{\kappa} (\vec{x}_q).
\end{equation}
If $\boldsymbol{\Theta}_{m,\kappa}(D) \cap \boldsymbol{\mu}_{\kappa}(D)$ was non-empty the lhs of \eqref{special G relationship_intro} would be strictly positive whereas the rhs of \eqref{special G relationship_intro} would be non-positive. An absurdity. Thus :
\begin{Lemma}
Fix $d \geq 1$, $\kappa \geq 1$. Then $\boldsymbol{\Theta}_{m,\kappa}(D) \subset \boldsymbol{\Theta}_{m,\alpha\kappa}(D) \subset \boldsymbol{\Theta}_{\alpha\kappa}(D)$, $\forall m \in \N$, and $\forall \alpha \in \N^*$.
\end{Lemma}
Simply because no counterexamples were found, we actually conjecture :

\begin{conjecture}
\label{conjecture000}
Fix $d=2$, $\kappa \geq 4$, $\kappa \in \N_e$. $\cup_{m \geq 0} \boldsymbol{\Theta}_{m,\kappa}(D) = \boldsymbol{\Theta}_{\kappa}(D)$.
\end{conjecture}

I don't understand $d=3$ nearly as well to submit a Conjecture like \ref{conjecture000} for it. It turns out it is very easy to find threshold energies in $\boldsymbol{\Theta}_{0,\kappa}(D)$. We prove :
\begin{Lemma}
\label{lemSUMcosINTRO} 
$\forall \ d,\kappa \in \N_e$, $\boldsymbol{\theta}_{0,\kappa} (D) := \left\{\prod_{1 \leq q \leq d} \cos( j_q \pi / \kappa) : (j_1,...,j_d) \in \{0,...,\kappa \}^d \right\} \subset \boldsymbol{\Theta}_{0,\kappa} (D)$.
\end{Lemma}
\begin{remark} 
This lemma supports the conjectures in relation to the band endpoints in Table \ref{tab:table1012sup}.
\end{remark}
Perhaps equality in Lemma \ref{lemSUMcosINTRO} holds; this is an open question.
Thresholds $\boldsymbol{\theta}_{0,\kappa} (D)$ in Lemma \ref{lemSUMcosINTRO} were already found in \cite{GM2}. Here we prove $\boldsymbol{\Theta}_{m,\kappa}(D) \neq \emptyset$, $\forall m \geq 0$, $\forall d \geq 2$, $\kappa \geq 4$ even. Thus, there are infinitely many thresholds for $d \geq 2$, $\kappa \geq 4$ even. This is a remarkable difference with the case of the dimension 1, or the case of $\kappa=2$ in any dimension (see \cite{GM2}) :

\begin{Lemma}
\label{lemSUMcosINTRO_00} 
Let $(d,\kappa) \in \{1\} \times \N_e \cup \N^* \times \{2\}$. Then $\boldsymbol{\theta}_{0,\kappa} (D) = \boldsymbol{\Theta}_{0,\kappa} (D) = \boldsymbol{\Theta}_{\kappa} (D)$.
\end{Lemma}

The interpolation setup of step (3) above is analogous to that in \cite{GM3}. We review it for clarity. Suppose $\E_{i_1}$ and $\E_{i_2}$ are consecutive thresholds, with $\E_{i_1} \in \boldsymbol{\Theta}_{m_1,\kappa}(D)$ and $\E_{i_2} \in \boldsymbol{\Theta}_{m_2, \kappa}(D)$ for some $m_1, m_2 \in \N$. Suppose the coordinate-wise energies are $(\vec{x}_q)_{q=0}^{m_1} \subset S_{\E_{i_1}}$ and $(\vec{y}_r)_{r=0}^{m_2} \subset S_{\E_{i_2}}$. Recall we have the assumption that the conjugate operator $\mathbb{A}$ is the \textit{same} $\forall E \in (\E_{i_1},\E_{i_2})$. Thus, while we want $G_{\kappa} > 0$ on $(\E_{i_1},\E_{i_2})$, a continuity argument implies that $G_{\kappa}$ is at best non-negative at the endpoints $\E_{i_1}$ and $\E_{i_2}$, due to \eqref{special G relationship_intro}. Also by continuity, $G_{\kappa}(\vec{x})$ must be a local minimum whenever  $G_{\kappa}(\vec{x}) = 0$ and $\vec{x}$ is an interior point of $S_{\E_{i_1}}$ or $S_{\E_{i_2}}$. Let $\mathrm{int}(\cdot)$ be the interior of a set. Thus we require :
\begin{equation}
\begin{cases}
\label{interpol_intro}
& G_{\kappa}(\vec{x}_q) = 0, \quad \forall 0 \leq q \leq m_1, \quad and \quad \grad G_{\kappa} (\vec{x}_q) = 0, \quad for \ \vec{x}_q \in \mathrm{int}(S_{\E_{i_1}}) \quad [\text{left}],\\
& G_{\kappa}(\vec{y}_r) = 0, \quad \forall 0 \leq r \leq m_2, \quad and \quad \grad G_{\kappa} (\vec{y}_r) = 0, \quad for \ \vec{y}_r \in \mathrm{int}(S_{\E_{i_2}})  \quad [\text{right}].
\end{cases}
\end{equation}
By \eqref{special G relationship_intro} conditions $G_{\kappa}(\vec{x}_{m_1}) = 0$ and $G_{\kappa}(\vec{y}_{m_2}) = 0$ are redundant. Constraints \eqref{interpol_intro} set up a system of linear equations to be solved for the coefficients $\rho_{j \kappa}$, i.e.\ we have polynomial interpolation. In order for the computer to numerically solve the system, multiples of $\kappa$, $\Sigma := \{ j_1 \kappa, j_2 \kappa, j_3 \kappa, ..., j_{\ell} \kappa \}$, need to be assumed. We choose $\Sigma$ essentially by trial and error but prioritize lower order polynomials to keep things as simple as possible. In other words, we loop over sets $\Sigma$ until we find an appropriate $\mathbb{A}$. Section \ref{sectionCHOICE} provides an illustration of $G_{\kappa}^E$ when the index set $\Sigma$ is inappropriately selected. Of course, if our assumption that the $\mathbb{A}$ is the same for all $E \in (\E_{i_1},\E_{i_2})$ is valid, then constraints \eqref{interpol_intro} are necessary but not necessarily sufficient in order to find an appropriate $\mathbb{A}$. As a rule of thumb, it makes sense to look for a conjugate operator $\mathbb{A}$ for which the number of coefficients $\rho_{j\kappa}$ corresponds to the number of constraints in \eqref{interpol_intro}. But unlike in \cite{GM3}, here we have not investigated if there is an appropriate conjugate operator $\mathbb{A}$ for which the number of non-zero coefficients $\rho_{j\kappa}$ must be greater than the number of constraints in \eqref{interpol_intro}. 


Our problem is harder as $d,\kappa$ increase. So we focus mostly on the dimension 2. Some of those results will carry over to $d \geq 3$. As for $\kappa$ we mostly limit the numerical illustrations and evidence to a handful of values. \textit{We always restrict our analysis to positive energies}, because $\boldsymbol{\mu}_{\kappa}(D) = - \boldsymbol{\mu}_{\kappa}(D)$, by Lemma \ref{Lemma_symmetryDelta}.

\textit{Until otherwise specified, we now focus exclusively on the \underline{dimension 2}}. For $\kappa \geq 4$, let
\begin{equation}
\label{setsJJ}
J_2 = J_2 (\kappa) := \left( \cos^2(\pi / \kappa), \cos(\pi / \kappa) \right), \quad J_1 = J_1(\kappa) := \left(\cos(\pi / \kappa), 1 \right).
\end{equation}
By Lemma \ref{lemSUMcosINTRO}, $\inf J_2, \sup J_2 = \inf J_1, \sup J_1 \in \boldsymbol{\theta}_{0,\kappa}(D)$. In \cite{GM2} we proved $J_1 \subset \boldsymbol{\mu}_{\kappa}(D)$, $\forall \kappa$ even. As for $J_2$ it was identified (numerically for the most part) as a gap between 2 bands of a.c.\ spectrum, but this was based on \eqref{LINEAR_combinationA} with $j=1$ only, see Table \ref{tab:table1012sup}. Let $\E_0 = \E_0 (\kappa) := \sup J_2$.

\begin{theorem} 
\label{thm_decreasing energy general}
Fix $\kappa \geq 4$, $\kappa \in \N_e$. There is a strictly decreasing sequence of energies $\{ \E_n \}_{n=0} ^{\infty} = \{ \E_n (\kappa) \}_{n=0} ^{\infty}$, which depends on $\kappa$, such that $\{ \E_n \} \subset J_2 \cap \boldsymbol{\Theta}_{\kappa}(D)$ and $\E_n \searrow \inf J_2 = \cos ^2(\pi / \kappa)$. Also, $\E_{2n-1}$ and $\E_{2n} \in \boldsymbol{\Theta}_{n, \kappa}(D)$, $\forall n \geq 1$.
\end{theorem}

The $\E_n$ are complicated numbers but exact solutions are sometimes attainable, see e.g.\ \eqref{exact_sols_k4} for exact solutions for $\kappa=4$. After graphing some numerical solutions for $\E_{2n}$, $1 \leq n \leq 100$, for $\kappa \in \{4,6\}$, we propose the same conjecture as in the case of $\Delta$ on the rate of convergence:

\begin{conjecture}
\label{conjecture12}
Let $\{ \E_n \}$ be the sequence in Theorem \ref{thm_decreasing energy general}. $\E_n - \inf J_2 = c(\kappa)/n^2 + o(1/n^2)$, $\forall \kappa \geq 4$, $\kappa \in \N_e$, where $c(\kappa)$ means a constant depending on $\kappa$.
\end{conjecture}

In section \ref{section_gen1} we state two Theorems and a Conjecture generalizing Theorem \ref{thm_decreasing energy general}. 

Unfortunately I was not successful in rigorously proving a Mourre estimate on any new interval. Nonetheless, based on numerical evidence given in sections \ref{appli_scheme_k2} and \ref{appli_scheme_k3} it looks like $J_2$ is the simplest of the gaps identified in \cite{GM2} to understand. Therefore we conjecture: 

\begin{conjecture}
\label{conjecture22}
Fix $\kappa \geq 4$, $\kappa \in \N_e$. Let $\{ \E_n \}$ be the sequence in Theorem \ref{thm_decreasing energy general}. For each interval $( \mathcal{E}_n , \mathcal{E}_{n-1} )$, $n \geq 1$, $\exists$ a conjugate operator $\mathbb{A}(n) = \sum_{q=1} ^{N(n)} \rho_{j_q \kappa} (n) A_{j_q \kappa}$, $A_{j_q \kappa} = \sum_{1 \leq i \leq 2} A_i (j_q,\kappa)$, such that the Mourre estimate \eqref{mourreEstimate123} holds with $\mathbb{A}(n)$, $\forall E \in ( \mathcal{E}_n , \mathcal{E}_{n-1} )$. $\mathbb{A}(n)$ is typically not unique. It can be chosen so that $N(n) = 2n$. In particular, $\{\E_n\} = J_2 \cap \boldsymbol{\Theta}_{\kappa}(D)$.
\end{conjecture}
$\kappa=4$ is the only value of $\kappa$ for which the closure of $J_2 \cup J_1$ equals $\sigma(D) \cap [0,1]$. Thus, if Conjecture \ref{conjecture22} is true, our problem is fully solved in the case of $\kappa=4$ (in dimension 2). But for $\kappa \geq 6$, Theorem \ref{thm_decreasing energy general} and Conjecture \ref{conjecture22}, together with the already existing results recorded in Table    \ref{tab:table1012sup}, do not paint a complete picture. For example, the above discussion does not address the situation on the interval $\simeq(0.25,0.5064)$, for $\kappa=6$. We make some progress in that direction, but things get even more complicated. In addition to \eqref{setsJJ}, for $\kappa \geq 6$ set 
$$J_3 = J_3 (\kappa) := \left( \cos(2 \pi / \kappa), \cos ^2(\pi / \kappa) \right).$$
$J_3$, $J_2$, $J_1$ are adjacent intervals. We looked only very briefly into the strange phenomenon regarding $J_3$, see section \ref{sec_Conjecture J3}. If the case of $D$ mirrors that of $\Delta$ it is plausible to conjecture: 

\begin{conjecture}
\label{conjecture_k344_2d}
For any $\kappa \geq 6$, $\kappa \in \N_e$, $J_3(\kappa) \subset \boldsymbol{\mu}_{\kappa}(D)$.
\end{conjecture}


In section \ref{section J3a} we recycle the proof of Theorem \ref{thm_decreasing energy general} to prove the existence of thresholds below $J_3$:

\begin{theorem} 
\label{thm_decreasing energy general_k3}
Fix $\kappa \geq 6$, $\kappa \in \N_e$. There is a strictly increasing sequence of energies $\{ \F_n \}_{n=1} ^{\infty}$ which depends on $\kappa$, such that $\{ \F_n \} \subset (\cos(\pi/ \kappa) \times \cos(2\pi/ \kappa), \cos(2\pi/ \kappa)) \cap \boldsymbol{\Theta}_{\kappa}(D)$ and $\F_n \nearrow \inf J_3 = \cos(2\pi/ \kappa)$. $\F_{2n-1}$, $\F_{2n} \in \boldsymbol{\Theta}_{n, \kappa}(D)$, $\forall n \geq 1$.
\end{theorem}

Figure \ref{fig:T3, increasing to 0.5} depicts the solutions $\F_1$, $\F_2$ and $\F_3$ for $\kappa=6$. Unfortunately I was only able to accurately numerically compute a couple of solutions $\F_n$ and so no conjecture on the rate of convergence of $\F_n$ will be formulated. In section \ref{section_gen2} two Theorems and a Conjecture generalizing Theorem \ref{thm_decreasing energy general_k3} for $\{ \F_n \}$ are mentioned.

If we work on the spectrum of $D$ starting from the middle (energy 0), and then move outwards, note that we expect 
$$J_{\text{last}} = J_{\text{last}} (\kappa) := \left(0,\cos^2((\kappa/2-1)\pi / \kappa)\right)$$
to belong to $\boldsymbol{\mu}_{\kappa}(D)$ (see Table \ref{tab:table1012sup} and especially \cite[Table XV]{GM2}). It therefore remains to better understand the nature of the spectrum on $ [\sup J_{\text{last}} (\kappa), \inf J_3 (\kappa)]$, minus the bands of a.c.\ spectrum that were already identified there.

Hopefully it will become clear from our examples and constructions that there are \textit{many more} thresholds $\in  [\sup J_{\text{last}} (\kappa), \inf J_3 (\kappa)]$ for $\kappa \geq 6$ in addition to the sequence $\{ \F_n \}$. For example, more such thresholds are given in Figures \ref{fig:example_2nd_gap}, \ref{fig:example_2nd_gap2} and \ref{fig:example_2nd_gap22} for $\kappa =6$, and we expect a bunch more to lie there. Here are open questions we find interesting :

$\bullet$ Of the thresholds $\in \boldsymbol{\theta}_{0,\kappa}(D)$, which ones are accumulation points, as a subset of $\boldsymbol{\Theta}_{\kappa}(D)$ ?

$\bullet$ Are there accumulation points $ \in \boldsymbol{\Theta}_{\kappa}(D) \setminus \boldsymbol{\theta}_{0,\kappa}(D)$ ?

$\bullet$ What are the rates of convergence to the accumulation points $ \in \boldsymbol{\Theta}_{\kappa}(D)$ ?

$\bullet$ Are there infinitely many accumulation points within $\boldsymbol{\Theta}_{\kappa}(D)$ ?

$\bullet$ Is there an interval $I \subset \sigma(D)$ for which $\boldsymbol{\Theta}_{\kappa}(D)$ is dense in $I$ ?

We do conjecture however (as in the case of $\Delta$):


\begin{conjecture}
\label{conjecture6666}
Fix $\kappa \geq 2$, $\kappa \in \N_e$. $\cup_{m \geq 0} \boldsymbol{\Theta}_{m,\kappa}(D)$ and $\boldsymbol{\Theta}_{\kappa}(D)$ are countable sets. 
\end{conjecture}


As far as the a.c.\ spectrum is concerned, we did not investigate the existence of $\boldsymbol{\mu}_{\kappa}(D)$ in $[\sup J_{\text{last}} (\kappa), \inf J_3(\kappa)]$, apart from what was already known. If the case of $D+V$ resembles that of $\Delta+V$, one might suspect this is the part of the spectrum that is much harder to crack. Challenges in trying to decode $\boldsymbol{\mu}_{\kappa}(D)$ were given in the introduction of \cite{GM3} and they apply here too. In spite of the challenges we do have an overall conjecture for the dimension 2:

\begin{conjecture}
\label{conjecture33}
Fix $d=2$, $\kappa \in \N_e$. Let $\E_{i_1}, \E_{i_2} \in \boldsymbol{\Theta}_{\kappa}(D)$ be two consecutive thresholds -- meaning that there aren't any other thresholds in between $\E_{i_1}$ and $\E_{i_2}$. Then there is a (finite?) linear combination $\mathbb{A} = \sum_{j=1} ^{N} \rho_{j \kappa} A_{j \kappa}$ such that the Mourre estimate \eqref{mourreEstimate123} holds with $\mathbb{A}$ for every energy $E \in ( \E_{i_1}, \E_{i_2} )$. In particular, in light of Theorem \ref{lapy305}, $\sigma_p(D+V)$ is locally finite on $( \E_{i_1}, \E_{i_2} )$, whereas the singular continuous spectrum of $D+V$ is void.
\end{conjecture}

\textit{We are done discussing $d=2$}. In higher dimensions we have only 1 general result: thresholds in dimension $d$ generate thresholds in dimension $d+1$, via scaling. Recall notation \eqref{def:std}. 
\begin{Lemma}
\label{shift_threshold_d}
$\{ \pm \cos(\frac{j \pi}{\kappa}) : 0 \leq j \leq \frac{\kappa}{2}-1 \} \times \boldsymbol{\Theta}_{m,\kappa}(D[d]) \subset \boldsymbol{\Theta}_{m,\kappa}(D[d+1])$, $\forall d \geq 1$, $\kappa \in \N_e$, $m \in \N$.
\end{Lemma}
For $(d,\kappa) \in \N^* \times \{2\}$, the inclusion in Lemma \ref{shift_threshold_d} is in fact equality, but perhaps there are values of $\kappa \geq 4$, $\kappa \in \N_e$, for which the inclusion in Lemma \ref{shift_threshold_d} is strict. 
Lemma \ref{shift_threshold_d} generalizes \cite[Lemma IV.2]{GM2}. 

We briefly treat the problem in \textit{dimension 3}. There is still considerable work to be done just to understand the case $\kappa=4$, especially on the interval $(\cos^3(\pi /4), \cos^2(\pi /4))$. \textit{Theorem \ref{thm3d_k2} and Conjecture \ref{conjecture11_3d} below are for $D$ in \underline{dimension $3$}}.
\begin{theorem}
\label{thm3d_k2} Fix $\kappa = 4$. We have :
\begin{itemize}
\item $(\cos(\pi/4),1) \subset \boldsymbol{\mu}_{\kappa}(D)$ (proved in \cite{GM2}). 
\item $0,\cos^3(\pi/4), \cos^2(\pi/4), \cos(\pi/4), 1 \in \boldsymbol{\theta}_{0,\kappa}(D)$ (Lemma \ref{lemSUMcosINTRO}).
\item Let $\{ \E_n = \E_n (\kappa=4) \}$ be the sequence in Theorem \ref{thm_decreasing energy general}. Applying Lemma \ref{shift_threshold_d} gives :
\item $\{\E_n\} \subset \left(\cos^2(\pi /4), \cos(\pi /4) \right) \cap \boldsymbol{\Theta}_{\kappa}(D)$, with $\E_n \searrow \cos^2(\pi /4)$, 
\item $\{\E_n \times \cos(\pi/4)\} \subset \left(\cos^3(\pi /4), \cos^2(\pi /4) \right) \cap \boldsymbol{\Theta}_{\kappa}(D)$, with $\E_n \times \cos(\pi/4) \searrow \cos^3(\pi /4)$.
\end{itemize}
\end{theorem}

Similarly to the case of $\Delta$ and $\kappa=2$, graphical evidence also suggests the following conjecture, although it is quite mysterious and surprising to me how and why it happens :
\begin{conjecture}
\label{conjecture11_3d}
Fix $\kappa =4$. Let $\{ \E_n = \E_n (\kappa=4) \}$ be the sequence in Theorem \ref{thm_decreasing energy general}. For each interval $( \mathcal{E}_n, \mathcal{E}_{n-1})$, $n \geq 1$,  the Mourre estimate \eqref{mourreEstimate123} holds with $\mathbb{A}(n) = \sum_{1 \leq q \leq N(n)} \rho_{j_q \kappa} (n) A_{j_q \kappa}$, $A_{j_q \kappa} = \sum_{1 \leq i \leq 3} A_i (j_q,\kappa)$, $\forall E \in ( \mathcal{E}_n, \mathcal{E}_{n-1})$, where the coefficients $\rho_{j_q \kappa}(n)$ are \textit{exactly} those used in the 2-dimensional case, see Conjecture \ref{conjecture22}. In particular $\{\E_n\}_{n=1}^{\infty} = J_2 \cap \boldsymbol{\Theta}_{\kappa}(D)$. 
\end{conjecture}

Conjecture \ref{conjecture11_3d} may extend to $\kappa \geq 6$, but we have not looked into it. The sequence in $(\cos^3(\pi /4), \cos^2(\pi /4))$ in Theorem \ref{thm3d_k2} is the only knowledge we have about this interval.

We conclude the introduction with a few comments. 

The formula for the $\omega_q$'s in this article are more involved than the one for the standard Laplacian in \cite{GM3}. This leads to a richer set of assumptions or ways these can be negative, see sections \ref{geo_construction} and \ref{revisit}.

In this article we discuss another type of threshold $\in \boldsymbol{\Theta}_{m,\kappa}(\Delta)$ which can happen in dimension 2 when a so-called \textit{alignment condition} is fulfilled, see section \ref{revisit}. This is an improvement over \cite{GM3} because, although they do appear in some graphs in \cite{GM4}, they were swept under the rug (especially the math behind them). In concordance with Conjecture \ref{conjecture22}, we believe these types of thresholds occur somewhere in $[\sup J_{\text{last}} (\kappa), \inf J_3(\kappa)]$. Unfortunately I was not successful in deriving a general formula for the $\omega_q$'s for these thresholds, see section \ref{revisit}. Furthermore, it is very likely that these types of thresholds occur only for $\kappa \geq 6$, but I don't have a proof for it.

Another improvement over \cite{GM3} is that in this article we better articulate the assumptions around systems \eqref{conjecture_system_conj_intro} and \eqref{conjecture_system_conj_intro_even} which are used to find thresholds in dimension 2. We believe this gives more clarity to the exposition. As a result however, the definition of the set $\mathfrak{T}_{n,\kappa}$ given here doesn't quite match the one given in \cite{GM3}, see section \ref{geo_construction}.

The thresholds found in this article raise the question about possible eigenvalues embedded in the continuous spectrum. I am not aware of an analysis of properties of such eigenfunctions for $D+V$, even for a long-range $V$ satisfying \eqref{generalLR condition} with $\kappa=1$. It appears to be an open question if this can be done using the Mourre estimate as in \cite{FH} (see also \cite{Ma}), or another technique.

We would also like to remind the reader that both $D$ and $\Delta$ converge to the continuous Laplacian in $\R^d$ in the norm resolvent sense, see \cite[Appendix A]{GM2} and \cite{NT}. This is another point that makes the study of $D$ worthwhile.

Finally, the reader is invited to consult the introduction of \cite{GM3} where additional relevant comments can be found and totally apply to this article too.

\noindent \textbf{Acknowledgements :} It is a pleasure to thank my former thesis advisor Sylvain Gol\'enia for conversations and Vojkan Jak\v{s}i\'c for encouraging me to study the Molchanov--Vainberg Laplacian. I also want to thank my great friend Laurent Beauregard for generously sharing programming ideas and always being there ready to pitch in.

\section{Basic properties and lemmas for the Chebyshev polynomials}

Let $T_{n}$ and $U_{n}$ be the Chebyshev polynomials of the first and second kind respectively of order $n \in \N$. They are defined by the formulas 
\begin{equation}
\label{def_T_U}
T_n(\cos(\theta)) = \cos(n \theta), \quad U_{n-1}( \cos(\theta)) = \sin(n \theta) / \sin(\theta), \quad \theta \in [-\pi,\pi], n \in \N.
\end{equation}
The parity of the polynomials $T_n$ and $U_n$ is the same as the parity of $n$. As we'll mainly use the Chebyshev polynomials $T_4$ and $T_6$ to illustrate our article we give their expressions :
\begin{equation}
\label{cheby_23}
\begin{aligned}
T_4(x) &= 8 x^4 - 8 x^2 + 1 \quad \text{and} \quad  T_6(x) &= 32 x^6 - 48 x^4 + 18 x^2 -1.
\end{aligned}
\end{equation}
It may be useful to be aware that $T_n(x) = T_n(y)$ can be factored into a product of straight lines or ellipses. For instance :
\begin{equation}
\label{cheby_2334}
\begin{aligned}
T_4(x) = T_4(y) &\Leftrightarrow (x-y)(x+y) (x^2+y^2-1)=0, \quad \text{and} \\
T_6(x) = T_6(y) & \Leftrightarrow (x - y) (x + y) (-3 + 4 x^2 - 4 x y + 4 y^2) (-3 + 4 x^2 + 4 x y + 4 y^2) =0.
\end{aligned}
\end{equation} 
The roots of $U_{n-1}$ are $\cos(l \pi / n)$, $1 \leq l \leq n-1$. We'll absolutely need a commutator $[ \cdot, \cdot]$ for \textit{functions}. For functions $f,g$ of real variables $x,y$, let 
\begin{equation}
\label{def_comm}
[f(x), g(y) ] := f(x) g(y) - f(y) g(x).
\end{equation}
\begin{remark}
The quantity $[f(x), g(y) ] / (x-y)$ is called \textit{Bezoutian} in the literature. Alternatively, \eqref{def_comm} ressembles a \textit{Wronskian}. Note this commutator satisfies $[f(x), g(y) ] = [g(y), f(x) ]$.
\end{remark}

Lemmas \ref{Tcos}, \ref{lemma_cos_alpha_beta}, \ref{variationsTk} and Corollaries \ref{corollaryEquiv} and \ref{basic_lemma} given below were already cited and proved in \cite{GM3}. They will play the same important role in this article ; in particular the corollaries are at the heart of our search for thresholds. 

\begin{Lemma}
\label{Tcos}
For $x,y \in [-1,1]$, $T_{\kappa}(x) = T_{\kappa}(y)$ if and only if $T_{\alpha \kappa}(x) = T_{\alpha \kappa}(y)$ for all $\alpha \in \N^*$. 
\end{Lemma}

\begin{Lemma}
\label{lemma_cos_alpha_beta}
Fix $\kappa \in \N^*$. If $\cos(\kappa \theta) = \cos(\kappa \phi)$ then 
$$ \sin(\kappa \phi ) \sin(2\kappa \theta) = \sin(\kappa \theta) \sin(2\kappa \phi) \Rightarrow \sin(\alpha \kappa \phi ) \sin(\beta \kappa \theta) = \sin(\alpha \kappa \theta) \sin(\beta \kappa \phi), \forall \alpha, \beta \in \N^*.$$
\end{Lemma}

\begin{corollary} 
\label{corollaryEquiv}
Let $\kappa \in \N^*$, $\kappa \geq 2$ be given. If $x,y \in \R$ are such that $U_{\kappa-1}(x)$, $U_{\kappa-1}(y) \neq 0$, then $T_{\kappa} (x) = T_{\kappa} (y) \Leftrightarrow [U_{\alpha \kappa-1}(x) , U_{\beta \kappa-1}(y) ] = 0, \forall \alpha, \beta \in \N^*$.
\end{corollary}

\begin{corollary}
\label{basic_lemma}
Let $\kappa \in \N^*$, $\kappa \geq 2$ be given. Let $x,y \in [-1,1]$. Then $[U_{\alpha \kappa-1}(x) , U_{\beta \kappa-1}(y) ]=0$ for all $\alpha, \beta \in \N^*$ if and only if $U_{\kappa-1}(x)=0$, or $U_{\kappa-1}(y)=0$, or $T_{\kappa}(x) = T_{\kappa}(y)$.
\end{corollary}

We also exploit the variations of $T_{\kappa}$ : 

\begin{Lemma}
\label{variationsTk}
Fix $\kappa \geq 1$. $T_{\kappa}([-1,1]) = [-1,1]$. $T_{\kappa}(1) = 1$, $T_{\kappa}(-1) = (-1)^{\kappa}$. The local extrema of $T_{\kappa}$ in $[-1,1]$ are located at $\cos (j \pi / \kappa)$, $0 \leq j \leq \kappa$. On $(\cos(j \pi / \kappa), \cos((j-1) \pi / \kappa))$, $j \in \{ 0, ..., \kappa\}$, $T_{\kappa}$ is strictly increasing if $j$ is odd and strictly decreasing if $j$ is even. 
\end{Lemma}

Finally, two other identities we'll exploit are :
\begin{equation}
\label{identity_derivative}
\frac{d}{dx} T_{\kappa}(x) = \kappa U_{\kappa-1}(x), \quad \text{and} \quad \frac{d}{dx} U_{\kappa-1}(x) = \frac{\kappa T_{\kappa}(x) - x \cdot U_{\kappa-1}(x)}{x^2-1}.
\end{equation}
A difference between this article and \cite{GM3} is that there the proofs didn't require the use the second identity in \eqref{identity_derivative}.

\section{Functional representation of the strict Mourre estimate for $D$ wrt.\ $\mathbb{A}$}
\label{stdLaplacianMourre}

Let $\mathcal{F} : \mathscr{H} \to L^2([-\pi,\pi]^d,d\xi)$ be the Fourier transform 
\begin{equation}
\normalsize
(\mathcal{F} u)(\xi) :=  (2\pi)^{-d/2} \sum \limits_{n \in \Z^d} u(n) e^{\i n \cdot \xi}, \quad  \xi=(\xi_1,\ldots,\xi_d).
\label{FourierTT}
\end{equation}
The commutator between $D$ and $A_{j\kappa}$, computed against compactly supported sequences, is $[D, \i A_{j\kappa}] =  \mathcal{F}^{-1} \left[  \sum_{i=1} ^d \sin(\xi_i) \sin(j\kappa \xi_i)\prod_{l \neq i} \cos(\xi_l)   \right] \mathcal{F}  = \sum_{i=1} ^d D \Delta_i ^{-1}(1-\Delta_i^2) U_{j\kappa-1} (\Delta_i)$, $\forall j\in\N^*$.
So $[D, \i A_{j\kappa}]$ extends to a bounded operator $[D, \i A_{j\kappa}]_{\circ}$. Let 
\begin{equation}
\label{def_m}
m(x) := 1-x^2.
\end{equation}
Fixing $E \in \sigma(D)$ allows us to remove the $x_d$ variable in \eqref{def:gE_intro}. We will often opt for that convention. Thus, consider the polynomial $g_{j\kappa}  ^E : [-1,1]^{d-1} \mapsto \R$,
\begin{equation}
\label{def:gE}
g_{j\kappa} ^E (x_1,...,x_{d-1}) :=  \sum_{i=1} ^{d-1} \frac{E}{x_i} m(x_i) U_{j\kappa - 1} (x_i) + \prod_{i=1} ^{d-1} x_i \cdot m\left(E/\prod_{i=1} ^{d-1} x_i \right) U_{j\kappa - 1} \left(E/ \prod_{i=1} ^{d-1} x_i \right) .
\end{equation}

\begin{Lemma}
\label{lemma_rootsy}
The roots of $m(x) U_{j \kappa-1}(x)$ are $\{\cos(l \pi / (j\kappa)) : 0 \leq l \leq j \kappa \}$. The intersection over $j \in \N^*$ of the latter set is $\{ \cos(l \pi / \kappa) : 0 \leq l \leq \kappa \}$ and these are roots of $m(x) U_{j \kappa-1}(x)$, $\forall j \in \N^*$. 
\end{Lemma}
If the linear combination of conjugate operators is $\mathbb{A} = \sum_{j \geq 1} \rho_{j\kappa} \cdot A_{j\kappa}$ set $G_{\kappa} ^E : [-1,1]^{d-1} \mapsto \R$,
\begin{equation}
\label{def:GGE}
G_{\kappa} ^E (x_1,...,x_{d-1}) :=  \sum_{j \geq 1} \rho_{j\kappa} \cdot g_{j\kappa} ^E (x_1,...,x_{d-1}) .
\end{equation}
\eqref{def:gE} and \eqref{def:GGE} are basically the same thing as \eqref{def:gE_intro} and \eqref{def:GGE_intro} but localized in energy $E$. Of course, $G_{\kappa} ^E$ depends on the choice of the coefficients $\rho_{j\kappa}$, but it is not indicated explicitly in the notation. Recall $S_E$ defined by \eqref{constE_MV} (constant energy surface). Note that $S_E$ is symmetric in all variables. So $S'_E := \left. S_E \right|_{\R^{d-1}}$ is unambiguously defined. The point is that $\left. G_{\kappa}^E \right|_{S'_E}$ is a functional representation of $1_{\{ E \} } (D) [D, \i \mathbb{A}]_{\circ} 1_{\{ E \} } (D)$. By functional calculus and continuity of the function $G_{\kappa} ^E$, $E \in \boldsymbol{\mu}_{\kappa} (D)$ if and only if $\left. G_{\kappa}^E \right|_{S'_E} > 0$.

We highlight specially the 2 and 3-dimensional cases as this is our main focus. In dimension 2, we adopt the simpler notation :
\begin{equation}
\label{def:gE22}
g_{j\kappa} ^E (x) = (E/x) \cdot m(x) U_{j\kappa - 1} (x) + x \cdot m(E/x) U_{j\kappa - 1} (E/x), \quad x \in [-1,-|E|] \cup [|E|,1].
\end{equation}
In dimension 3, we adopt the simpler notation :
\begin{equation}
\label{def:gE33}
g_{j\kappa} ^E (x,y) = (E/x) \cdot m(x) U_{j\kappa - 1} (x) + (E/y) \cdot m(y) U_{j\kappa - 1} (y)  + xy \cdot m(E/(xy)) U_{j\kappa - 1} (E/(xy)),
\end{equation}
$y \in [-1,-|E|] \cup [|E|,1]$ and $x \in [-1,-|E/y|] \cup [|E/y|,1]$.

The proof of the following Lemma follows directly from the definition of $g_{j \kappa} ^E(x)$.
\begin{Lemma} 
\label{even_pol}
Fix $d=2$. If $\kappa \in \N_e$, then $G_{\kappa} ^E(x)$ is an even polynomial for any coefficients $\rho_{j\kappa}$.
\end{Lemma}

\begin{Lemma}
\label{symmetry_minima_E/2_mv}
Let $d=2$, $E>0$, $\kappa \in \N^*$. Then $\frac{d}{dx} g_{j\kappa}^E (\sqrt{E}) = 0$ for all $j \in \N^*$. In particular $\frac{d}{dx} G_{\kappa}^E(\sqrt{E}) = 0$ for any choice of coefficients $\rho_{j \kappa}$. 
\end{Lemma}
\begin{proof}
Straightforwardly from \eqref{def:gE22}.
\qed
\end{proof}

In \cite{GM3}, the function $g_{j \kappa} ^E(x)$, and hence $G_{\kappa} ^E(x)$, had a nice visual symmetry property, namely, it was symmetric wrt.\ the axis $x = E/2$, see \cite[Lemma 3.4]{GM3}. In this article, we still have a pseudo-symmetry property, which is that 
\begin{equation}
\label{mult_symmetry_property}
g_{j \kappa} ^E (\sqrt{E} \cdot t) = g_{j \kappa}^E (\sqrt{E}/t), \quad \forall j \in \N^*, \quad \text{which entails} \quad G_{\kappa} ^E (\sqrt{E} \cdot t) = G_{\kappa}^E (\sqrt{E}/t), \quad \forall t \in \R.
\end{equation}
We will refer to this symmetry property as a \textit{multiplicative} \textit{symmetry wrt.\ the axis} $x = \sqrt{E}$. Needless to say that this is a highly questionable wording.

\section{Generalities about the sets  \texorpdfstring{$\boldsymbol{\mu}_{\kappa} (D)$, $\boldsymbol{\Theta}_{\kappa} (D)$ and $\boldsymbol{\Theta}_{m,\kappa} (D)$}{TEXT}}

For the proofs of the Lemmas \ref{Lemma_symmetryDelta}, \ref{lemSUMcosINTRO} and \ref{shift_threshold_d} given just below, we revert back to the notation \eqref{def:gE_intro} instead of \eqref{def:gE}. Thanks to Lemmas \ref{Lemma_symmetryDelta} and \ref{Lemma_symmetryDelta_88} we may focus only on positive energies in this article. 

\begin{Lemma}
\label{Lemma_symmetryDelta}
$\forall d\geq1$, $\forall \kappa \in \N_e$, $\boldsymbol{\mu}_{\kappa} (D) = - \boldsymbol{\mu}_{\kappa} (D)$. Taking complements, $\boldsymbol{\Theta}_{\kappa} (D) = - \boldsymbol{\Theta}_{\kappa} (D)$.
\end{Lemma}
\begin{proof}
First note that $S_{-E} = \bigcup_{j=1}^d \{ \lambda_j(x_1, ..., x_d) : (x_1, ..., x_d) \in S_E\},$ where $\lambda_j (x_1, ..., x_d) = (f_1,...,f_d)$, with $f_i = x_i$ if $i \neq j$ and $f_j = -x_j$. Thus, from \eqref{def:gE_intro}, and the fact that the $U_{j\kappa-1}(\cdot)$ are odd functions, we have that 
$g_{j\kappa} (x_1, ..., x_d)=  - g_{j\kappa} (\lambda_j(x_1, ..., x_d))$ whenever $(x_1, ..., x_d) \in S_E$.
\qed
\end{proof}

\begin{Lemma}
\label{Lemma_symmetryDelta_88}
For any $d \in \N^*$, for any $\kappa \in \N_e$, any $m \in \N$, $\boldsymbol{\Theta}_{m,\kappa} (D) = - \boldsymbol{\Theta}_{m,\kappa} (D)$.
\end{Lemma}


The proof of Lemma \ref{Lemma_symmetryDelta_88} is like that of Lemma \ref{Lemma_symmetryDelta} and follows from the definition of $\boldsymbol{\Theta}_{m,\kappa} (D)$. 

\noindent \textit{Proof of Lemma \ref{lemSUMcosINTRO}}.
Let $E= \prod_{q=1} ^d x_q$, $x_q = \cos( j_q \pi / \kappa)$. Then $g_{j\kappa}(x_1, ..., x_d) = 0$, $\forall j \in \N^*$.
\qed

\noindent \textit{Proof of Lemma \ref{shift_threshold_d}}.
Let $E \in  \boldsymbol{\Theta}_{m,\kappa} (D[d])$, with $E = \prod_{l=1} ^d x_{q,l}$ for $0 \leq q \leq m$. Set $\E = E \times \cos(l\pi / \kappa)$, $0 \leq l \leq \kappa$, but $l \neq \kappa/2$. Then $\forall j \in \N^*$:

$$g_{j \kappa} (\vec{x}_m, \cos(l\pi / \kappa)) = g_{j \kappa} (\vec{x}_m) =  \sum _{q =0} ^{m-1} \omega_q \cdot g_{j \kappa} (\vec{x}_q) = \sum _{q =0} ^{m-1} \omega_q \cdot g_{j \kappa} (\vec{x}_q, \cos(l\pi / \kappa)).$$
This implies $\E \in \boldsymbol{\Theta}_{m,\kappa} (D[d+1])$.
\qed

\section{A geometric construction to find thresholds in dimension 2}
\label{geo_construction}

The idea below -- systems \eqref{conjecture_system_conj_intro} and \eqref{conjecture_system_conj_intro_even} -- is our bread and butter to find thresholds $\in \boldsymbol{\Theta}_{m,\kappa} (\Delta)$ in dimension 2. In this section we discuss properties of  solutions to these systems and explain how to construct them graphically. In sections \ref{section J2} and \ref{section J3a} we apply the idea and prove the existence and uniqueness of solutions under certain additional conditions. 

For subsections \ref{sub_odd} and \ref{sub_even}, fix $\kappa \in \N_e$, $\kappa \geq 4$, and $n \in \N^*$. Consider real variables $\E_n, X_{0,n}, X_{1,n}, ..., X_{n,n}, X_{n+1,n}$. 

\subsection{The case of $n$ odd.}
\label{sub_odd}

For $n \in \N^*$ \underline{odd} define $\mathfrak{T}_{n,\kappa}$ to be the set of \textit{strictly positive} $\E_n = \E_n (\kappa)$ such that the system of $n+3$ equations in $n+3$ unknowns 
\begin{equation}
\label{conjecture_system_conj_intro}
\begin{cases}
T_{\kappa}(X_{q,n}) = T_{\kappa}(X_{n-q,n}), \quad \forall \ q = 0,1,..., (n-1)/2&  \\
X_{n-q,n} = \E_n / X_{1+q,n}, \quad \forall \ q = -1, 0,..., (n-1)/2  & 
\end{cases}
\end{equation}
has a solution satisfying 
\begin{align}
\label{o1}
& X_{q,n} \in (-1,1) \setminus \{0\}, \quad \forall \ 0 \leq q \leq (n-1)/2, \\ 
\label{o2}
& T'_{\kappa}(X_{q,n}) \neq 0, \quad \forall \ 0 \leq q \leq (n-1)/2, \\
\label{o3}
& X_{n+1,n} \in \{ \cos(j \pi / \kappa) : 0 \leq j \leq \kappa \} \setminus \{0\}.
\end{align}
To be clear the $\E_n$ and $X_{i,n}$'s depend on both $\kappa$ and $n$. When $n$ is fixed and no confusion can arise we will simply write $X_i$ instead of $X_{i,n}$. We start with a useful observation, which follows immediately from the variations of $T_{\kappa}$, cf.\ Lemma \ref{variationsTk} :
\begin{Lemma}
\label{lem_sym_odd}
Suppose system \eqref{conjecture_system_conj_intro} and \eqref{o1}--\eqref{o3} hold. Then
\begin{align}
\label{o1bis} 
& X_{n-q,n} \in (-1,1) \setminus \{0\}, \quad \forall \ 0 \leq q \leq (n-1)/2, \\
\label{o2bis}
& T'_{\kappa}(X_{n-q,n}) \neq 0, \quad \forall \ 0 \leq q \leq (n-1)/2.
\end{align}
\end{Lemma}

It is possible however (and will be the case sometimes) that $X_{n+1,n} \in \{\pm 1\}$ or $T'_{\kappa}(X_{n+1,n}) = 0$. 

\begin{remark} System \eqref{conjecture_system_conj_intro} cannot admit a solution if $\E_n <0$, because the second line of \eqref{conjecture_system_conj_intro} implies $X_{(n+1)/2} ^2 = \E_n$. Thus, for $\E_n < 0$, we define $\mathfrak{T}_{n,\kappa}$ to be the set of strictly negative $\E_n = \E_n (\kappa)$ such that the system
\begin{equation}
\label{conjecture_system_conj_intro_22}
\begin{cases}
T_{\kappa}(X_{q,n}) = T_{\kappa}(X_{n-q,n}), \quad \forall \ q = 0,1,..., (n-1)/2&  \\
X_{n-q,n} = \E_n / X_{1+q,n}, \quad \forall \ q = -1, 0,..., (n-3)/2  & \\
X_{(n+1)/2} ^2 = -\E_n & 
\end{cases}
\end{equation}
has a solution satisfying \eqref{o1}--\eqref{o3}. Of course, Lemma \ref{lem_sym_odd} holds for solutions with $\E_n <0$.
\end{remark}
Given Lemma \ref{Lemma_symmetryTT} below we allow ourselves to focus only on positive energies from here on.

\begin{Lemma}
\label{Lemma_symmetryTT}
$\mathfrak{T}_{n, \kappa} = - \mathfrak{T}_{n, \kappa}$, for all $n \in \N_o$, $\kappa \geq 4$, $\kappa \in \N_e$.
\end{Lemma}

\begin{proof}
We are in dimension 2. Note that $X_q \in S_{\E_n} \Leftrightarrow -X_q \in S_{-\E_n}$. Fix $\E_n >0$. $\E_n$ and $X_q$, $0 \leq q \leq n+1$ solve \eqref{conjecture_system_conj_intro} and satisfy \eqref{o1}--\eqref{o3} if and only if $-\E_n$, $-X_q$, $0 \leq q \leq (n-1)/2$ and $X_q$, $(n+1)/2 \leq q \leq n+1$ solve \eqref{conjecture_system_conj_intro_22} and satisfy \eqref{o1}--\eqref{o3}. Note we used the parity of $T_{\kappa}$.
\qed
\end{proof}

\begin{figure}[htb]
  \centering
   \includegraphics[scale=0.197]{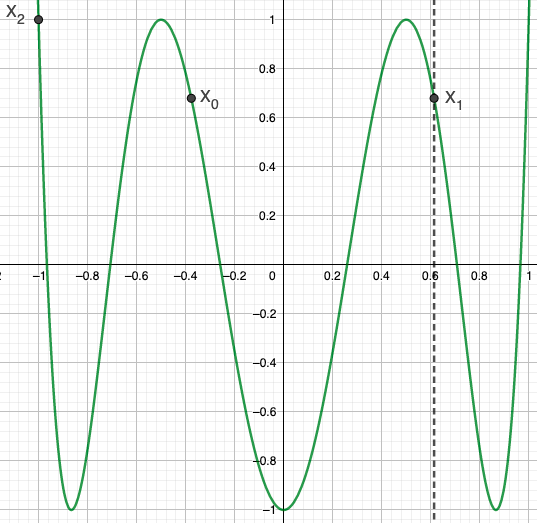}
 \includegraphics[scale=0.174]{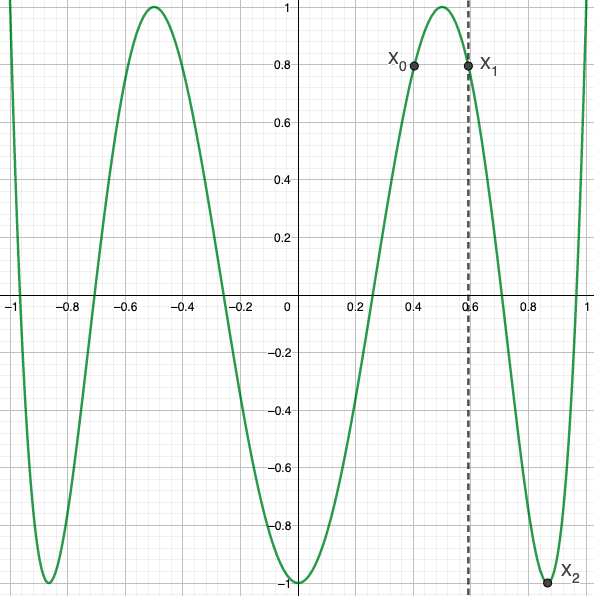}
  \includegraphics[scale=0.18]{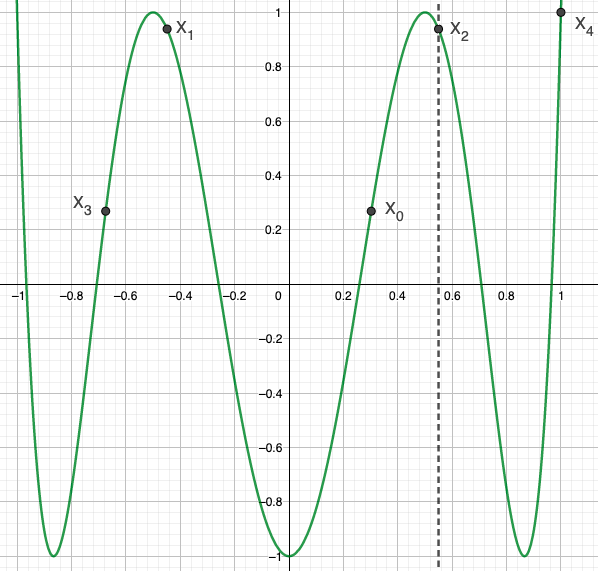}
  \includegraphics[scale=0.18]{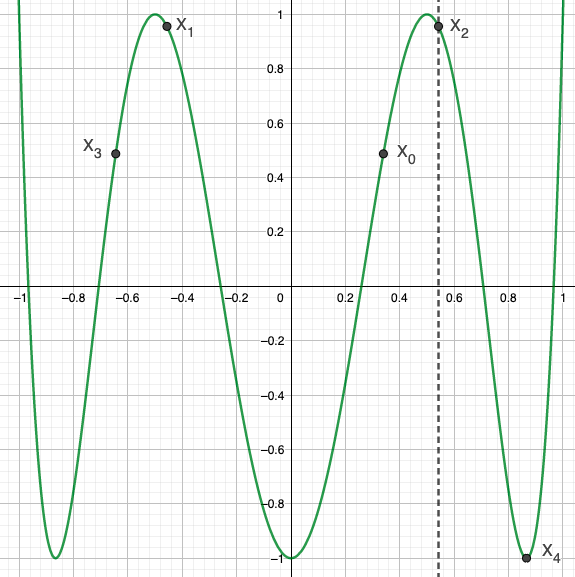}
\caption{$T_{\kappa=6}(x)$. Thresholds $\in \boldsymbol{\Theta}_{m,\kappa} (D)$ in 2nd gap = $(\cos^2(2\pi / 6), \cos(2\pi / 6))$. Left to right : $\E_1 \simeq 0.37677$, $\E_1 \simeq 0.34918$, $\E_3 \simeq 0.30236$, $\E_3 \simeq 0.29420$}
\label{fig:example_2nd_gap}
\end{figure}

We discuss the a priori non-uniqueness of the solutions to the system. There may be several solutions to system \eqref{conjecture_system_conj_intro} satisfying \eqref{o1}--\eqref{o3} for fixed $n$ and $\kappa$. This is because $T_{\kappa}$ is not an injective function on $[-1,1]$. So for given $x$, $T_{\kappa}(x) = T_{\kappa}(y)$ has several solutions and this in turn means there is generally an abundance of solutions. Another reason why there may be several solutions is because there are several options for $X_{n+1}$. So for example, the $1st$ and $2nd$ graphs in Figure \ref{fig:example_2nd_gap} are solutions to the system for $(\kappa,n) = (6,1)$ whereas the $3rd$ and $4th$ graphs in that Figure illustrate solutions to the system for $(\kappa,n) = (6,3)$. These are examples where the energy solutions $\E_n(\kappa)$ are different. In Figure \ref{fig:test_T3k3firsty_versions} we display examples of solutions to the system for given $(\kappa,n)$ for which the energy $\E_n(\kappa)$ is the same but the configuration of the $X_q$'s is different. System \eqref{conjecture_system_conj_intro} together with \eqref{o1}--\eqref{o3} is therefore not enough to guarantee a unique solution.


The following Proposition establishes a link between system \eqref{conjecture_system_conj_intro} and the definition of $\boldsymbol{\Theta}_{m,\kappa} (D)$:

\begin{proposition} (odd terms)
\label{prop1_intro}
Fix $\kappa \in \N_e$, $\kappa \geq 4$, and let $n \in \N_o$ be given. Let $(X_q)_{q=0}^{n+1}$, $\mathcal{E}_n (\kappa)$ be a solution such that $\mathcal{E}_n (\kappa) \in \mathfrak{T}_{n,\kappa}$, with $\E_n >0$. Then for all $j \in \N^*$,
\begin{equation}
\label{problem1_tosolve_again_101_88}
g_{j \kappa} ^{\mathcal{E}_n}  ( X_{\frac{n+1}{2}} ) = \sum _{q =0} ^{\frac{n-1}{2}} \omega_q \cdot g_{j \kappa} ^{\mathcal{E}_n} ( X_{q}), \quad \omega_q = 2 \cdot (-1)^{\frac{n-1}{2}-q} \cdot \frac{\prod _{p=q} ^{\frac{n-1}{2}} X_{p}}{\prod _{p'=\frac{n+1}{2}} ^{n-q} X_{p'}} \cdot \frac{\prod _{p=\frac{n+1}{2}} ^{n-q} m(X_p) U_{\kappa-1} (X_p)}{\prod _{p'=q} ^{\frac{n-1}{2}} m(X_{p'}) U_{\kappa-1} (X_{p'})}.
\end{equation}
\end{proposition}

\begin{remark}
The formula for $\omega_q$ in \eqref{problem1_tosolve_again_101_88} is slightly more complicated than the one for the standard Laplacian, see \cite[Proposition 5.1]{GM3}.
\end{remark}
Note that the $\omega_q$ are independent of $j$ and well-defined thanks to \eqref{o1}, \eqref{o2} and Lemma \ref{lem_sym_odd} (note the use of the first identity in \eqref{identity_derivative}). Recall that for our energies to belong to $\boldsymbol{\Theta}_{m,\kappa} (D)$ we want $\omega_q \leq 0$. To this end we introduce 3 additional assumptions to be considered separately. \\

\noindent \underline{\blue{Additional Assumption for $n$ odd :}}
\begin{enumerate}[label=\textbf{AO.\arabic*}]
\item \label{ao1} $X_q \cdot X_{n-q}> 0$ and $T'_{\kappa}(X_{q}) \cdot T'_{\kappa}(X_{n-q}) < 0$ for $0 \leq q \leq (n-1)/2$.
\item \label{ao2} $T'_{\kappa}(X_{q}) \cdot T'_{\kappa}(X_{n-q})  > 0$ and $X_q \cdot X_{n-q} < 0$ for $0 \leq q \leq (n-1)/2$.
\item \label{ao3} ("\textit{mix and match"}) there are disjoint $I_1,I_2 \subset \{ q : 0 \leq q \leq (n-1)/2 \}$ such that $I_1 \sqcup I_2  = \{ q : 0 \leq q \leq (n-1)/2 \}$, and $X_q \cdot X_{n-q}> 0$, $T'_{\kappa}(X_{q}) \cdot T'_{\kappa}(X_{n-q}) < 0$ for $q \in I_1$, and $T'_{\kappa}(X_{q}) \cdot T'_{\kappa}(X_{n-q})  > 0$, $X_q \cdot X_{n-q} < 0$ for $q \in I_2$.
\end{enumerate}

Note that assumption \ref{ao1} (resp.\ \ref{ao2}) is a special case of \ref{ao3} where $I_2 = \emptyset$ (resp.\ $I_1 = \emptyset$). The following Corollary completes the link between system \eqref{conjecture_system_conj_intro} and $\boldsymbol{\Theta}_{m,\kappa} (D)$.

\begin{corollary}
\label{imp_cor_odd}
Let $(X_q)_{q=0}^{n+1}$, $\mathcal{E}_n (\kappa)$ be exactly as in Proposition \ref{prop1_intro} so that \eqref{problem1_tosolve_again_101_88} holds. If \ref{ao3} further holds, then the $\omega_q$ are all strictly negative and so $\E_n(\kappa) \in \boldsymbol{\Theta}_{(n+1)/2,\kappa} (D[d=2])$.
\end{corollary}

\begin{remark}
\label{RR1}
Assumptions \ref{ao1} and \ref{ao3} are likely not suitable if $\E_n <0$, because then $X_{1+q}$ and $X_{n-q}$ need to have opposite signs as per the second line of \eqref{conjecture_system_conj_intro}. If $\E_n >0$ however, then we have examples to illustrate \ref{ao1}, \ref{ao2} and \ref{ao3}, see for instance Figures \ref{fig:test_T3k3firsty} and \ref{fig:test_T3k3firsty_versions}.
\end{remark}

Certain graphs in Figures \ref{fig:test_T3k3firsty} and \ref{fig:test_T3k3firsty_versions} illustrate solutions to system \eqref{conjecture_system_conj_intro} for $\kappa =4$, and satisfying \eqref{o1}--\eqref{o3}. The $1st$, $3rd$ and $5th$ graphs in Figure \ref{fig:test_T3k3firsty} satisfy assumption \ref{ao1} with $X_n$, $X_{n-q} > 0$. The $3rd$ graph in Figure \ref{fig:test_T3k3firsty_versions} satisfies assumption \ref{ao1} with $X_n$, $X_{n-q} < 0$. The $1st$ and $6th$ graphs in Figure \ref{fig:test_T3k3firsty_versions} satisfy assumption \ref{ao2}, whereas the $5th$ graph in Figure \ref{fig:test_T3k3firsty_versions} satisfies assumption \ref{ao3}.

\begin{figure}[htb]
  \centering
   \includegraphics[scale=0.16]{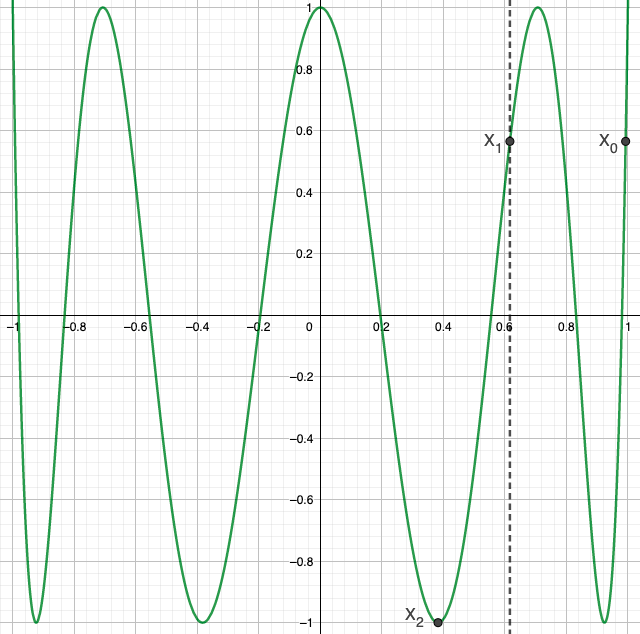}
 \includegraphics[scale=0.16]{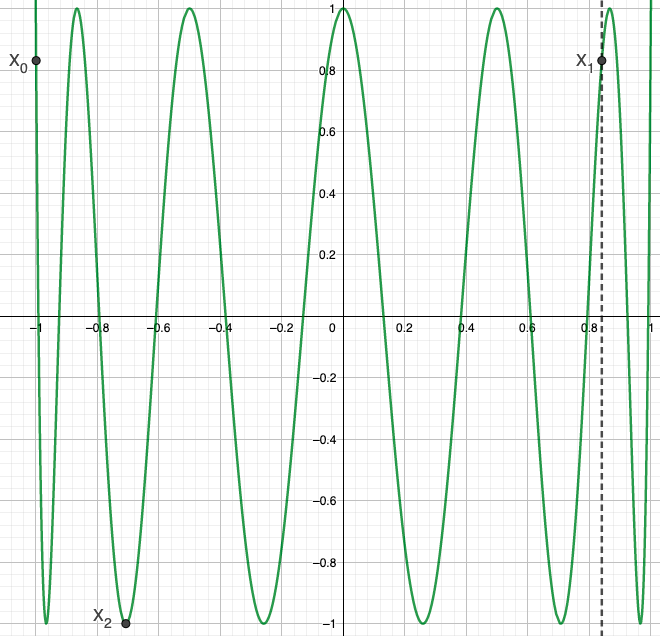}
\caption{Graph of $T_{\kappa}(x)$. Solutions to system \eqref{conjecture_system_conj_intro} satisfying \eqref{o1}--\eqref{o3}. Left : $(\kappa,n)=(8,1)$, $\E_1 \simeq 0.37987$. Right : $(\kappa,n)=(12,1)$, $\E_1 \simeq 0.70625$. However unclear if these energies $\E_1$ belong to $\boldsymbol{\Theta}_{\kappa}(D[d=2])$ or $\boldsymbol{\mu}_{\kappa}(D[d=2])$}
\label{fig:example_unclear}
\end{figure}

\begin{remark}
\label{RR2}
It is unclear to me if there are other assumptions that can ensure $\omega_q \leq 0$ (but see the formulas and discussion in section \ref{revisit})
\end{remark}
\begin{remark}
\label{RR3}
It is easy to graphically build solutions to \eqref{conjecture_system_conj_intro} that satisfy \eqref{o1}--\eqref{o3}, but don't satisfy \ref{ao3}. Such examples are given in Figure \ref{fig:example_unclear}. Note that in these examples \eqref{problem1_tosolve_again_101_88} holds perfectly well, although
$$\omega_0 = 2 \frac{X_0}{X_1} \cdot \frac{m(X_1)}{m(X_0)} \cdot \frac{U_{\kappa-1}(X_1)}{U_{\kappa-1}(X_0)} > 0$$
and so by definition those solutions $\E_1$ do not belong to $\boldsymbol{\Theta}_{1,\kappa}(D[d=2])$. After crossing this information with bands identified in \cite[Table XV]{GM2} it remains unclear to me if these energies belong to $\boldsymbol{\Theta}_{\kappa}(D[d=2])$ or $\boldsymbol{\mu}_{\kappa}(D[d=2])$. I am not quite happy with these examples because I was hoping to find an example of a solution to \eqref{conjecture_system_conj_intro} that satisfies \eqref{o1}--\eqref{o3}, but not \ref{ao3}, and then further be able to confidently confirm that the solution belongs to $\boldsymbol{\mu}_{\kappa}(D[d=2])$. In other words, I haven't been able to disprove the possibility of a threshold energy for which a linear combination of the form \eqref{problem1_tosolve_again_101_88} holds with an $\omega_q >0$. Maybe there is something more to understand here.
\end{remark}


\subsection{The case of $n$ even.}
\label{sub_even}

We move on with the case of $n$ even. For $n \in \N^*$ \underline{even} define $\mathfrak{T}_{n,\kappa}$ to be the set of non-zero $\E_n$ such that the system of $n+3$ equations in $n+3$ unknowns 

\begin{equation}
\label{conjecture_system_conj_intro_even}
\begin{cases}
T_{\kappa}(X_{q,n}) = T_{\kappa}(X_{n-q,n}), \quad \forall \ q = 0,1,..., n/2-1 &  \\
X_{n-q,n} = \E_n / X_{1+q,n}, \quad \forall \ q = -1, 0,..., n/2 -1  & 
\end{cases}
\end{equation}
has a solution satisfying 
\begin{align}
\label{o11}
& X_{q,n} \in (-1,1) \setminus \{0\}, \quad \forall \ 0 \leq q \leq n/2-1, \\
\label{o22}
& T'_{\kappa}(X_{q,n}) \neq 0, \quad \forall \ 0 \leq q \leq n/2-1, \\
\label{o33}
& X_{n+1,n} \in \{ \cos(j \pi / \kappa) : 0 \leq j \leq \kappa \} \setminus \{0\}, \\
\label{o44}
& X_{n/2,n} \in \{ \cos(j \pi / \kappa) : 0 \leq j \leq \kappa \} \setminus \{0\}. 
\end{align}
Again, the $\E_n$ and $X_{i,n}$'s depend on both $\kappa$ and $n$. When $n$ is fixed and no confusion can arise we will simply write $X_i$ instead of $X_{i,n}$. The following observation follows immediately from the variations of $T_{\kappa}$, cf.\ Lemma \ref{variationsTk} :
\begin{Lemma}
\label{lem_sym_even}
Suppose system \eqref{conjecture_system_conj_intro_even} and \eqref{o11}--\eqref{o44} hold. Then
\begin{align}
\label{o11bis} 
& X_{n-q,n} \in (-1,1) \setminus \{0\}, \quad \forall \ 0 \leq q \leq n/2-1, \\
\label{o22bis}
& T'_{\kappa}(X_{n-q,n}) \neq 0, \quad \forall \ 0 \leq q \leq n/2-1.
\end{align}
\end{Lemma}

It is possible however (and will be the case sometimes) that $X_{n/2,n} \in \{ \pm 1\}$, or $X_{n+1,n} \in \{ \pm 1\}$, or $T'_{\kappa}(X_{n/2,n}) = 0$, or $T'_{\kappa}(X_{n+1,n}) = 0$. Unlike in the $n$ odd case, here we don't have to make a distinction between positive or negative energy solutions. Furthermore, we also have :
\begin{Lemma}
\label{Lemma_symmetryTT_neven}
$\mathfrak{T}_{n, \kappa} = - \mathfrak{T}_{n, \kappa}$, for all $n \in \N_e$, $\kappa \geq 4$, $\kappa \in \N_e$.
\end{Lemma}
The following Proposition establishes a link between system \eqref{conjecture_system_conj_intro_even} and $\boldsymbol{\Theta}_{m,\kappa} (D)$:

\begin{proposition} (even terms)
\label{prop1_intro_even}
Fix $\kappa \in \N_e$, $ \kappa \geq 4$, and let $n \in \N_e$ be given. Let $(X_q)_{q=0}^{n+1}$, $\mathcal{E}_n (\kappa)$ be a solution such that $\mathcal{E}_n (\kappa) \in \mathfrak{T}_{n,\kappa}$. Then for all $j \in \N^*$,
\begin{equation}
\label{problem1_tosolve_again_101even_88}
g_{j \kappa} ^{\mathcal{E}_n}  ( X_{\frac{n}{2}} ) = \sum _{q =0} ^{\frac{n}{2}-1} \omega_q \cdot g_{j \kappa} ^{\mathcal{E}_n}  ( X_{q}), \quad \omega_q =  (-1)^{\frac{n}{2}-1-q} \cdot  \frac{\prod _{p=q} ^{\frac{n}{2}-1}  X_{p}}{\prod _{p'=\frac{n}{2}+1} ^{n-q} X_{p'}}  \cdot \frac{\prod _{p=\frac{n}{2}+1} ^{n-q} m(X_p) U_{\kappa-1} (X_p)}{\prod _{p'=q} ^{\frac{n}{2}-1} m(X_{p'}) U_{\kappa-1} (X_{p'})}.
\end{equation}
\end{proposition}

Note that the $\omega_q$ are independent of $j$ and well-defined thanks to \eqref{o11}, \eqref{o22}, and Lemma \ref{lem_sym_even}. As in the case of $n$ odd we consider 3 additional assumptions to be considered separately. \\

\noindent \underline{\blue{Additional Assumption for $n$ even :}}
\begin{enumerate}[label=\textbf{AE.\arabic*}]
\item \label{ae1} $X_q \cdot X_{n-q} > 0$ and $T'_{\kappa}(X_{q}) \cdot T'_{\kappa}(X_{n-q}) < 0$ for $0 \leq q \leq n/2-1$.
\item \label{ae2} $T'_{\kappa}(X_{q}) \cdot T'_{\kappa}(X_{n-q}) > 0$ and $X_q \cdot X_{n-q} < 0$ for $0 \leq q \leq n/2-1$.
\item \label{ae3} ("\textit{mix and match"}) there are disjoint $I_1,I_2 \subset \{ q : 0 \leq q \leq n/2-1 \}$ such that $I_1 \sqcup I_2  = \{ q : 0 \leq q \leq n/2-1 \}$, and $X_q \cdot X_{n-q}> 0$, $T'_{\kappa}(X_{q}) \cdot T'_{\kappa}(X_{n-q}) < 0$ for $q \in I_1$, and $T'_{\kappa}(X_{q}) \cdot T'_{\kappa}(X_{n-q})  > 0$, $X_q \cdot X_{n-q} < 0$ for $q \in I_2$.
\end{enumerate}

Note that assumption \ref{ae1} (resp.\ \ref{ae2}) is a special case of \ref{ae3} where $I_2 = \emptyset$ (resp.\ $I_1 = \emptyset$). The following Corollary completes the link between system \eqref{conjecture_system_conj_intro} and the definition of $\boldsymbol{\Theta}_{m,\kappa} (D)$.

\begin{corollary}
\label{imp_cor_even}
Let $(X_q)_{q=0}^{n+1}$, $\mathcal{E}_n (\kappa)$ be exactly as in Proposition \ref{prop1_intro_even} so that \eqref{problem1_tosolve_again_101even_88} holds. If \ref{ae3} further holds, then the $\omega_q$ are all strictly negative in which case $\E_n(\kappa) \in \boldsymbol{\Theta}_{n/2,\kappa} (D[d=2])$.
\end{corollary}

Remarks \ref{RR1}--\ref{RR3} apply to the $n$ even case as well. 


\subsection{Expressing the energy $\E_n$ as the solution to a single equation with one unknown}

It is possible to express the thresholds $\E_n$ in systems \eqref{conjecture_system_conj_intro} and \eqref{conjecture_system_conj_intro_even} as solutions to a single equation with a finite number of continued fractions. It's a matter of knowing which branch of $T_{\kappa} ^{-1} (x)$ to choose from. It is easier to explain with an example : 

\begin{example}
Fix $n=3$ and let $X_{n+1} = X_4$ be given as per \eqref{o2}. Then according to system \eqref{conjecture_system_conj_intro}, $\E_3 (\kappa)$ is the solution to the following equation :
\begin{align*}
\E_3  &= X_4 \cdot X_0 = X_4 \cdot T_{\kappa} ^{-1} T_{\kappa} (X_3) =  X_4 \cdot T_{\kappa} ^{-1} T_{\kappa} (\E_3 / X_1)   \\
&= X_4 \cdot T_{\kappa} ^{-1} T_{\kappa} \left(\E_3 / T_{\kappa} ^{-1} T_{\kappa}  (X_2) \right) =  X_4 \cdot T_{\kappa} ^{-1} T_{\kappa} \left(\E_3 / T_{\kappa} ^{-1} T_{\kappa}  \left( \sqrt{\E_3}  \right) \right).
\end{align*}
\end{example}
\normalsize
The advantage of having 1 equation with 1 unknown over a system of equations with several unknowns is that it is easier (in our opinion) to solve numerically. Proposition \ref{lem1_sequence3344} is one (of many) applications of this idea.

\subsection{The systems as dynamical graphical constructions}
\label{subby}

Behind systems \eqref{conjecture_system_conj_intro} and \eqref{conjecture_system_conj_intro_even} is a simple graphical construction and interpretation which is the topic of this subsection. We always assume $\E_n >0$. We start with 2 key remarks/observations.

\begin{remark} 
\label{axis_sym}
As per system \eqref{conjecture_system_conj_intro} (the $n$ odd case), if $\E_n >0$ then $X_{(n+1)/2} = \sqrt{\E_n}$ always holds and this helps for a graphical construction.
\end{remark}

\begin{remark}
\label{obs_symmetry}
Consider equation $X_{n-q} = \E_n / X_{1+q}$ in \eqref{conjecture_system_conj_intro} or \eqref{conjecture_system_conj_intro_even}. For simplicity assume $\E_n, X_{1+q}, X_{n-q} >0$. Then $\min(X_{n-q}, X_{1+q}) \leq \sqrt{\E_n} \leq \max(X_{n-q}, X_{1+q})$. This is very helpful to bear in mind to construct and interpret solutions. 
\end{remark} 

Figures \ref{fig:example_2nd_gap}, \ref{fig:test_T3k3firsty} and \ref{fig:test_T3k3firsty_versions}, for example, illustrate threshold solutions $\in \boldsymbol{\Theta}_{m,\kappa} (D)$. The vertical dotted line is the axis of multiplicative symmetry $x = \sqrt{\E_n}$. The key observation when looking at these graphs is that every point $(X_q, T_{\kappa}(X_q))$ always satisfies 2 crucial conditions :
\begin{enumerate}[label*=\arabic*.]
\item \textit{a multiplicative symmetry condition} : each $X_q$ is the multiplicative symmetric of another point $X_r$ wrt.\ the axis $x= \sqrt{\E_n}$, namely $X_r = X_{n-q+1}$. This is Remark \ref{obs_symmetry}.
\item \textit{a level condition} : each $X_q$ satisfies at least one of the following 3 conditions :
\begin{enumerate}[label*=\arabic*.]
\item $m(X_q) =0$, or
\item $U_{\kappa-1}(X_q) = 0$ (equivalently $T'_{\kappa}(X_q)=0$), or 
\item $\exists X_r$ such that $T_{\kappa}(X_q) = T_{\kappa}(X_r)$, namely $X_r = X_{n-q}$, and furthermore the $X_q$ and $X_r$ satisfy either \ref{ao1} or \ref{ao2} (respectively \ref{ae1} or \ref{ae2}). 
\end{enumerate}
\end{enumerate}

The symmetry and level conditions set the rules of the game to construct valid threshold solutions. Possible constructions are as follows :

\noindent \underline{\textbf{Algorithm for $n$ odd}} -- system \eqref{conjecture_system_conj_intro} along with conditions \eqref{o1}--\eqref{o3} : 
\begin{enumerate}
\item Fix $\kappa \geq 4$, $\kappa \in \N_e$, and plot the Chebyshev polynomial $T_{\kappa}(x)$.
\item Initialize energy $E$ to a certain value such that $\sqrt{E} \in [-1,1] \setminus \{ \cos( j \pi / \kappa) : 0 \leq j \leq \kappa \}$, and draw the vertical axis of multiplicative symmetry $x=\sqrt{E}$.
\item Place the first point $(X_{q_0}, T_{\kappa}(X_{q_0}))$ such that $X_{q_0} = \sqrt{E}$.
\item Place $(X_{q_1}, T_{\kappa}(X_{q_1}))$, $(X_{q_2}, T_{\kappa}(X_{q_2}))$, ..., $(X_{q_{n+1}}, T_{\kappa}(X_{q_{n+1}}))$ by alternating between applying the level condition 2.3 wrt.\ the last point constructed, and then the multiplicative symmetry condition wrt.\ to the last point constructed.
\item Finally calibrate $E$ in such a way that the $x$-coordinate of the last point constructed, $X_{q_{n+1}}$, also satisfies a level condition 2.1 or 2.2. Upon calibration, $E = \E_n$.
\end{enumerate}

\noindent \underline{\textbf{Algorithm for $n$ even}} -- system \eqref{conjecture_system_conj_intro_even} along with conditions \eqref{o11}--\eqref{o44} : 
\begin{enumerate}
\item Fix $\kappa \geq 4$, $\kappa \in \N_e$, and plot the Chebyshev polynomial $T_{\kappa}(x)$.
\item Initialize energy $E$ to a certain value and draw the vertical axis of multiplicative symmetry $x=\sqrt{E}$.
\item Place a first point $(X_{q_0}, T_{\kappa}(X_{q_0}))$ such that $X_{q_0}$ satisfies a level condition 2.1 or 2.2.
\item Place $X_{q_1}$ as the multiplicative symmetric of $X_{q_0}$.
\item Place $(X_{q_2}, T_{\kappa}(X_{q_2}))$, $(X_{q_3}, T_{\kappa}(X_{q_3}))$, ..., $(X_{q_{n+1}}, T_{\kappa}(X_{q_{n+1}}))$ by alternating between applying the level condition 2.3 wrt.\ the last point constructed, and then the multiplicative symmetry condition wrt.\ to the last point constructed.
\item Finally calibrate $E$ in such a way that the $x$-coordinate of the last point constructed, $X_{q_{n+1}}$, also satisfies a level condition 2.1 or 2.2. Upon calibration, $E = \E_n$.
\end{enumerate}

Note that the proposed Algorithms are dynamical constructions : the positioning of all the $X_{q_i}$'s depends on the value of $E$ (with the exception of $X_{q_0}$ in the $n$ even case). In the final step when $E$ is adjusted, all the $X_{q_i}$'s migrate (with the exception of $X_{q_0}$ that stays put in the $n$ even case). Adjusting $E$ preserves the symmetry and level conditions.

\subsection{Other key formulas -- to be used to setup the linear interpolation}
\label{other_key_formulas}

The following Lemma, whose proof immediately follows from the definitions, will be impactful when we do polynomial interpolation.
\begin{Lemma}
\label{obvious_lem}
Fix $d=2$. Let $\E_n, X_{0}, ..., X_{n+1}$ be any solution to system \eqref{conjecture_system_conj_intro}, respectively system \eqref{conjecture_system_conj_intro_even}. Then for $q= -1,0,...,(n-1)/2$, respectively $q= -1,0,...,n/2-1$, 
\begin{equation}
\label{key_obvious}
g_{j \kappa} ^{\E_n} (X_{1+q}) = g_{j \kappa} ^{\E_n} (X_{n-q}), \quad \text{for all} \  j \in \N^*.
\end{equation}
In particular, $G_{\kappa} ^{\E_n} (X_{1+q}) = G_{\kappa} ^{\E_n} (X_{n-q})$ for any choice of coefficients $\rho_{j \kappa}$.
\end{Lemma}

The following Lemma is less obvious but equally impactful for the polynomial interpolation.
\begin{Lemma}
\label{Lemma_not_trivial} Let $\E_n, X_{0}, ..., X_{n+1}$ be any solution to system \eqref{conjecture_system_conj_intro}, respectively system \eqref{conjecture_system_conj_intro_even}. Then for $q= -1,0,...,(n-1)/2$, respectively $q= -1,0,...,n/2-1$, 
\begin{equation}
\label{key_negative}
X_{1+q} \cdot \frac{d}{dx} g_{j \kappa} ^{\E_n} (X_{1+q}) = - X_{n-q} \cdot \frac{d}{dx} g_{j \kappa} ^{\E_n} (X_{n-q}), \quad \text{for all} \  j \in \N^*.
\end{equation}
In particular this implies $X_{1+q} \cdot \frac{d}{dx} G_{\kappa} ^{\E_n} (X_{1+q}) = - X_{n-q} \cdot \frac{d}{dx} G_{\kappa} ^{\E_n} (X_{n-q})$ for any choice of coefficients $\rho_{j \kappa}$, and this for all $q$'s in the aforementioned range.
\end{Lemma}
\begin{remark}
In \cite{GM3} calculations were simpler and instead of \eqref{key_negative} we had more simply $\frac{d}{dx} g_{j \kappa} ^{\E_n} (X_{1+q}) = -\frac{d}{dx} g_{j \kappa} ^{\E_n} (X_{n-q})$ for all $j \in \N^*$.
\end{remark}
\begin{remark}
Assumptions \eqref{o1}--\eqref{o3}, respectively \eqref{o11}--\eqref{o44}, are not needed to prove Lemmas \ref{obvious_lem} and \ref{Lemma_not_trivial}.
\end{remark}
\begin{proof}
First things first, 
\begin{equation}
\begin{aligned}
\label{FFF}
\frac{d}{dx} g_{j \kappa} ^E (x) & = -\frac{E}{x^2} \cdot m(x) U_{j \kappa-1}(x)  - 2E \cdot U_{j \kappa-1}(x) + \frac{E}{x} \cdot m(x) U'_{j \kappa-1}(x) \\
& \quad + m(E/x)  U_{j \kappa-1}(E/x) + 2\frac{E^2}{x^2} \cdot U_{j \kappa-1}(E/x) - \frac{E}{x} \cdot m(E/x) U'_{j \kappa-1}(E/x).
\end{aligned}
\end{equation}
We focus on what comes after the $=$ sign. Group the first two terms on the first row of \eqref{FFF} together after having expanded $m(x)$ into $1-x^2$, and similarly for the first two terms on the second row of \eqref{FFF} ; for the terms with $U'$ apply the second identity in \eqref{identity_derivative}. Thus :

\begin{equation}
\begin{aligned}
\label{FFF2}
\frac{d}{dx} g_{j \kappa} ^E (x) & = -\frac{E}{x^2} U_{j \kappa-1}(x)  - E \cdot U_{j \kappa-1}(x) + E \cdot U_{j \kappa-1}(x) - j\kappa \frac{E}{x} \cdot T_{j \kappa}(x)  \\
& \quad + U_{j \kappa-1}(E/x) + \frac{E^2}{x^2} \cdot U_{j \kappa-1}(E/x) - \frac{E^2}{x^2} \cdot U_{j \kappa}(E/x) + j\kappa \frac{E}{x} \cdot T_{j \kappa}(E/x) \\
& = -\frac{E}{x^2} U_{j \kappa-1}(x) - j\kappa \frac{E}{x} \cdot T_{j \kappa}(x) + U_{j \kappa-1}(E/x) + j\kappa \frac{E}{x} \cdot T_{j \kappa}(E/x).
\end{aligned}
\end{equation}
Evaluate $\frac{d}{dx} g_{j \kappa} ^E (x)$ at $x = X_{1+q}$, $X_{n-q}$ using the third row of \eqref{FFF2}. Recalling $\E_n = X_{1+q} \cdot X_{n-q}$:
$$ X_{1+q} \cdot \frac{d}{dx}  g_{j \kappa} ^{\E_n} (X_{1+q}) = - X_{n-q} \cdot U_{j \kappa-1}(X_{1+q})  + X_{1+q} \cdot U_{j \kappa-1}(X_{n-q}) + j\kappa \E_n  \left(T_{j \kappa}(X_{n-q}) - T_{j \kappa}(X_{1+q})\right),$$
$$ X_{n-q} \cdot \frac{d}{dx} g_{j \kappa} ^{\E_n} (X_{n-q}) = - X_{1+q} \cdot U_{j \kappa-1}(X_{n-q})  + X_{n-q} \cdot U_{j \kappa-1}(X_{1+q}) + j\kappa \E_n  \left(T_{j \kappa}(X_{1+q}) - T_{j \kappa}(X_{n-q})\right).$$
The last two lines are the negative of each other.
\qed
\end{proof}

\subsection{Proofs of Propositions \ref{prop1_intro}, \ref{prop1_intro_even} and Corollaries \ref{imp_cor_odd} and \ref{imp_cor_even}} \hfill\\

\noindent \textit{Proof of Proposition \ref{prop1_intro}.} Recall $\E_n >0$ is assumed. From \eqref{conjecture_system_conj_intro}, $X_{(n+1)/2} = \sqrt{\mathcal{E}_n }  \Rightarrow g_{j \kappa} ^{\mathcal{E}_n} ( X_{(n+1)/2} ) = 2 X_{(n+1)/2} \cdot m(X_{(n+1)/2}) U_{j \kappa-1}( X_{(n+1)/2})$. For $q \geq 1$, let 
\begin{equation*}
\begin{aligned}
\tilde{\omega}_q &:= \frac{\omega_q}{2} \cdot \prod _{r=(n+1)/2} ^{n} X_r \cdot \prod _{r=0} ^{(n-1)/2} m(X_{r}) U_{\kappa-1} (X_{r}) \\
&=  (-1)^{\frac{n-1}{2}-q} \cdot \prod _{r=n-q+1} ^{n} X_r \cdot \prod _{r=0} ^{q-1} m(X_{r}) U_{\kappa-1} (X_{r}) \cdot \prod _{p=q} ^{(n-1)/2} X_p \cdot \prod _{p=(n+1)/2} ^{n-q} m(X_p) U_{\kappa-1} (X_p) .
\end{aligned}
\end{equation*} 
In order for what follows to apply to the cases $n=1,3$ as well, interpret $\prod _{\alpha} ^{\beta} =0$ and $\sum_{\alpha} ^{\beta} =0$ whenever $\beta < \alpha$. Multiplying \eqref{problem1_tosolve_again_101_88} throughout by $2^{-1} \prod _{r=(n+1)/2} ^{n} X_r \cdot   \prod _{r=0} ^{(n-1)/2} m(X_{r}) U_{\kappa-1} (X_{r}) $ shows that \eqref{problem1_tosolve_again_101_88} is equivalent to
\begin{equation}
\label{problem1_tosolve_again_202}
\begin{aligned}
&X_{(n+1)/2} \cdot m(X_{(n+1)/2}) U_{j \kappa-1}( X_{(n+1)/2}) \cdot \prod _{r=(n+1)/2} ^{n} X_r \cdot \prod _{r=0} ^{(n-1)/2} m(X_{r}) U_{\kappa-1} (X_{r}) \\
& \quad = (-1)^{(n-1)/2} \cdot X_{n+1} \cdot m(X_0) U_{j\kappa-1} (X_0) \cdot \prod _{p=0} ^{(n-1)/2} X_p \cdot \prod _{p=(n+1)/2} ^{n} m(X_p) U_{\kappa-1} (X_p)  \\
& \quad \quad +  \sum _{q =1} ^{(n-1)/2} \tilde{\omega}_q \cdot X_{n-q+1} \cdot m(X_q) U_{j\kappa-1} (X_q) + \sum _{q =1} ^{(n-1)/2} \tilde{\omega}_q \cdot X_q \cdot  m(X_{n-q+1}) U_{j\kappa-1} (X_{n-q+1}) \\
& \quad = (-1)^{(n-1)/2} \cdot X_{n+1} \cdot m(X_0) U_{j\kappa-1} (X_0) \cdot \prod _{p=0} ^{(n-1)/2} X_p \cdot \prod _{p=(n+1)/2} ^{n} m(X_p) U_{\kappa-1} (X_p)  \\
& \quad \quad + X_{(n+3)/2} \cdot X_{(n-1)/2} \cdot m(X_{(n-1)/2}) U_{j\kappa-1} (X_{(n-1)/2}) \cdot m(X_{(n+1)/2}) U_{\kappa-1} (X_{(n+1)/2}) \times \\
& \quad \quad \quad \quad \times \prod _{r=(n+3)/2} ^{n} X_r \cdot \prod _{p'=0} ^{(n-3)/2} m(X_{p'}) U_{\kappa-1} (X_{p'}) \\
& \quad \quad - (-1)^{(n-1)/2} \cdot X_1 \cdot X_n \cdot m(X_0) U_{\kappa-1} (X_0) \cdot m(X_n) U_{j\kappa-1} (X_n) \cdot \prod_{p=(n+1)/2} ^{n-1} m(X_p) U_{\kappa-1} (X_p) \times \\
& \quad \quad \quad \quad \times \prod _{p=1} ^{(n-1)/2} X_p \cdot \prod _{p=(n+1)/2} ^{n-1} m(X_p) U_{\kappa-1} (X_p)  \\
& \quad \quad +  \sum _{q =1} ^{(n-3)/2} \tilde{\omega}_q \cdot X_{n-q+1} \cdot m(X_q) U_{j\kappa-1} (X_q) + \sum _{q =2} ^{(n-1)/2} \tilde{\omega}_q \cdot X_{q} \cdot m(X_{n-q+1}) U_{j\kappa-1} (X_{n-q+1}).
\end{aligned}
\end{equation}
We focus on what comes after the last = sign. We apply Corollary \ref{corollaryEquiv}. The 2nd term on the rhs of \eqref{problem1_tosolve_again_202} equals the lone term on the lhs of \eqref{problem1_tosolve_again_202}. The 3rd term on the rhs of \eqref{problem1_tosolve_again_202} cancels the 1st term on the rhs of \eqref{problem1_tosolve_again_202}. Finally the 2 sums at the very end of the rhs of \eqref{problem1_tosolve_again_202} cancel each other ; specifically, the $q^{th}$ term in the first sum equals
\footnotesize
$$(-1)^{(n-1)/2-q} \cdot X_{n-q+1} \cdot m(X_q) U_{j\kappa-1} (X_q) \cdot  \prod _{p=(n+1)/2} ^{n-q} m(X_p) U_{\kappa-1} (X_p) \cdot \prod _{p'=0} ^{q-1} m(X_{p'}) U_{\kappa-1} (X_{p'}) \cdot \prod _{r=n-q+1} ^{n} X_r \cdot \prod _{p=q} ^{(n-1)/2} X_p$$
\normalsize
and it cancels the $q+1^{th}$ term in the second sum which equals
\footnotesize
$$-(-1)^{(n-1)/2-q} \cdot X_{q+1} \cdot m(X_{n-q}) U_{j\kappa-1} (X_{n-q}) \cdot \prod _{p=(n+1)/2} ^{n-q-1} m(X_p) U_{\kappa-1} (X_p) \cdot \prod _{p'=0} ^{q} m(X_{p'}) U_{\kappa-1} (X_{p'})  \cdot \prod _{r=n-q} ^{n} X_r \cdot \prod _{p=q+1} ^{(n-1)/2} X_p ,$$
\normalsize
and this for $q=1,2,...,(n-3)/2$.
\qed

\noindent \textit{Proof of Corollary \ref{imp_cor_odd}.} 
We want to prove that $\omega_q <0$ under assumption \ref{ao3}. Assume $I_1, I_2$ are a partition of $\{0,1,...,(n-1)/2\}$ as in assumption \ref{ao3}. Fix $q$, $0 \leq q \leq (n-1)/2$. If $\lambda(q) := | \{ q' \in I_1 : q' \geq q\} |$ (the cardinality of the set), then $| \{ q' \in I_2 : q' \geq q\} | = (n-1)/2-q+1-\lambda(q)$. 
Because $X_q, X_{n-q} \in (-1,1)$, $\forall q=0,..., (n-1)/2$, $m(X_q), m(X_{n-q}) \in (0,1)$, and so the sign of $\omega_q$ is that of $(-1)^{(n-1)/2-q} \times (-1)^{(n-1)/2-q + 1 -\lambda(q)} \times (-1)^{\lambda(q)} = (-1)^n = -1$ since $n$ is odd.
\qed

\vspace{0.5cm}
\noindent \textit{Proofs of Proposition \ref{prop1_intro_even} and Corollary \ref{imp_cor_even}.} Very similar to that of Proposition \ref{prop1_intro} and Corollary \ref{imp_cor_odd} so we leave it to the reader. 
\qed

\section{A decreasing sequence of thresholds in \texorpdfstring{$J_2(\kappa) :=  (\cos ^2 (\pi / \kappa ) , \cos(\pi / \kappa) )$}{TEXT}}
\label{section J2}

This entire section is in dimension 2. Using ideas of section \ref{geo_construction} we prove the existence of a sequence of threshold energies $\in J_2 :=  (\cos ^2 (\pi / \kappa ) , \cos(\pi / \kappa) )$. Theorem \ref{thm_decreasing energy general} is a consequence of Propositions \ref{prop1}, \ref{prop2}, \ref{prop_interlace}, \ref{prop_convergeence}. We state these Propositions now, and prove them at the end of the section. A justification for Conjecture \ref{conjecture12} is also given.

\begin{proposition} (odd terms)
\label{prop1}
Fix $\kappa \geq 4$, $\kappa \in \N_e$, and let $n \in \N_o$ be given. System \eqref{conjecture_system_conj_intro} admits a unique solution satisfying $\mathcal{E}_n \in J_2(\kappa)$ and
\begin{equation}
\label{chain_1}
\mathcal{E}_n = X_{0} < X_{1} < X_{2} < ... <  X_{(n-1)/2} < \cos(\pi / \kappa) < X_{(n+1)/2} < ... < X_{n} < X_{n+1} := 1.
\end{equation}
This solution satisfies \eqref{o1}, \eqref{o2} and \eqref{o3}, and so $\E_n \in \mathfrak{T}_{n,\kappa}$. Furthermore it satisfies \eqref{ao1} and so $\mathcal{E}_n \in \boldsymbol{\Theta}_{(n+1)/2,\kappa} (D[d=2])$.
\end{proposition} 

\begin{proposition}  (even terms)
\label{prop2}
Fix $\kappa \geq 4$, $\kappa \in \N_e$, and let $n \in \N_e$ be given. System \eqref{conjecture_system_conj_intro_even} admits a unique solution satisfying $\mathcal{E}_n \in J_2(\kappa)$ and 
\begin{equation}
\label{chain_1even}
\mathcal{E}_n = X_{0} < X_{1} < X_{2} < ... <  X_{n/2} = \cos(\pi / \kappa) < X_{n/2+1} < ... < X_{n} < X_{n+1} := 1.
\end{equation}
This solution satisfies \eqref{o11}, \eqref{o22}, \eqref{o33} and \eqref{o44}, and so $\E_n \in \mathfrak{T}_{n,\kappa}$. Furthermore it satisfies \eqref{ae1} and so $\mathcal{E}_n \in \boldsymbol{\Theta}_{n/2,\kappa} (D[d=2])$.

\end{proposition} 

\begin{figure}[htb]
  \centering
 \includegraphics[scale=0.112]{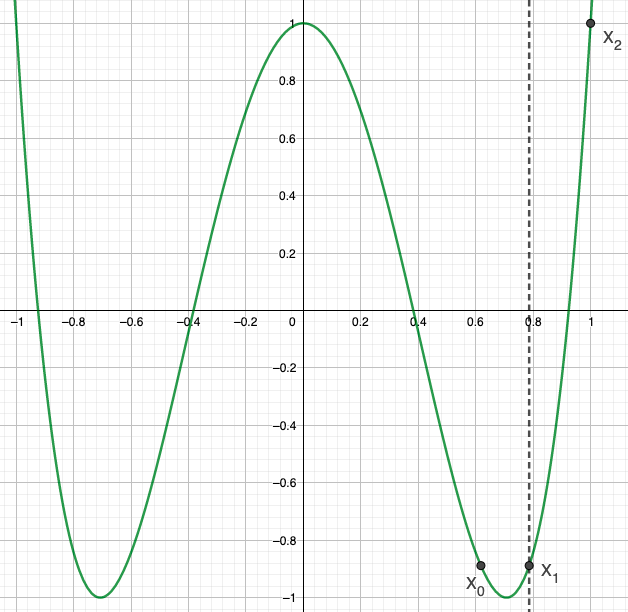}
  \includegraphics[scale=0.112]{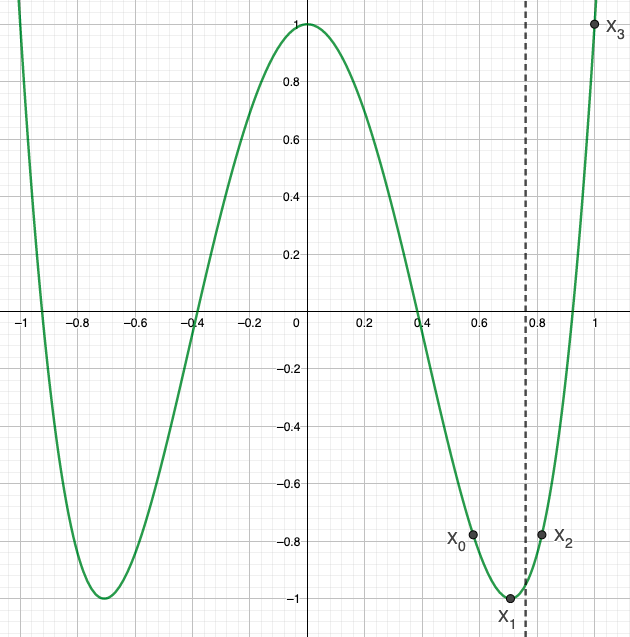}
  \includegraphics[scale=0.112]{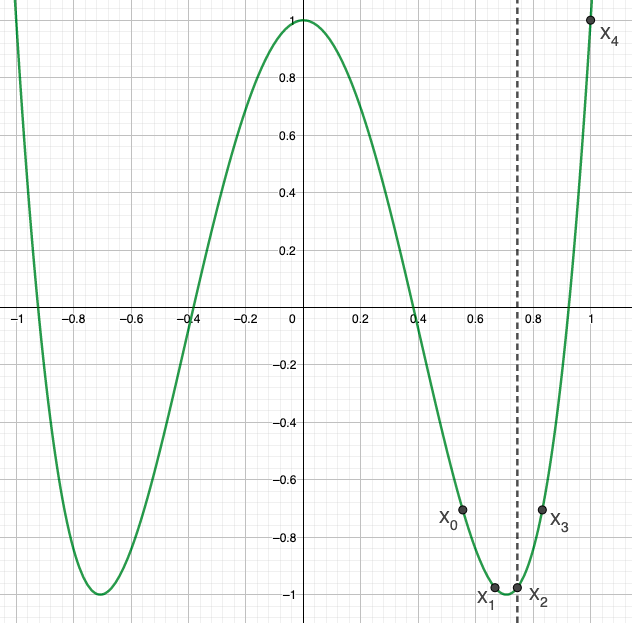}
    \includegraphics[scale=0.112]{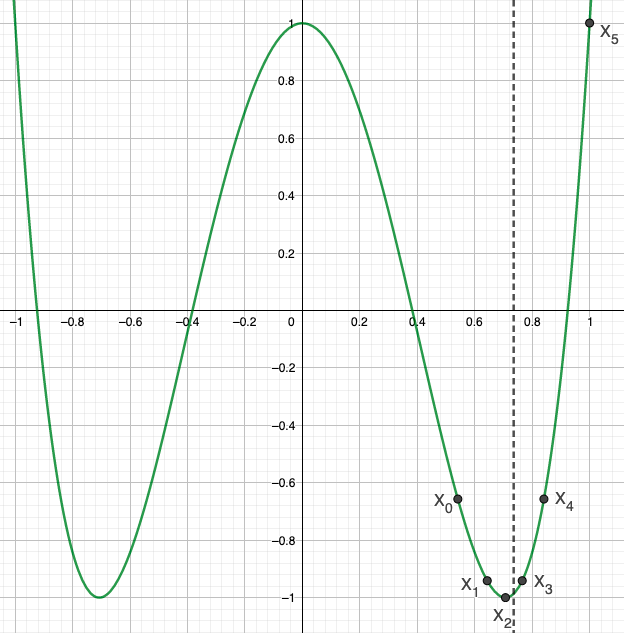}
      \includegraphics[scale=0.112]{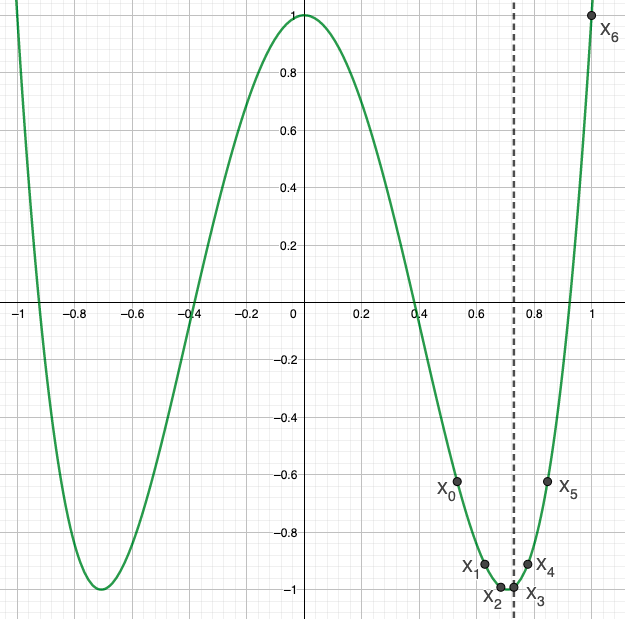}
      \includegraphics[scale=0.112]{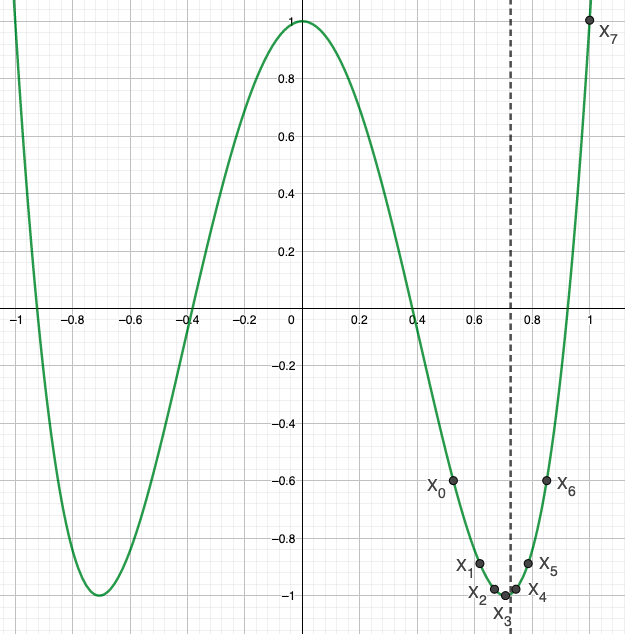}
\caption{$T_{\kappa=4}(x)$. Solutions $\E_n$ in Propositions \ref{prop1} and \ref{prop2}. Left to right : $\E_1 \simeq 0.61803$, $\E_2 \simeq 0.57735$, $\E_3 \simeq 0.55496$, $\E_4 \simeq 0.54120$, $\E_5 \simeq 0.53208$, $\E_6 \simeq 0.52573$.}
\label{fig:test_T3k3firsty}
\end{figure}

Figure \ref{fig:test_T3k3firsty} illustrates the solutions in Propositions \ref{prop1} and \ref{prop2} for $1 \leq n \leq 6$ and $\kappa=4$. This is a rare ocurrence where we have a few exact solutions. They are :

\begin{equation}
\label{exact_sols_k4}
\begin{aligned}
\E_1 &= (5^{1/2}-1)/2 \simeq 0.61803, \text{(so-called golden ratio conjugate)}, \\
\E_2 &= 3^{-1/2} \simeq 0.57735, \\
\E_3 &= 2/3 - 3^{-1} 7^{1/2}  \cos \left( 3^{-1} \arctan(3 \sqrt{3}) \right) + (7/3)^{1/2} \sin \left( 3^{-1} \arctan (3\sqrt{3}) \right) \simeq 0.55496, \\
\E_4 &= (1-2^{-1/2})^{1/2} \simeq 0.54120, \\
\E_5 &= 2 \cos \left( 2\pi / 9 \right) -1 \simeq 0.53208, \\
\E_6 &= ((5-5^{1/2})/10)^{1/2} \simeq 052573, \\
\E_7 &\simeq 0.52111, \\
\E_8 &= (2-3^{1/2})^{1/2} \simeq 0.51764.
\end{aligned}
\end{equation}

\begin{figure}[htb]
  \centering
 \includegraphics[scale=0.105]{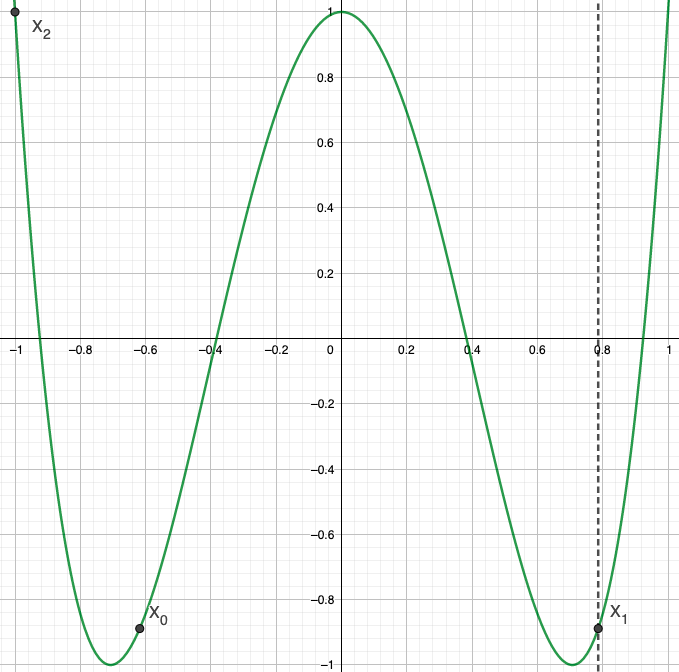}
  \includegraphics[scale=0.105]{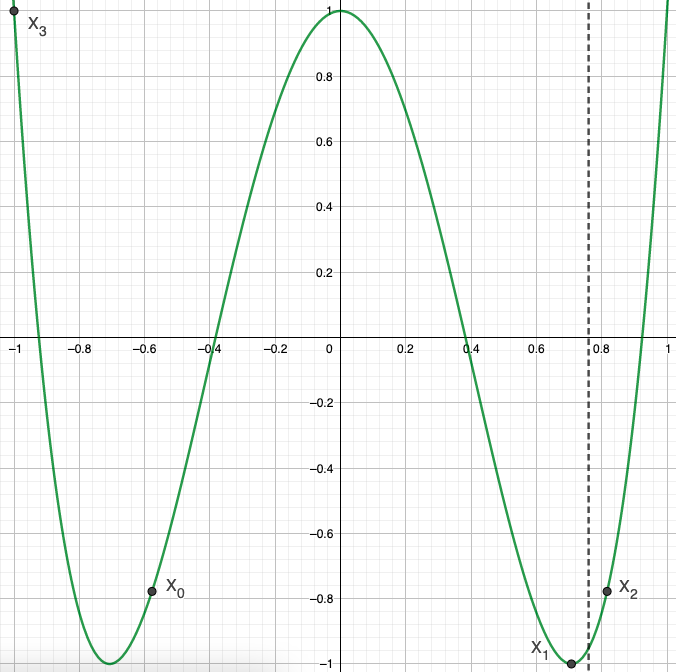}
  \includegraphics[scale=0.105]{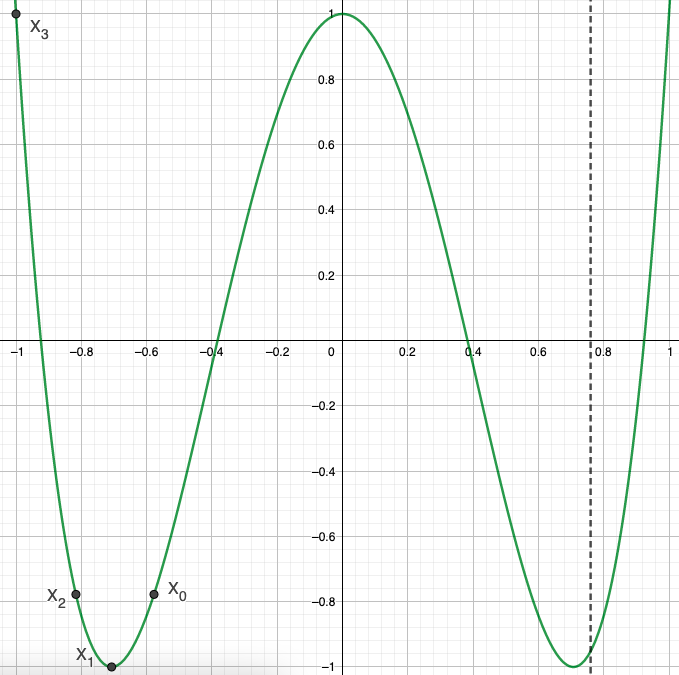}
   \includegraphics[scale=0.105]{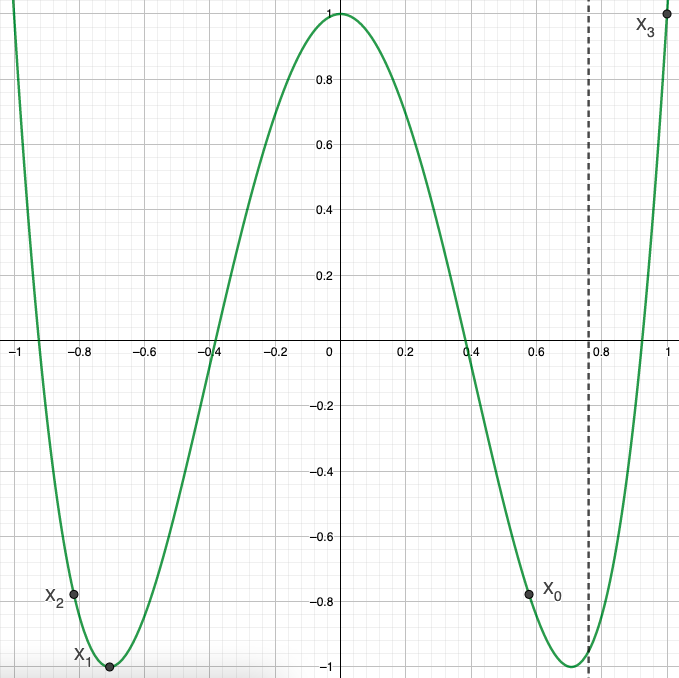}
      \includegraphics[scale=0.105]{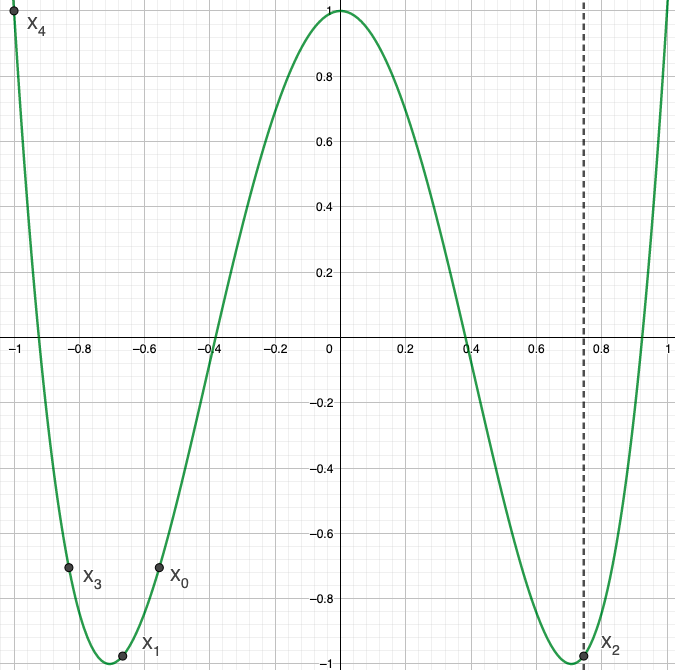}
      \includegraphics[scale=0.105]{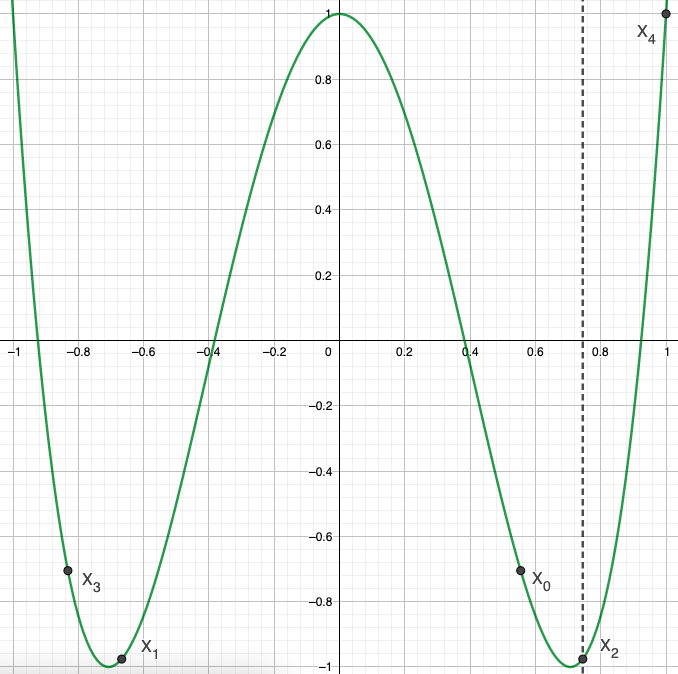}
\caption{$T_{\kappa=4}(x)$. Same solutions $\E_n$ as in Figure \ref{fig:test_T3k3firsty} but using different configurations for the $X_q$'s. Left to right : $\E_1 \simeq 0.61803$, $\E_2 \simeq 0.57735$, $\E_2 \simeq 0.57735$, $\E_2 \simeq 0.57735$, $\E_3 \simeq 0.55496$, $\E_3 \simeq 0.55496$}
\label{fig:test_T3k3firsty_versions}
\end{figure}

\begin{figure}[htb]
  \centering
 \includegraphics[scale=0.113]{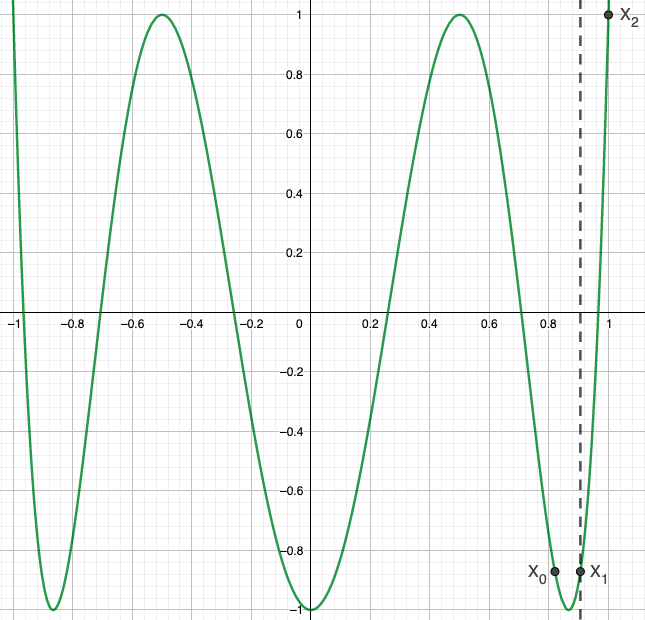}
  \includegraphics[scale=0.113]{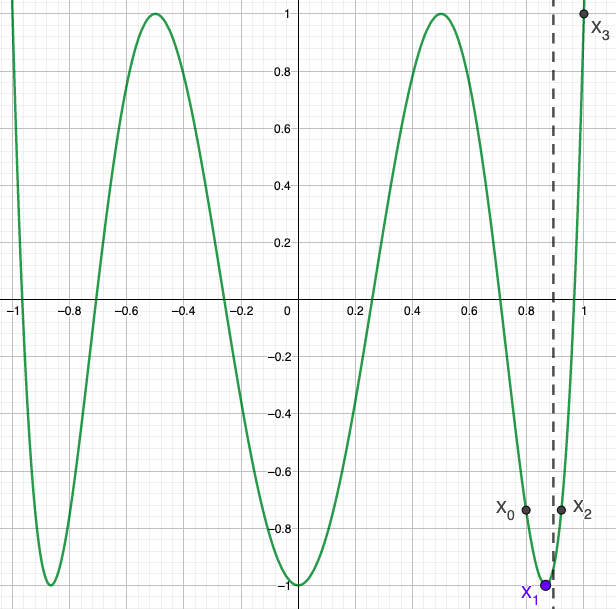}
  \includegraphics[scale=0.113]{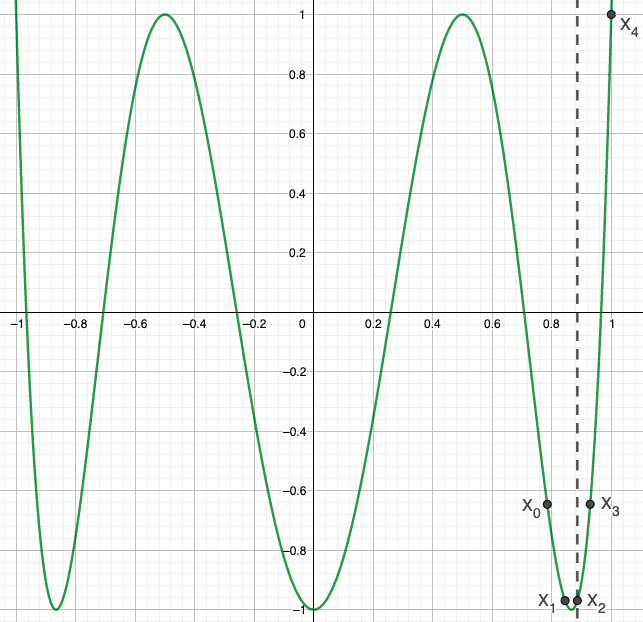}
    \includegraphics[scale=0.113]{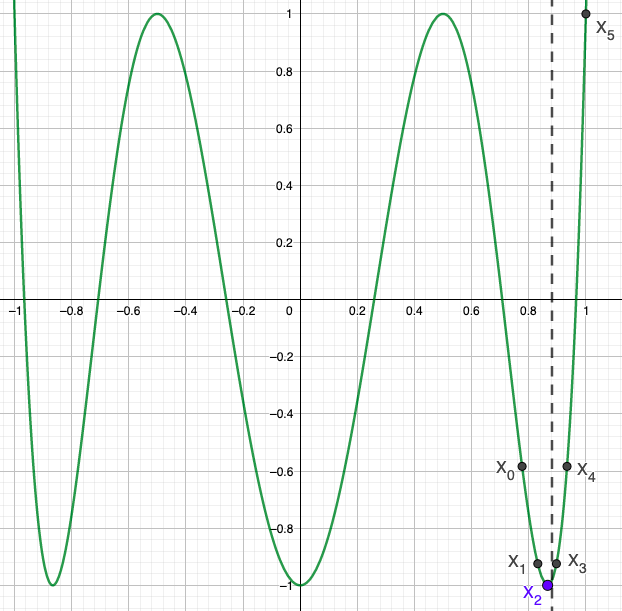}
      \includegraphics[scale=0.113]{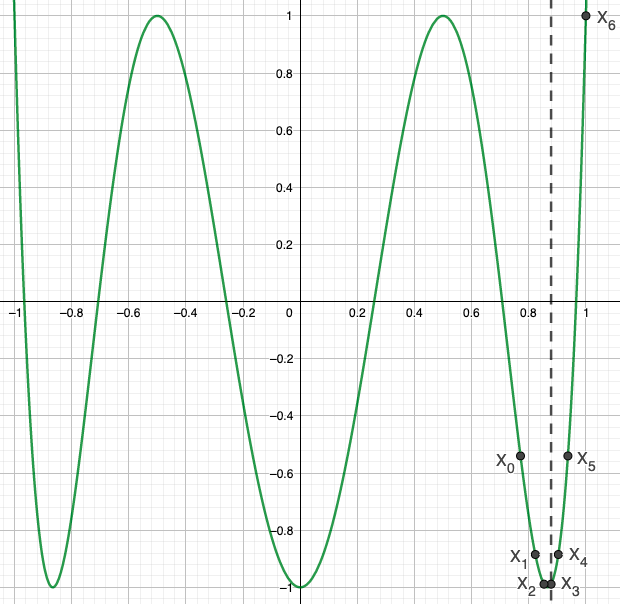}
      \includegraphics[scale=0.113]{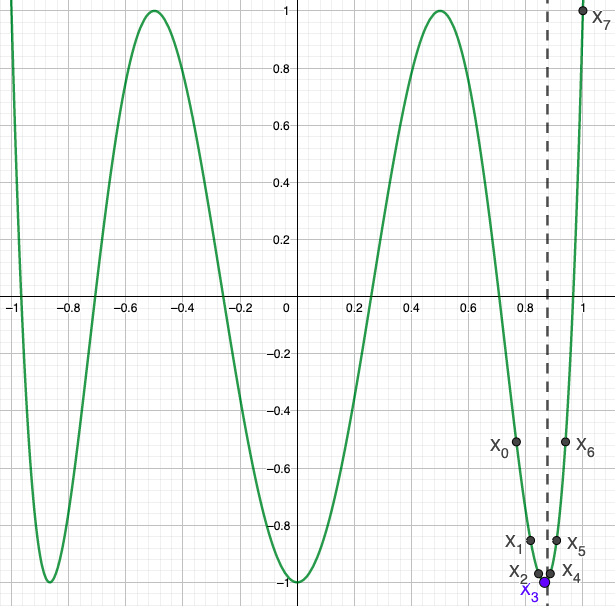}
\caption{$T_{\kappa=6}(x)$. Solutions $\E_n$ in Propositions \ref{prop1} and \ref{prop2}. Left to right : $\E_1 \simeq 0.82011$, $\E_2 \simeq 0.79770$, $\E_3 \simeq 0.78481$, $\E_4 \simeq 0.77662$, $\E_5 \simeq 0.77106$, $\E_6 \simeq 0.76710$.}
\label{fig:test_T6k3firsty}
\end{figure}

Figure \ref{fig:test_T3k3firsty_versions} depicts other configurations of the $X_q$'s that give the same threshold energy solutions as in the first graphs of Figure \ref{fig:test_T3k3firsty}. Figure \ref{fig:test_T6k3firsty} illustrates the solutions in Propositions \ref{prop1} and \ref{prop2} for $1 \leq n \leq 6$ and $\kappa=6$. 

\begin{proposition} 
\label{prop_interlace}
Fix $\kappa \geq 4$, $\kappa \in \N_e$. The odd and even energy solutions $\mathcal{E}_n$ of Propositions \ref{prop1} and \ref{prop2} interlace and are a strictly decreasing sequence : $\mathcal{E}_{n+2} < \mathcal{E}_{n+1} < \mathcal{E}_{n}$, $\forall n\in \N^*$. 
\end{proposition}

\begin{proposition}
\label{prop_convergeence} 
Let $\E_n$ be the solutions in Propositions \ref{prop1} and \ref{prop2}. Then $\E_n \searrow \cos ^2 (\pi / \kappa)$. 
\end{proposition}

Let us express the solutions $\E_n$ as solutions to a single equation for $\kappa=4,6$. To do this we need to select the appropriate branches of $T_4 ^{-1}$ and $T_6 ^{-1}$. Let
\begin{equation*}
\begin{cases}
f_E : x \mapsto \sqrt{1-(E/x)^2} &  \text{if} \ \kappa = 4, \\
f_E : x \mapsto \frac{(E/x) + \sqrt{3}\sqrt{1-(E/x)^2} }{2} & \text{if} \ \kappa = 6. 
\end{cases}
\end{equation*}

\begin{figure}[htb]
  \centering
\includegraphics[scale=0.34]{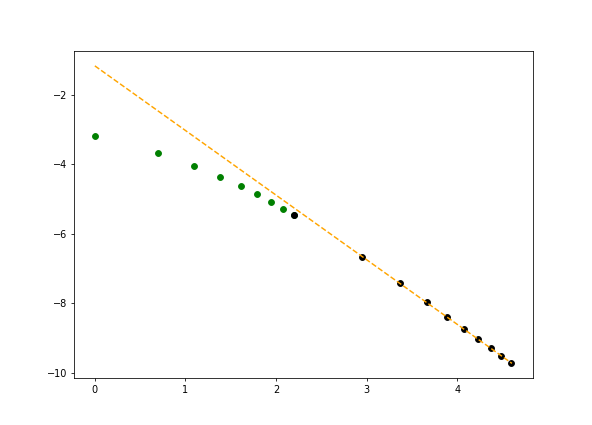}
\includegraphics[scale=0.3]{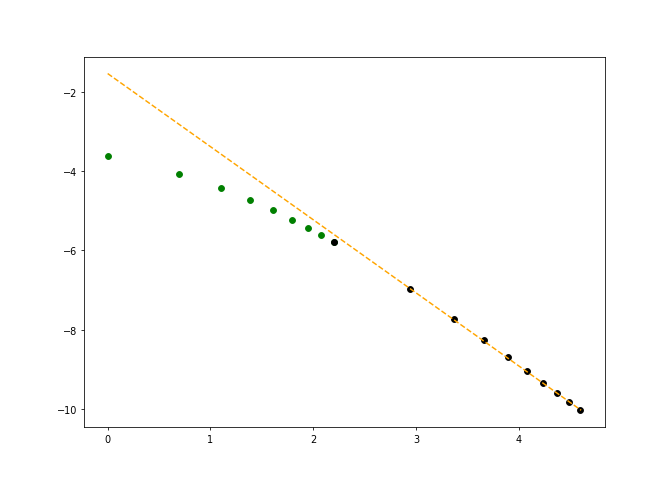}
\caption{Graphs with $\log(n)$ on $x$-axis and $\log(\E_{2n}(\kappa) - \cos ^2(\pi / \kappa))$ on $y$-axis. Left : $\kappa=4$ ; Right : $\kappa=6$. Green dots are $1 \leq n \leq 9$ ; black dots are $n= 10, 20, 30, ...,100$. Orange line is trend line based on linear regression of black dots.}
\label{fig:test_T3k3_converge}
\end{figure}

\begin{proposition}
\label{lem1_sequence3344}
Let $\kappa \in \{4,6\}$. Even terms: the $\E_{2n}$ in Theorem \ref{thm_decreasing energy general} are solution to the equation 
$\E_{2n} = f_{\E_{2n}}^{(n)}\left(\cos(\pi / \kappa) \right)$, $n \in \N^*$ ($f_{\E_{2n}}$ composed with itself $n$ times evaluated at $\cos(\pi / \kappa)$). Odd terms: the $\E_{2n-1}$ in Theorem \ref{thm_decreasing energy general} are solution to the equation $\E_{2n-1} = f_{\E_{2n-1}}^{(n)}(\sqrt{\E_{2n-1}})$. 
\end{proposition}

For example, the case of $\kappa=4$ illustrates well how the continued fractions show up. For $(\kappa,n) = (4,6)$ we have:
$$\E_6 ^2 = 1 - \frac{\E_6 ^2}{1- \frac{\E_6 ^2}{1 - \frac{\E_6 ^2}{\cos^2(\pi/4)}}}.$$

Figure \ref{fig:test_T3k3_converge} illustrates solutions $\E_{2n}(\kappa)$ of Proposition \ref{lem1_sequence3344} for $\kappa = 4,6$. In these graphs the equation of the trend line is $y = -1.872x -1.081$ and $y = -1.853x -1.457$ respectively. In particular the slope is close to $-2$, and this is our rationale behind Conjecture \ref{conjecture12}. Unfortunately I was not able to get Python to efficiently compute many more terms as in the case of the standard Laplacian in order to better approximate the trend line, but if the case of $D$ were to mirror that of $\Delta$ we should expect the trend line to steepen (closer to $-2$) as more points are plotted.

We now sequentially give the missing proofs of the aforementioned results in this section. We begin with a remark :

\begin{remark}
\label{r:ppcroissant}
Thanks to Lemma \ref{variationsTk}, 
\begin{enumerate}
\item Given $\cos(2\pi /\kappa)< a<b< \cos(\pi /\kappa)$, there exist unique
$\cos(\pi /\kappa) < b'< a'<1$ such that $T_\kappa(a)=T_\kappa(a')>T_\kappa(b')=T_\kappa(b)$. 
\item Given $\cos(\pi /\kappa) < b'< a'<1$, there exist unique
$\cos(2\pi /\kappa)< a<b< \cos(\pi /\kappa)$, such that $T_\kappa(a')=T_\kappa(a)>T_\kappa(b)=T_\kappa(b')$.
\end{enumerate}
Moreover, $(a,b)$ depends bi-continuously on $(a',b')$.
\end{remark}

\vspace{1cm}

\noindent \textit{Proof of Proposition \ref{prop1}.} 
We implement the dynamical algorithm for $n$ odd of section \ref{geo_construction}. Initialize energy $E$ to $E = E(\alpha) = \cos ^2(\pi / \kappa)+\alpha$ with $\alpha \in \mathbf{A}_{max} := (0, \cos(\pi/\kappa) - \cos^2(\pi/\kappa))$. First, by Remark \ref{obs_symmetry} we know $X_{(n+1)/2} = \sqrt{E} = \sqrt{\cos^2(\pi / \kappa) + \alpha} \in (\cos(\pi / \kappa), 1 )$. In particular when $\alpha \searrow 0_+$, note that $X_{(n+1)/2} \searrow \cos(\pi / \kappa)_+$. 
Now, up to a smaller $\alpha$ still within $\mathbf{A}_{max}$, we construct inductively and continuously in $\alpha$ all of the remaining $X_q = X_q (\alpha)$, by checking all the constraints of \eqref{conjecture_system_conj_intro}, \eqref{o1} -- \eqref{o3}, but with the exception of \eqref{o2}, i.e.\ $X_{n+1} = 1 \Leftrightarrow X_0 = E$. 

$\bullet$ $X_{(n-1)/2}$ is determined in $(\cos(2\pi / \kappa), \cos(\pi / \kappa) )$ so that $T_{\kappa}(X_{(n-1)/2}) = T_{\kappa}(X_{(n+1)/2})$. In particular, $X_{(n-1)/2}  < \cos(\pi / \kappa) < X_{(n+1)/2}$. Note that $X_{(n-1)/2} \nearrow \cos(\pi / \kappa)_-$, as $\alpha \searrow 0_+$.

$\bullet$ $X_{(n+3)/2}$ is the multiplicative symmetric of $X_{(n-1)/2}$ wrt.\ $X_{(n+1)/2} = \sqrt{E}$. So $X_{(n+1)/2} < X_{(n+3)/2}$. Up to a smaller $\alpha$ possibly, $X_{(n+3)/2} \in (\cos(\pi/\kappa),1)$. As $\alpha \searrow 0_+$, $X_{(n+3)/2} \searrow \cos(\pi / \kappa)_{+}$. 

$\bullet$ As per Remark \ref{r:ppcroissant}, $\exists ! X_{(n-3)/2} \in (\cos(2\pi / \kappa), \cos(\pi/\kappa))$ such that $ X_{(n-3)/2} <X_{(n-1)/2}$ and
$T_{\kappa}(X_{(n-1)/2}) = T_{\kappa}(X_{(n+1)/2}) < T_{\kappa}(X_{(n+3)/2})$. Again, $X_{(n-3)/2}\nearrow \cos(\pi / \kappa)_{-}$, as $\alpha \searrow 0_+$. 

$\bullet$ $X_{(n+5)/2}$ is the multiplicative symmetric of $X_{(n-3)/2}$ wrt.\ $X_{(n+1)/2} = \sqrt{E}$. Up to a smaller $\alpha$ possibly, $X_{(n+5)/2} \in (\cos(\pi / \kappa), 1)$. Since $X_{(n-3)/2} < X_{(n-1)/2}$, we infer $X_{(n+3)/2} < X_{(n+5)/2}$.  In particular, $T_{\kappa}(X_{(n-3)/2}) = T_{\kappa}(X_{(n+3)/2}) < T_{\kappa}(X_{(n+5)/2})$. Once more, $X_{(n+5)/2} \searrow \cos(\pi / \kappa)_{+}$, as $\alpha \searrow 0_+$. 

$\bullet$ We continue this ping pong game inductively till all of the $X_q =X_q(\alpha)$, $q=0,...,n+1$, have been defined. Note that the last step of the ping pong game was to place $X_{n+1}$ in such a way that it is the multiplicative symmetric of $X_0$ wrt.\  $X_{(n+1)/2} = \sqrt{E}$ (2nd line of \eqref{conjecture_system_conj_intro}). Now we consider the set $\mathbf{A}_n$ ($\mathbf{A}$ depends on $n$) of all the positive $\alpha$'s that allow a construction verifying :
\begin{equation}
\label{chain_reveal}
E \leq X_{0} < X_{1} < X_{2} < ... <  X_{(n-1)/2} < \cos(\pi / \kappa) <X_{(n+1)/2} < ... < X_{n} < X_{n+1} \leq 1.
\end{equation}
$\mathbf{A}_n \subset \mathbf{A}_{max}$ since if $\alpha \geq \cos(\pi / \kappa)-\cos^2(\pi / \kappa)$, $X_{n+1} = E / X_0 \geq \cos(\pi / \kappa) / X_0 \geq 1$. This observation will imply that $\E_n \in J_2 (\kappa)$ when the proof is over. As a side note, it is not hard to see that $\mathbf{A}_{n+2} \subset \mathbf{A}_{n}$ ; later in this section we prove $\cap_{n \in \N, \ n \ \text{odd}} \mathbf{A}_{n} = \emptyset$. It remains to argue that there is a unique $\alpha^* \in \mathbf{A}_n$ such that $X_{n+1}(\alpha^*) = 1 \Leftrightarrow X_0(\alpha^*) = E(\alpha^*)$. First, note that by construction, the chain of strict inequalities in \eqref{chain_reveal} remains valid as $\alpha$ increases in $\mathbf{A}_n$. Second, note that $X_{(n+1)/2} = \sqrt{E} = \sqrt{\cos^2(\pi / \kappa) + \alpha} < X_{n+1}$ and so $X_{n+1} \nearrow +\infty$ as $\alpha \nearrow + \infty$. Moreover, $X_{n+1}$ is strictly increasing for $\alpha \in \mathbf{A}_{n}$. Thirdly, and finally, note that by construction, as $\alpha$ increases in $\mathbf{A}_n$, $X_{n+1}$ must reach $1$ before $X_0$ reaches $\cos(2\pi / \kappa)$. This is because $T_{\kappa}(X_0) = T_{\kappa}(X_n) < T_{\kappa}(X_{n+1})$. Another way to see this is to argue by contradiction. If $X_0$ were to reach $\cos(2\pi / \kappa)$ before $X_{n+1}$ reaches $1$, then 
$$\cos^2(\pi / \kappa) / \cos(2\pi / \kappa) \leq E / \cos(2\pi / \kappa) = E / X_0 = X_{n+1} < 1 \Rightarrow \cos^2(\pi / \kappa) < \cos(2\pi / \kappa),$$
which is a false statement. Thus, $\exists ! \alpha^*$ s.t.\ $X_{n+1}(\alpha^*) = 1$. The energy solution $\E_n$ is $E(\alpha^*)$. 
\qed


\vspace{1cm}


\noindent \textit{Proof of Proposition \ref{prop2}.} 

We implement the dynamical algorithm for $n$ even of section \ref{geo_construction}. The main difference is that this time $X_{n/2} := \cos(\pi / \kappa)$. 
It implies that the values $X_0, X_1$, ..., $X_{n/2-1}$ will belong to $(\cos(2\pi /\kappa), \cos(\pi /\kappa))$, whereas the values $X_{n/2+1}, X_{n/2+2}$, ..., $X_{n}$ will belong to $(\cos(\pi /\kappa), 1)$. $X_{n+1}$ will be placed ultimately so that it equals 1.

Initialize energy $E$ to $E = E(\alpha) = \cos^2(\pi / \kappa)+\alpha$ with $\alpha \in \mathbf{A}_{max} := (0, \cos(\pi/\kappa)-\cos^2(\pi/\kappa))$. First, $\sqrt{E} = \sqrt{\cos^2(\pi / \kappa) + \alpha} \in (\cos(\pi / \kappa), 1 )$. 
Now, up to a smaller $\alpha$ still within $\mathbf{A}_{max}$, we construct inductively and continuously in $\alpha$ all of the remaining $X_q = X_q (\alpha)$, by checking all the constraints of \eqref{conjecture_system_conj_intro_even}, \eqref{o11} -- \eqref{o33}, but with the exception of the $X_{n+1}$ condition in \eqref{o22}, i.e.\ $X_{n+1} = 1 \Leftrightarrow X_0 = E$.

$\bullet$ $X_{n/2+1}$ is the multiplicative symmetric of $X_{n/2}$ wrt.\ $\sqrt{E}$. So $X_{n/2} < \sqrt{E} < X_{n/2+1} = \cos(\pi / \kappa) + \alpha / \cos(\pi / \kappa)$. As per Remark \ref{r:ppcroissant}, $X_{n/2-1}$ is constructed in $(\cos(2\pi / \kappa), \cos(\pi / \kappa) )$ so that $T_\kappa( X_{n/2-1}) = T_\kappa (X_{n/2+1})$. We turn to $X_{n/2+2}$ which is is the multiplicative symmetric of $X_{n/2-1}$ wrt.\ $\sqrt{E}$. Up to a smaller $\alpha$ possibly, $X_{n/2+2} \in (\cos(\pi/\kappa), 1)$. 
As per Remark \ref{r:ppcroissant}, there is a unique $X_{n/2-2} \in (\cos(2\pi / \kappa), \cos(\pi/\kappa))$ such that $ X_{n/2-2} <X_{n/2-1}$ and
$T_{\kappa}(X_{n/2-2}) = T_{\kappa}(X_{n/2+2}) > T_{\kappa}(X_{n/2+1})$.

$\bullet$ We continue this ping pong game inductively till all of the $X_q =X_q(\alpha)$, $q=0,...,n+1$, have been defined. Note that the last step of the ping pong game was to place $X_{n+1}$ in such a way that it is the multiplicative symmetric of $X_0$ wrt.\  $\sqrt{E}$ (2nd line of \eqref{conjecture_system_conj_intro_even}). Now we consider the set $\mathbf{A}_n$ ($\mathbf{A}$ depends on $n$) of all the positive $\alpha$'s that allow a construction verifying :

\begin{equation}
\label{chain_1_reveal}
E \leq X_{0} < X_{1} < X_{2} < ... <  X_{n/2} = \cos(\pi / \kappa) < X_{n/2+1} < ... < X_{n} < X_{n+1} \leq 1.
\end{equation} 
$\mathbf{A}_n \subset \mathbf{A}_{max}$ since if $\alpha \geq \cos(\pi / \kappa)-\cos^2(\pi / \kappa)$, $X_{n+1} = E / X_0 \geq  \cos(\pi / \kappa) / X_0 \geq 1$. As a side note, it is not hard to see that $\mathbf{A}_{n+2} \subset \mathbf{A}_{n}$ ; later in this section we prove $\cap_{n \in \N, \ n \ \text{even}} \mathbf{A}_{n} = \emptyset$. It remains to argue that there is a unique $\alpha^* \in \mathbf{A}_n$ such that $X_{n+1}(\alpha^*) = 1 \Leftrightarrow X_0(\alpha^*) = E(\alpha^*)$. First, note that by construction, the chain of strict inequalities in \eqref{chain_1_reveal} remains valid as $\alpha$ increases in $\mathbf{A}_n$. Second, note that $\sqrt{E} = \sqrt{\cos^2(\pi / \kappa) + \alpha} < X_{n/2+1} < X_{n+1}$ and so $X_{n+1} \nearrow +\infty$ as $\alpha \nearrow + \infty$. Moreover, $X_{n+1}$ is strictly increasing for $\alpha \in \mathbf{A}_{n}$. Thirdly, and finally, note that by construction, as $\alpha$ increases in $\mathbf{A}_n$, $X_{n+1}$ must reach $1$ before $X_0$ reaches $\cos(2\pi / \kappa)$ (see the previous proof for the argument). Thus, $\exists ! \alpha^*$ s.t.\ $X_{n+1}(\alpha^*) = 1$. The energy solution $\E_n$ is $E(\alpha^*)$.
\qed

\vspace{1cm}

\noindent \textit{Proof of Proposition \ref{prop_interlace}.} 

Fix $n$ odd. So $n+1$ is even. Fix $E:= \min(\mathcal{E}_n, \E_{n+1})$ (we suppose at this point that we don't know which of the 2 energies is smaller) with $\E_n$ and $\E_{n+1}$ determined as in the proofs of Propositions \ref{prop1} and \ref{prop2} respectively. The construction gives $(X_{i,n} (E) )_{i=0} ^{n+1}$ satisfying \eqref{chain_reveal} and $(X_{i,n+1} (E) )_{i=0} ^{n+2}$ satisfying \eqref{chain_1_reveal}. By the choice of $E$ we either have $X_{n+1,n} (E) = 1$ or $X_{n+2,n+1} (E) = 1$. This is to be determined. Starting from the bottom of the well we see that :
\[X_{(n-1)/2,n} (E) <X_{(n+1)/2,n+1}(E) = \cos(\pi / \kappa) < \sqrt{E} = X_{(n+1)/2,n} (E) < X_{(n+1)/2+1,n+1}(E).\]
By the ping pong game that ensues, and using Remark \ref{r:ppcroissant}, we inductively infer 
\[X_{(n+1)/2+q,n} (E) < X_{(n+1)/2+q+1,n+1}(E), \quad  \mbox{ for } q=0,1,...,(n+1)/2.\]
So $X_{n+1,n} (E) < X_{n+2,n+1}(E)$. It must be therefore that $X_{n+2,n+1}(E) = 1$ and so $E = \E_{n+1} \leq \E_n$. Furthermore, $X_{n+1,n} (E) < X_{n+2,n+1}(E)$ implies $\mathcal{E}_{n+1} < \mathcal{E}_{n}$.

Fix $n$ even. So $n+1$ is odd. We proceed with the same setup as before. Fix $E:= \min(\mathcal{E}_n, \E_{n+1})$ with $\E_n$ and $\E_{n+1}$ determined as in the proofs of Propositions \ref{prop2} and \ref{prop1} respectively. The construction gives $(X_{i,n} (E) )_{i=0} ^{n+1}$ satisfying \eqref{chain_1_reveal} and $(X_{i,n+1} (E) )_{i=0} ^{n+2}$ satisfying \eqref{chain_reveal}. By the choice of $E$ we either have $X_{n+1,n} (E) = 1$ or $X_{n+2,n+1} (E) = 1$. This is to be determined. Starting from the bottom of the well we see that :
\[X_{n/2,n+1}(E) < X_{n/2,n} (E) = \cos(\pi / \kappa) < \sqrt{E} = X_{n/2+1,n+1}(E) < X_{n/2+1,n} (E) \]
By the ping pong game that ensues, and using Remark \ref{r:ppcroissant}, we inductively infer 
\[X_{n/2+q,n} (E) < X_{n/2+q+1,n+1}(E), \quad  \mbox{ for } q=0,1,...,n/2+1.\]
So $X_{n+1,n} (E) < X_{n+2,n+1}(E)$. It must be therefore that $X_{n+2,n+1}(E) = 1$ and so $E = \E_{n+1} \leq \E_n$. Furthermore, $X_{n+1,n} (E) < X_{n+2,n+1}(E)$ implies $\mathcal{E}_{n+1} < \mathcal{E}_{n}$.
\qed

Finally, to prove Proposition \ref{prop_convergeence}, we'll start with a Lemma which characterizes a geometric property of the graph of $T_{\kappa}$ :

\begin{Lemma}
\label{l:croissdiff}
Let $\kappa\geq 2$. If $\cos(2\pi/\kappa) < a<\cos(\pi/\kappa)<b < 1$ are such that $T_\kappa(a) = T_\kappa(b)$, then 
\begin{equation}
\label{distanceLeftGreaterthanRight}
\cos(\pi/\kappa)-a>b-\cos(\pi/\kappa).
\end{equation}
\end{Lemma}
We refer to \cite{GM3} for the proof of Lemma \eqref{l:croissdiff}.

\begin{Lemma} \label{lemfunction}
Fix $\kappa \in \N^*$. For all $E \in (\cos^2(\pi / \kappa) , 1]$, $2(\sqrt{E} - \cos(\pi / \kappa) ) > E - \cos^2(\pi / \kappa)$
\end{Lemma}
\begin{proof}
Let $f(E) = 2(\sqrt{E} - \cos(\pi / \kappa) ) - (E - \cos^2(\pi / \kappa))$. Then $f(\cos^2(\pi / \kappa)) =0$ and $f'(E) >0$ for $E \in [\cos^2(\pi / \kappa) , 1)$ This implies $f(E) > 0$ for $E \in (\cos^2(\pi / \kappa) , 1]$.
\qed
\end{proof}

\begin{Lemma} \label{lemfunction2}
Let $E >0$ be given. Then for every $t \in (0, \sqrt{E})$ there is a unique $t' >0$ such that $E = (\sqrt{E} - t)(\sqrt{E} + t')$. Moreover $t' > t$ always holds.
\end{Lemma}
\begin{proof}
$E = (\sqrt{E} - t)(\sqrt{E} + t') \Leftrightarrow \sqrt{E}(t'-t) = t \cdot t' \Rightarrow t' > t$. Existence and uniqueness of $t'$ is straightforward.
\qed
\end{proof}


\noindent \textit{Proof of Proposition \ref{prop_convergeence}.}  

By Proposition \ref{prop_interlace}, and since $\E_n> \cos ^2(\pi / \kappa)$, $\exists \ell$ such that $\E_n \searrow \ell \geq  \cos ^2(\pi / \kappa)$. It is enough to show that
$\E_{2n+1} \to  \cos ^2(\pi / \kappa)$, $n \in \N^*$. Therefore, we suppose that $n$ is odd. We proceed by contradiction. Suppose $2\epsilon:= \ell-  \cos ^2 (\pi / \kappa)>0$. Recall $X_{(n+1)/2} (\E_n) = \sqrt{\E_n}$. Then Lemma \ref{lemfunction} implies that for all $n \geq 1$ and odd, $X_{(n+1)/2} (\E_n) - \cos(\pi/\kappa) > \epsilon$. Choose $n$ odd large enough so that $n \epsilon > 1$. 

By Lemma \ref{l:croissdiff},  $\cos(\pi/\kappa) - X_{(n-1)/2} (\mathcal{E}_n) > X_{(n+1)/2} (\mathcal{E}_n) - \cos(\pi/\kappa) > \epsilon$. Next, 
since $X_{(n-1)/2} (\E_n)$ and $X_{(n+3)/2} (\E_n)$ are multiplicative symmetrics wrt.\ $\sqrt{\mathcal{E}_n}$, we have by Lemma \ref{lemfunction2} that 
$$\E_n = X_{(n-1)/2} (\E_n) \cdot X_{(n+3)/2} (\E_n), \quad X_{(n-1)/2} (\E_n) = (\sqrt{\E_n} - t), \ X_{(n+3)/2} (\E_n) = (\sqrt{\E_n} + t')$$ 
with $t' > t$ and so
\begin{equation*}
\begin{aligned}
X_{(n+3)/2} (\E_n) - \cos(\pi/\kappa) &= X_{(n+3)/2} (\E_n) - \sqrt{\E_n} +\sqrt{\E_n} - \cos(\pi/\kappa) \\
&> t' + \epsilon > t + \epsilon \\
& = \sqrt{\E_n} - X_{(n-1)/2} (\E_n) + \epsilon \\
&> 3\epsilon.
\end{aligned}
\end{equation*}
Again, apply Lemma \ref{l:croissdiff} to get $\cos(\pi/\kappa) - X_{(n-3)/2} (\mathcal{E}_n) > 3 \epsilon$. Continuing in this way, we end up with 
$X_{(n+q)/2} (\E_n) - \cos(\pi/\kappa) > q \epsilon$ for $q=1,3,5,...,n$. But $X_{n,n} (\mathcal{E}_n) > n\epsilon > 1$ is absurd. We conclude that $\ell= \cos^2(\pi / \kappa)$. \qed

\section{A generalization of section \ref{section J2} : a sequence $\E_n \searrow \cos ^2 (j\pi/ \kappa)$}
\label{section_gen1}

The construction used to get a sequence in the right-most well of $T_{\kappa}(x)$ in section \ref{section J2} is not specific to the right-most well. One can build a similar sequence in other wells of $T_{\kappa}(x)$.

\subsection{Decreasing sequence in upright well, $j$ odd}
\begin{figure}[htb]
  \centering
 \includegraphics[scale=0.115]{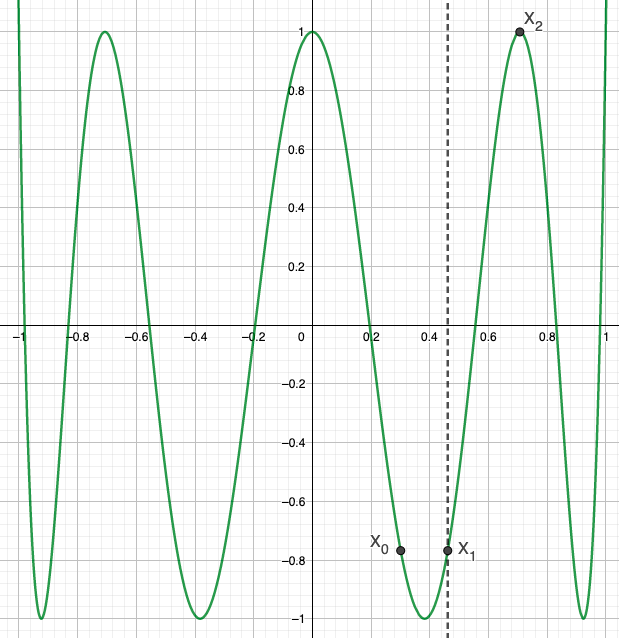}
 \includegraphics[scale=0.115]{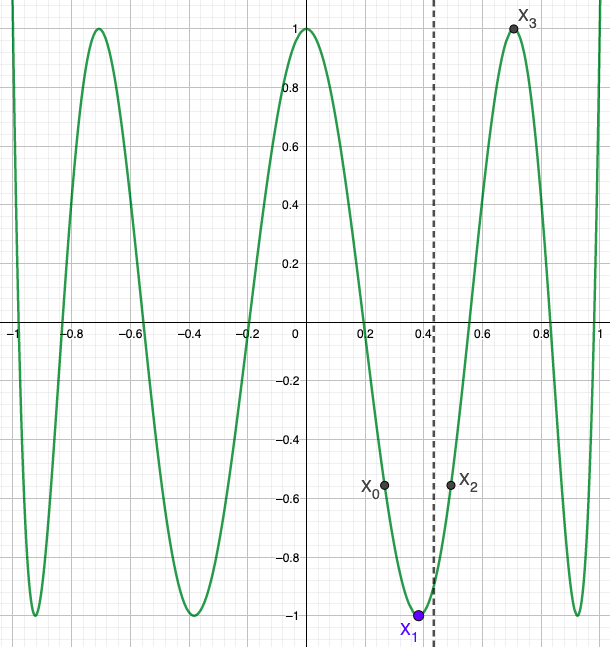}
 \includegraphics[scale=0.115]{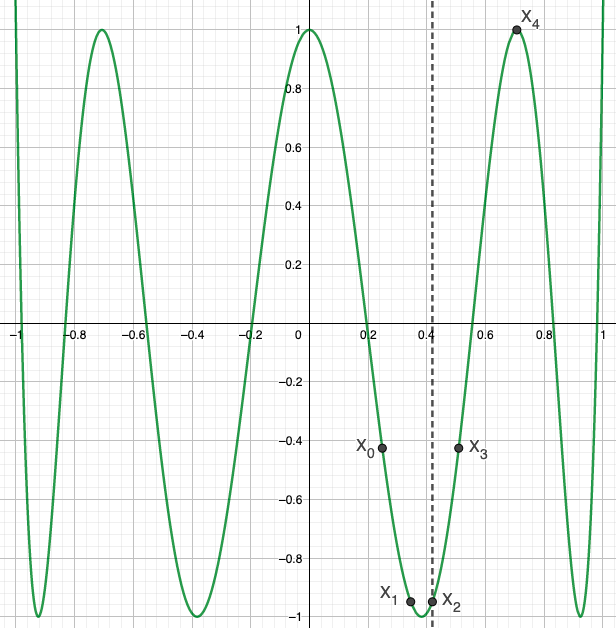}
  \includegraphics[scale=0.115]{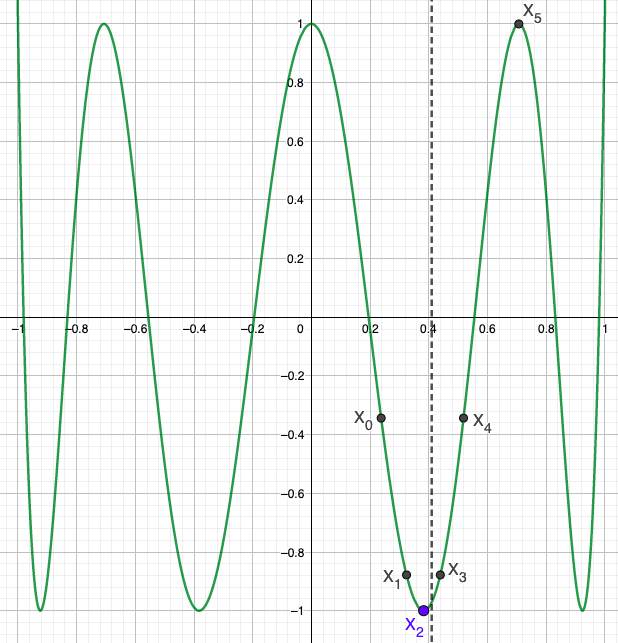}
   \includegraphics[scale=0.115]{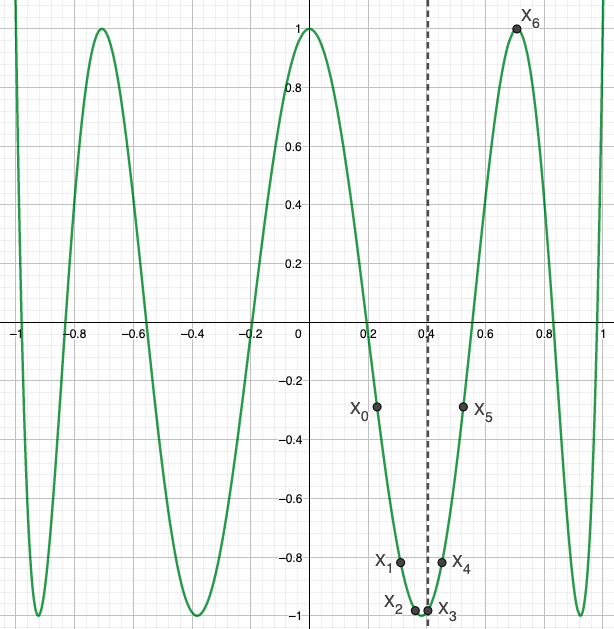}
  \includegraphics[scale=0.115]{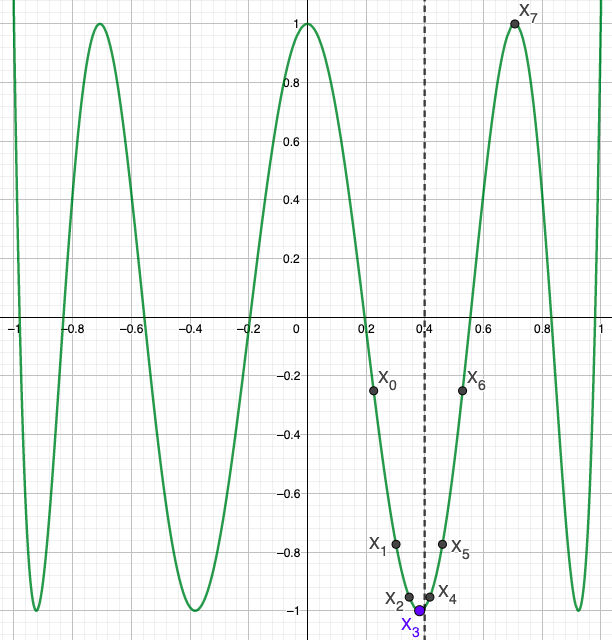}
\caption{$T_{\kappa=8}(x)$. Solutions $\E_n$. Left to right : $\E_1 \simeq 0.21289$, $\E_2 \simeq 0.18861$, $\E_3 \simeq 0.17584$, $\E_4 \simeq 0.16820$, $\E_5 \simeq 0.16325$, $\E_6 \simeq 0.15983$.}
\label{fig:test_T8decreasingNormalWell}
\end{figure}

Figure \ref{fig:test_T8decreasingNormalWell} illustrates a decreasing sequence $\E_n \searrow \cos^2(j\pi / \kappa)$, for $\kappa = 8$, $j=3$. Note that the dotted line $x=\sqrt{\E_n}$ is to the right of the minimum $x= \cos(j\pi / \kappa)$ but converges to it.

Thus, we propose a generalization of Theorem \ref{thm_decreasing energy general} :
\begin{theorem} 
\label{thm_averagekk111}
Fix $\kappa \geq 4$, $\kappa \in \N_e$. Fix $1 \leq j \leq \floor*{\kappa/2}$, $j$ odd. There is a sequence $\{ \E_n \}_{n=1} ^{\infty}$, which depends on $\kappa$, s.t.\ $\{ \E_n \} \subset \left(\cos^2(j\pi/ \kappa), \cos((j-1)\pi / \kappa) \times \cos(j\pi/ \kappa) \right) \cap \boldsymbol{\Theta}_{\kappa}(D)$, and $\E_{n+2} < \E_{n+1} < \E_n$, $\forall n \in \N^*$. Also, $\E_{2n-1}$, $\E_{2n} \in \boldsymbol{\Theta}_{n, \kappa}(D)$, $\forall n \geq 1$, and 
\begin{equation}
\label{chain_109}
\mathcal{E}_n / \cos((j-1)\pi / \kappa) = X_{0} < X_{1} < ... < X_{n} < X_{n+1} := \cos((j-1)\pi / \kappa).
\end{equation}
\end{theorem}

\subsection{Decreasing sequence in upside down well, $j$ even}

\begin{figure}[htb]
  \centering
 \includegraphics[scale=0.115]{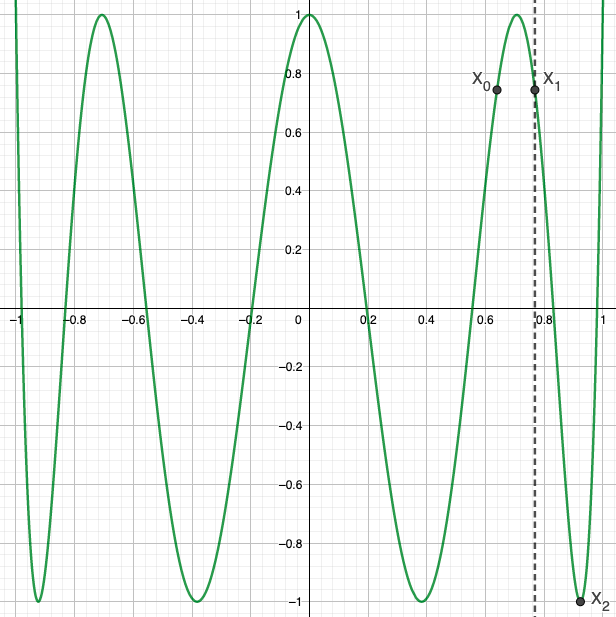}
 \includegraphics[scale=0.115]{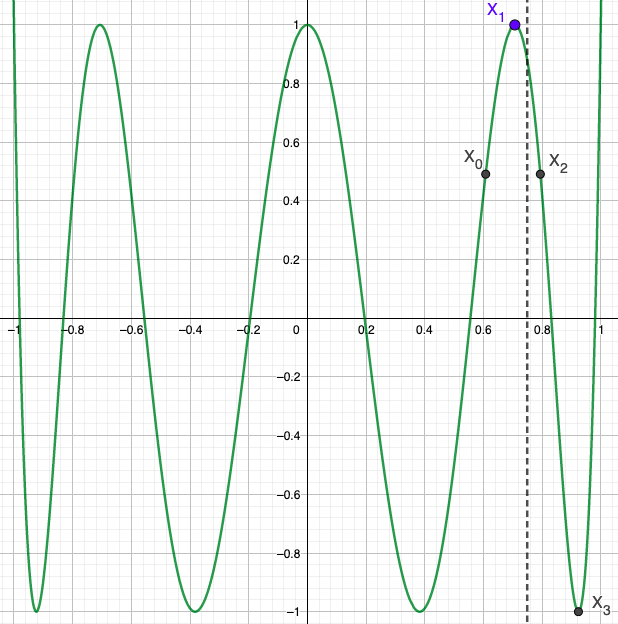}
 \includegraphics[scale=0.115]{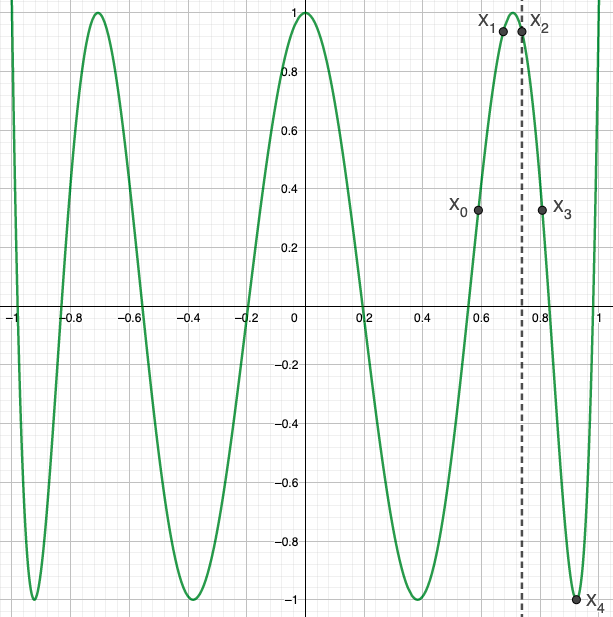}
  \includegraphics[scale=0.115]{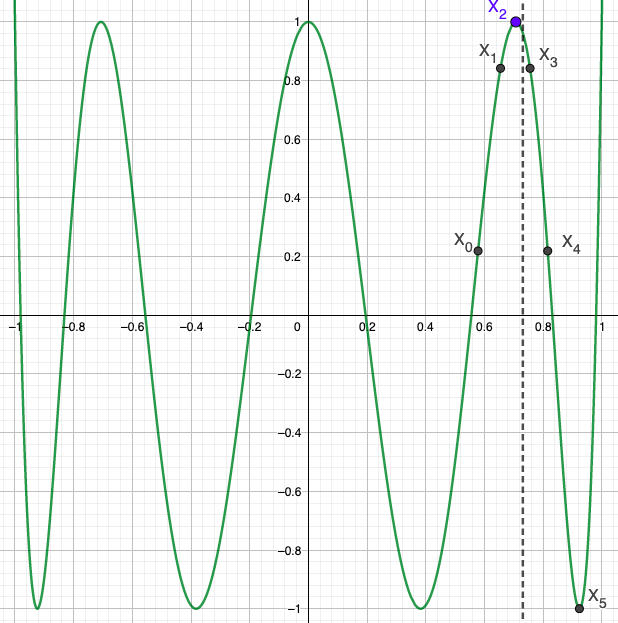}
    \includegraphics[scale=0.115]{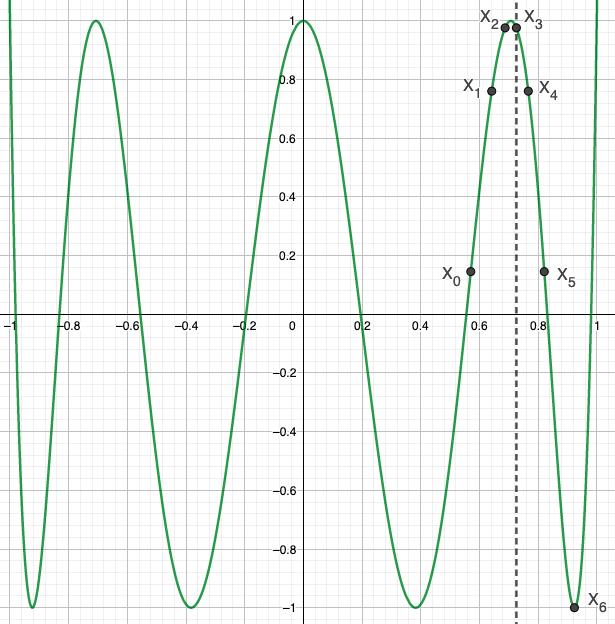}
  \includegraphics[scale=0.115]{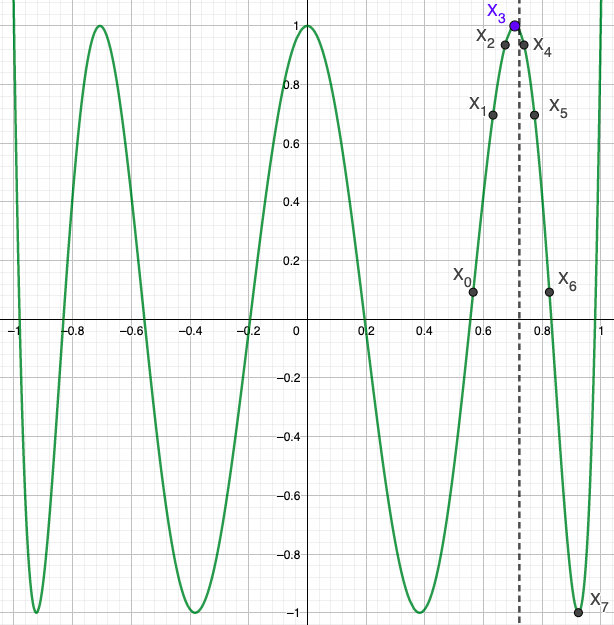}
\caption{$T_{\kappa=8}(x)$. Solutions $\E_n$. Left to right : $\E_1 \simeq 0.59091$, $\E_2 \simeq 0.56152$, $\E_3 \simeq 0.54484$, $\E_4 \simeq 0.53432$, $\E_5 \simeq 0.52720$, $\E_6 \simeq 0.52212$.}
\label{fig:test_T8decreasingUpsideDownWell}
\end{figure}

For $j$ even, the well is upside down. Figure \ref{fig:test_T8decreasingUpsideDownWell} illustrates a decreasing sequence $\E_n \searrow \cos ^2(j\pi / \kappa)$, for $\kappa = 8$, $j=2$. Note that the dotted line $x=\sqrt{\E_n}$ is to the right of the maximum $x= \cos(j\pi / \kappa)$ but converges to it.

Thus, we propose a generalization of Theorem \ref{thm_decreasing energy general} :
\begin{theorem} 
\label{thm_averagekk112}
Fix $\kappa \geq 4$, $\kappa \in \N_e$. Fix $2 \leq j \leq \floor*{\kappa/2}$, $j$ even. There is a sequence $\{ \E_n \}_{n=1} ^{\infty}$, which depends on $\kappa$, s.t.\ $\{ \E_n \} \subset \left(\cos^2 (j\pi/ \kappa), \cos((j-1)\pi / \kappa) \times \cos(j\pi/ \kappa) \right) \cap \boldsymbol{\Theta}_{\kappa}(D)$, and $\E_{n+2} < \E_{n+1} < \E_n$, $\forall n \in \N^*$. Also, $\E_{2n-1}$, $\E_{2n} \in \boldsymbol{\Theta}_{n, \kappa}(D)$, $\forall n \geq 1$, and 
\begin{equation}
\label{chain_1099}
\mathcal{E}_n / \cos((j-1)\pi / \kappa) = X_{0} < X_{1} < ... < X_{n} < X_{n+1} := \cos((j-1)\pi / \kappa).
\end{equation}
\end{theorem}

\subsection{A comment on the proofs of these Theorems and a Conjecture on the limit}

To prove Theorems \ref{thm_averagekk111} and \ref{thm_averagekk112} one needs to adapt the proofs of Propositions \ref{prop1}, \ref{prop2} and \ref{prop_interlace}. The adaptation of these Propositions is straightforward. 
Moreover, to see why $\{ \E_n \} \subset (\cos^2(j\pi/ \kappa), \cos((j-1)\pi / \kappa) \times \cos(j\pi/ \kappa))$, note that by the construction 
$$\cos^2(j\pi/ \kappa) < X_{(n+1)/2} ^2 = \E_n, \quad  X_0 \times X_{n+1} = \E_n <  \cos(j\pi/ \kappa) \times \cos((j-1)\pi / \kappa).$$

As for the limit we conjecture :
\begin{conjecture}
\label{conj_JJJ}
Let $\{\E_n\}$ be the sequence in Theorems \ref{thm_averagekk111} and \ref{thm_averagekk112}. Then $\E_n \searrow \cos ^2(j\pi/ \kappa)$.
\end{conjecture}

As explained in \cite{GM3} we don't know how to adapt the proof of Lemma \ref{l:croissdiff} in order to prove Conjecture \ref{conj_JJJ}. We do conjecture :
\begin{conjecture}
\label{l:croissdiffJJ}
Let $\kappa\geq 2$, $\kappa \in \N_e$. Fix $1 \leq j \leq \floor*{\kappa/2}$. If $\cos((j+1)\pi/\kappa)<a<\cos(j\pi/\kappa)<b<\cos((j-1)\pi/\kappa)$ are such that $T_\kappa(a) = T_\kappa(b)$, then 
\begin{equation}
\label{distanceLeftGreaterthanRightJJ}
\cos(j\pi/\kappa)-a>b-\cos(j\pi/\kappa).
\end{equation}
\end{conjecture}
If Conjecture \ref{l:croissdiffJJ} holds, Conjecture \ref{conj_JJJ} should follow directly.



\section{An increasing sequence of thresholds below $J_3(\kappa) :=\left( \cos(2 \pi / \kappa), \cos^2(\pi / \kappa) \right)$}
\label{section J3a}

This entire section is in dimension 2. We prove the existence of a sequence of threshold energies $\F_n = \F_n (\kappa) \nearrow \inf J_3 (\kappa)$. This section proves Theorem \ref{thm_decreasing energy general_k3} for $\{ \F_n \}$.

\begin{proposition} (odd terms)
\label{prop_k=3_odd_888}
Fix $\kappa \geq 6$, $\kappa \in \N_e$, and $n \in \N_o$. System \eqref{conjecture_system_conj_intro} has a unique solution (denoted $\F_n$ instead of $\E_n$) such that $\F_n \in (\cos(\pi / \kappa) \times \cos(2\pi / \kappa), \cos(2\pi / \kappa))$ and
\begin{equation}
\label{chain_1_k3_888}
\cos(2\pi / \kappa) =: X_{n+1} < X_{n}  < ... <  X_{(n+1)/2} = \sqrt{\F_n }  < X_{(n-1)/2} < ...  < X_0 < 1,
\end{equation}
This solution satisfies \eqref{o1}, \eqref{o2} and \eqref{o3}, and so $\F_n \in \mathfrak{T}_{n,\kappa}$. Furthermore it satisfies \eqref{ao1} and so $\F_n \in \boldsymbol{\Theta}_{(n+1)/2,\kappa} (D[d=2])$.
\end{proposition} 

\begin{proposition} (even terms)
\label{prop_k=3_even_888}
Fix $\kappa \geq 6$, $\kappa \in \N_e$, and $n \in \N_e$. System \eqref{conjecture_system_conj_intro_even} has a unique solution (denoted $\F_n$ instead of $\E_n$) such that $\F_n \in (\cos(\pi / \kappa) \times \cos(2\pi / \kappa), \cos(2\pi / \kappa))$ and
\begin{equation}
\label{chain_1_k3_888even}
\cos(2\pi / \kappa) =: X_{n+1} < X_{n}  < ... < X_{n/2} = \cos(\pi / \kappa) < X_{n/2-1} < ... < X_0 < 1.
\end{equation}
This solution satisfies \eqref{o11}--\eqref{o44} and so $\F_n \in \mathfrak{T}_{n,\kappa}$. Furthermore it satisfies \eqref{ae1} and so $\F_n \in \boldsymbol{\Theta}_{n/2,\kappa} (D[d=2])$.
\end{proposition}

Figure \ref{fig:T3, increasing to 0.5} illustrates the solutions in Propositions \ref{prop_k=3_odd_888} and \ref{prop_k=3_even_888} for $1 \leq n \leq 3$ and $\kappa=6$. Unfortunately the $X_q$'s were accumulating so quickly to $\cos(2\pi/6)$ and $1$ that I didn't find a legible way of graphing the cases $n \geq 4$.

\begin{figure}[H]
  \centering
    \includegraphics[scale=0.195]{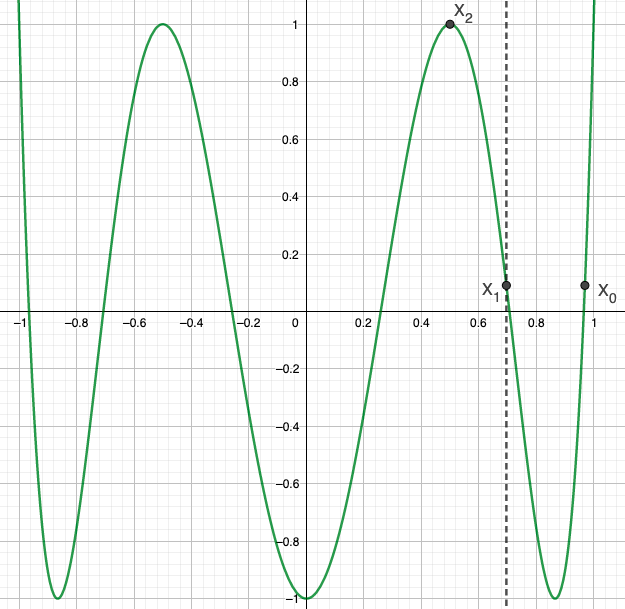}
    \includegraphics[scale=0.16]{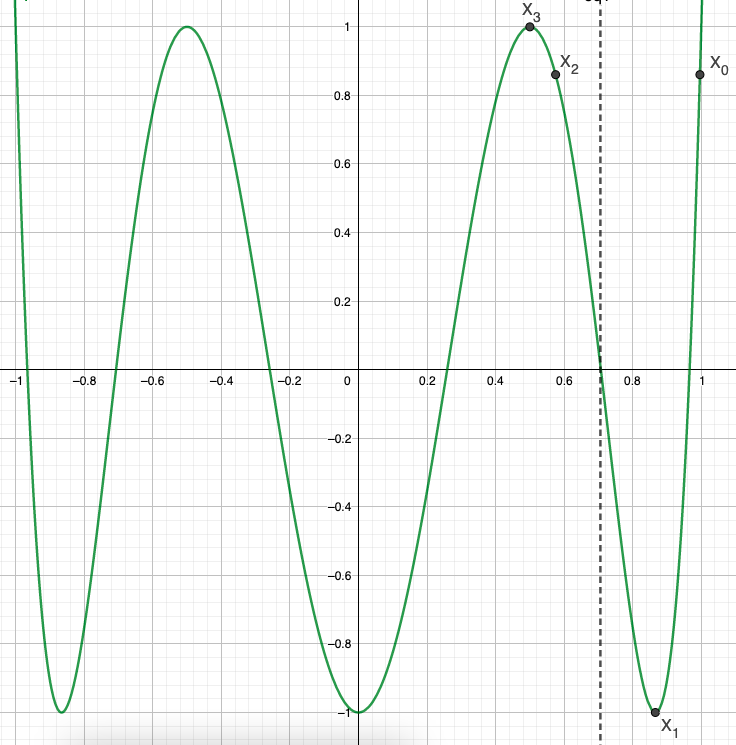}
    \includegraphics[scale=0.193]{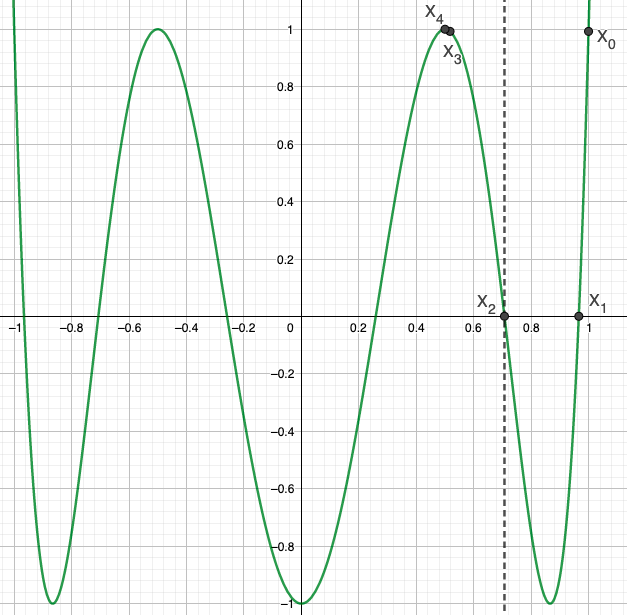}
\caption{$T_{\kappa=6}(x)$. Solution $\F_n$ of Propositions \ref{prop_k=3_odd_888} and \ref{prop_k=3_even_888}. Left to right : $\F_1 \simeq 0.48487$ ; $\F_2 \simeq 0.49802$ ; $\F_3 \simeq 0.49990$}
\label{fig:T3, increasing to 0.5}
\end{figure}


\begin{proposition} 
\label{prop_interlace_888}
Fix $\kappa \geq 6$, $\kappa \in \N_e$. The odd and even energy solutions $\mathcal{F}_n$ of Propositions \ref{prop_k=3_odd_888} and \ref{prop_k=3_even_888} interlace and are a strictly increasing sequence : $\mathcal{F}_{n} < \mathcal{F}_{n+1} < \mathcal{F}_{n+2}$, $\forall n\in \N^*$. Also, $\F_n \nearrow \inf J_3 := \cos(2\pi / \kappa)$. 
\end{proposition} 

\section{A generalization of section \ref{section J3a} : a sequence $\F_n \nearrow \cos((j-1)\pi/ \kappa)\times \cos((j+1)\pi/ \kappa)$}
\label{section_gen2}

Again, the construction used to get a sequence in the right-most well of $T_{\kappa}(x)$ in section \ref{section J3a} is not specific to the right-most well. One can build a similar sequence in other wells. In this section we get an increasing sequence $\F_n \nearrow \cos((j-1)\pi/ \kappa) \times \cos((j+1)\pi/ \kappa)$.

\subsection{Increasing sequence in upright well, $j$ odd}

\begin{theorem} 
\label{thm_averagekk113}
Fix $\kappa \geq 4$, $\kappa \in \N_e$. Fix $1 \leq j \leq \floor*{\kappa/2}$, $j$ odd. There is a sequence $\{ \F_n \}_{n=1} ^{\infty}$, which depends on $\kappa$, s.t.\ $\{ \F_n \} \subset (\cos(j\pi/ \kappa) \times \cos((j+1)\pi/ \kappa), \cos((j-1)\pi / \kappa) \times \cos((j+1)\pi/ \kappa)) \cap \boldsymbol{\Theta}_{\kappa}(D)$, and $\F_n < \F_{n+1} < \F_{n+2}$, $\forall n \in \N^*$. Also, $\F_{2n-1}$, $\F_{2n} \in \boldsymbol{\Theta}_{n, \kappa}(D)$, $\forall n \geq 1$, and 
\begin{equation}
\label{chain_1093}
\cos((j+1)\pi / \kappa) =: X_{n+1} < X_{n} < ... < X_{1} < X_{0} < \cos((j-1)\pi / \kappa).
\end{equation}
\end{theorem}

\subsection{Increasing sequence in upside down well, $j$ even}
\begin{theorem} 
\label{thm_averagekk1145}
Fix $\kappa \geq 4$, $\kappa \in \N_e$. Fix $1 \leq j \leq \floor*{\kappa/2}$, $j$ even. There is a sequence $\{ \F_n \}_{n=1} ^{\infty}$, which depends on $\kappa$, s.t.\ $\{ \F_n \} \subset (\cos(j\pi/ \kappa) \times \cos((j+1)\pi/ \kappa), \cos((j-1)\pi / \kappa) \times \cos((j+1)\pi/ \kappa)) \cap \boldsymbol{\Theta}_{\kappa}(D)$, and $\F_n < \F_{n+1} < \F_{n+2}$, $\forall n \in \N^*$. Also, $\F_{2n-1}$, $\F_{2n} \in \boldsymbol{\Theta}_{n, \kappa}(D)$, $\forall n \geq 1$, and 
\begin{equation}
\label{chain_10935}
\cos((j+1)\pi / \kappa) =: X_{n+1} < X_{n} < ... < X_{1} < X_{0} < \cos((j-1)\pi / \kappa).
\end{equation}
\end{theorem}

\begin{remark} The proofs of Theorems \ref{thm_averagekk113} and \ref{thm_averagekk1145} are analogous to those of Propositions \ref{prop1}, \ref{prop2} and \ref{prop_interlace}. Moreover, to see why $\{ \F_n \} \subset (\cos(j\pi/ \kappa) \times \cos((j+1)\pi/ \kappa), \cos((j-1)\pi / \kappa) \times \cos((j+1)\pi/ \kappa))$, note that by the construction 
$$\cos(j\pi/ \kappa) \times \cos((j+1)\pi/ \kappa) < X_0 \times X_{n+1} = \E_n < \cos((j-1)\pi / \kappa) \times \cos((j+1)\pi/ \kappa).$$
\end{remark}

\subsection{Conjecture on the limit}

We don't have a proof for the following Conjecture :

\begin{conjecture} 
Let $\{ \F_n \}$ be the sequence in Theorems \ref{thm_averagekk113} and \ref{thm_averagekk1145}. Then 
$$\F_n \nearrow \cos((j-1)\pi / \kappa) \times \cos((j+1)\pi/ \kappa).$$
\end{conjecture}

\section{The geometric construction of section \ref{geo_construction} revisited : the alignment condition}
\label{revisit}

We briefly revisit section \ref{geo_construction}. It is possible to construct even more thresholds if we tweak the assumption on $X_{n+1}$. This is the topic of this section. Instead of assumption \eqref{o3} (respectively \eqref{o33}), we require 
\begin{equation}
\label{o3bisounours}
T_{\kappa}(X_{n+1}) = T_{\kappa}(X_{p}) \quad \text{for some} \ p \in \{ 0 \leq p \leq (n-1)/2 \},\end{equation}
respectively
\begin{equation}
\label{o33bisounours}
T_{\kappa}(X_{n+1}) = T_{\kappa}(X_{p}) \quad \text{for some} \ p \in \{ 0 \leq p \leq n/2-1 \}.
\end{equation}

For $n \in \N_o$ define $\mathfrak{T}_{n,\kappa} ^+$ to be the set of non-zero $\E_n (\kappa)$ such that the system \eqref{conjecture_system_conj_intro} has a solution satisfying \eqref{o1}, \eqref{o2} and \eqref{o3bisounours} (definition for $\E_n >0$). For $n \in \N_e$ define $\mathfrak{T}_{n,\kappa} ^+ $ to be the set of non-zero $\E_n (\kappa)$ such that the system \eqref{conjecture_system_conj_intro_even} has a solution satisfying \eqref{o11}, \eqref{o22}, \eqref{o44} and \eqref{o33bisounours}.

\begin{figure}[htb]
  \centering
   \includegraphics[scale=0.195]{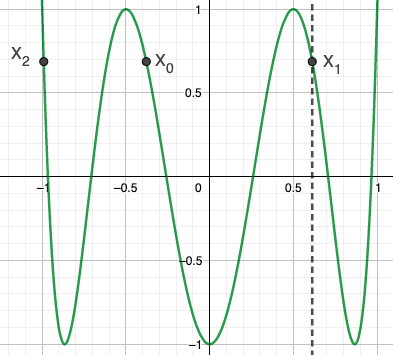}
 \includegraphics[scale=0.122]{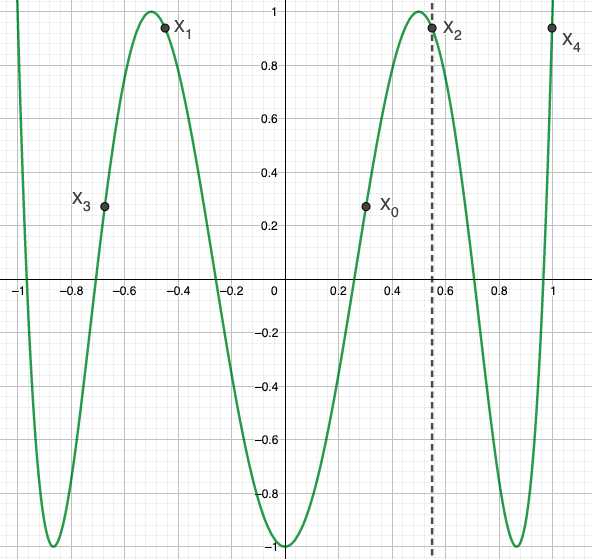}
  \includegraphics[scale=0.11]{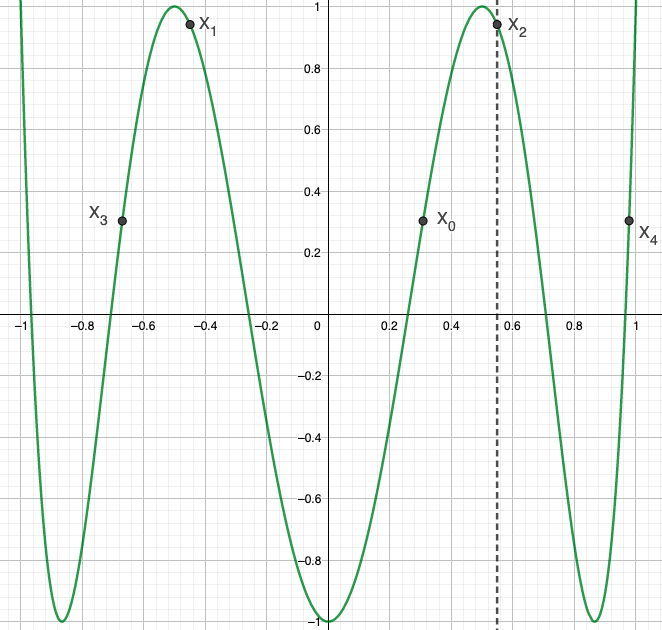}
  \includegraphics[scale=0.103]{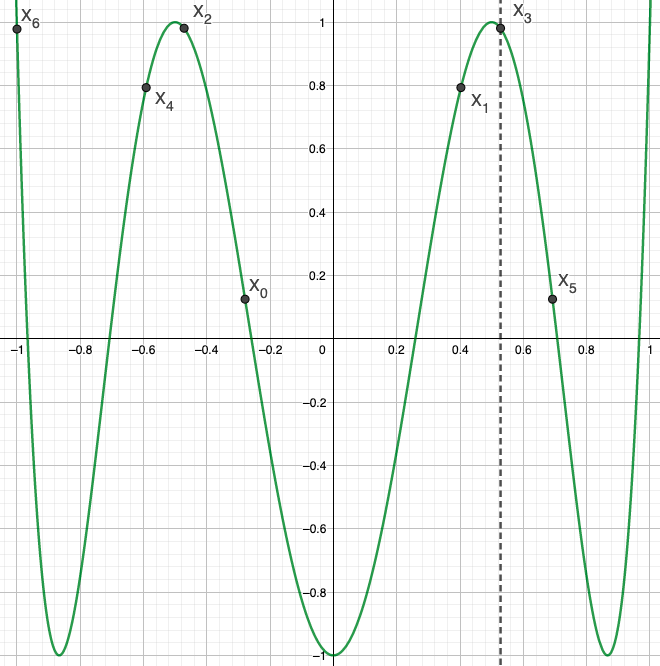}
    \includegraphics[scale=0.1]{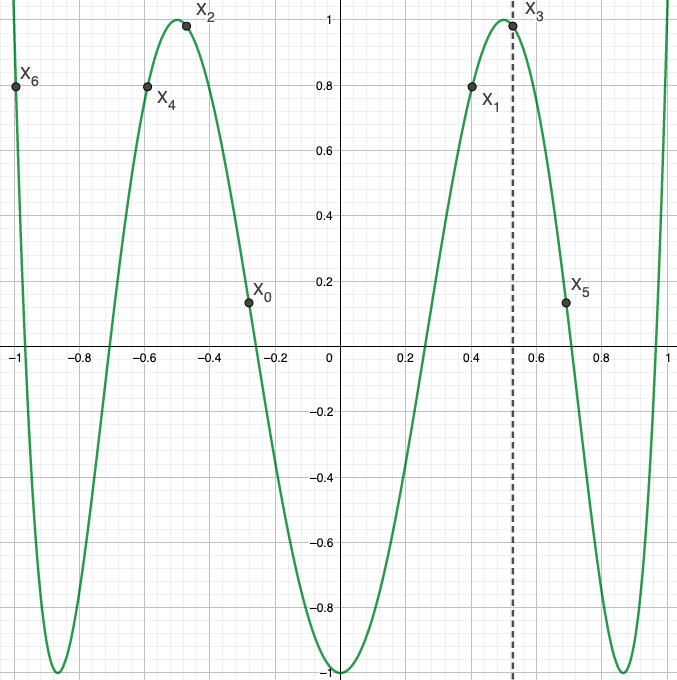}
        \includegraphics[scale=0.1]{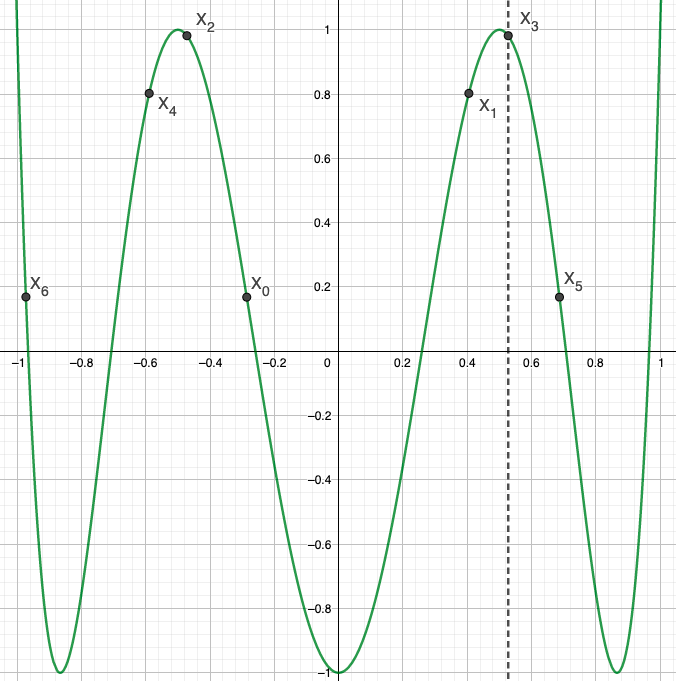}
\caption{$T_{\kappa=6}(x)$. Threshold solutions for which $X_{n+1}$ satisfies the alignment condition \eqref{o3bisounours}. Left to right : $\E_1 = 0.375$, $p=0$ ; $\E_3 \simeq 0.30227$, $p=1$ ; $\E_3 \simeq 0.30118$, $p=0$ ; $\E_5 \simeq 0.27878$, $p=2$ ; $\E_5 \simeq 0.27867$, $p=1$ ; $\E_5 \simeq 0.27819$, $p=0$}
\label{fig:example_2nd_gap2}
\end{figure}

\begin{figure}[htb]
  \centering
   \includegraphics[scale=0.129]{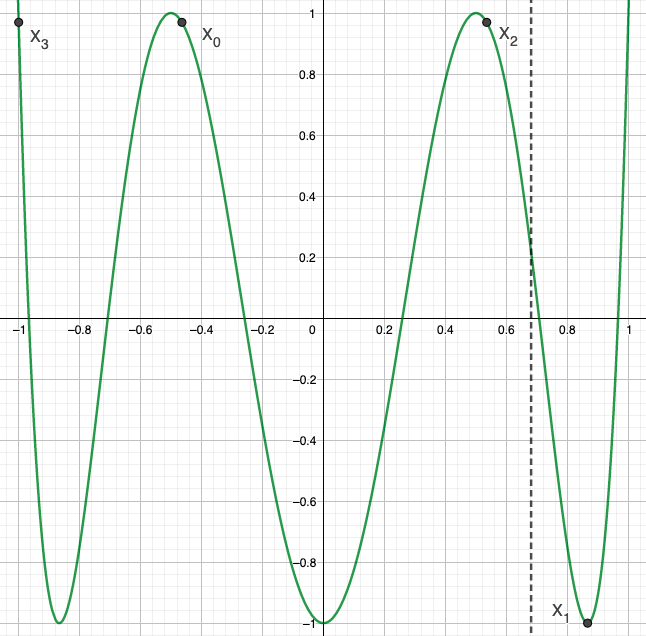}
 \includegraphics[scale=0.124]{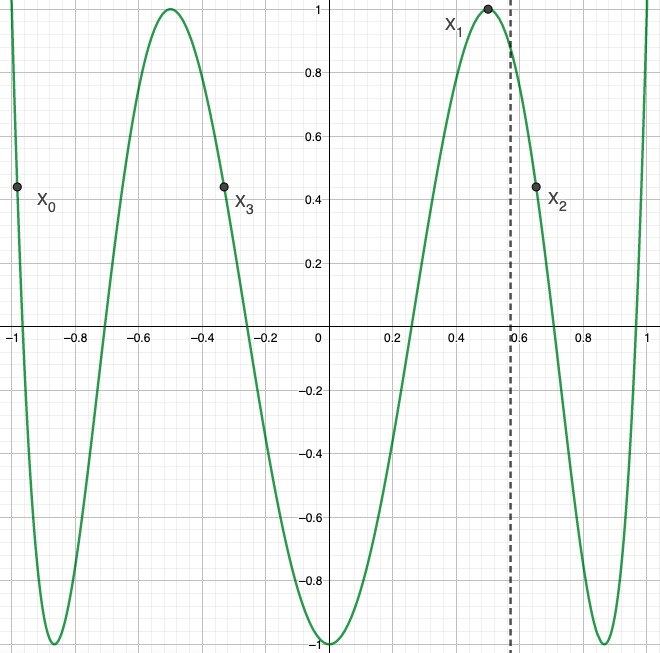}
  \includegraphics[scale=0.13]{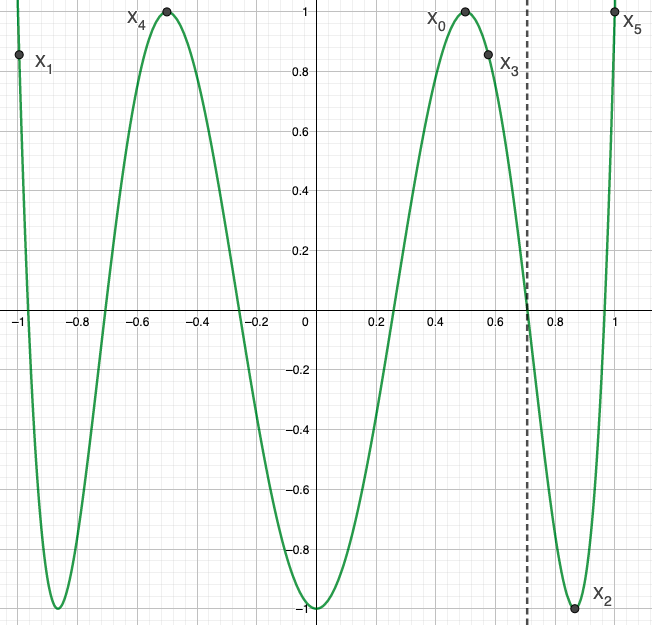}
  \includegraphics[scale=0.127]{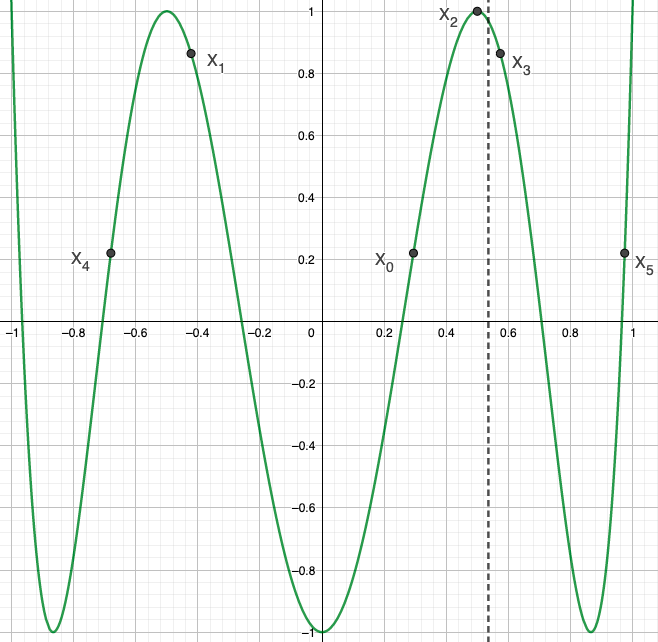}
    \includegraphics[scale=0.124]{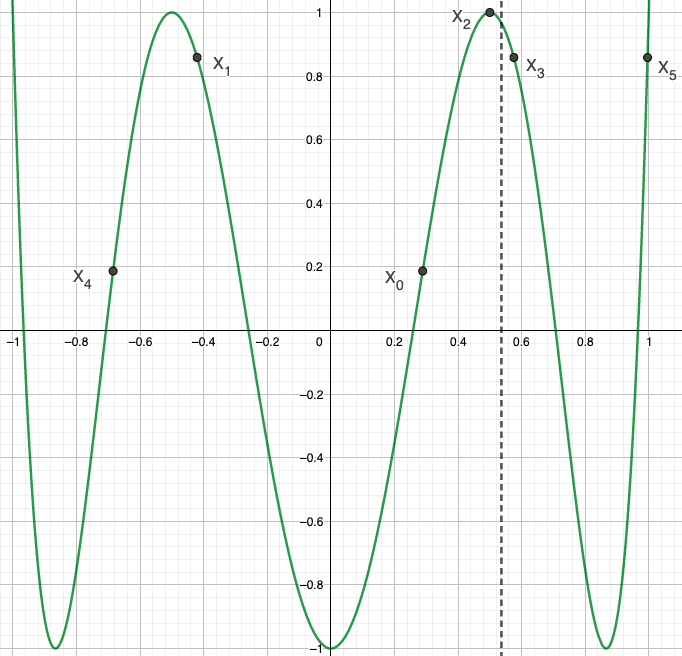}
\caption{$T_{\kappa=6}(x)$. Threshold solutions for which $X_{n+1}$ satisfies the alignment condition \eqref{o33bisounours}. Left to right : $\E_2 \simeq 0.46353$, $p=0$ ; $\E_2 \simeq 0.32569$, $p=0$ ; $\E_4 \simeq 0.49898$, $p=0$ ; $\E_4 \simeq 0.28709$, $p=0$ ; $\E_4 \simeq 0.28780$, $p=1$}
\label{fig:example_2nd_gap22}
\end{figure}

Figures \ref{fig:example_2nd_gap2} and \ref{fig:example_2nd_gap22} depict geometric constructions of such thresholds for $n=1,3,5$ and $n=2,4$ respectively. Note that we have chosen $\kappa=6$ for the illustrations and that in this case all these additional threshold energies lie in the second gap, namely $(\cos^2(2\pi / 6), \cos(2\pi / 6)) = (1/4, 1/2)$. 

\begin{proposition} ($n=1,2,3,4$)
\label{prop1_intro_go}
Fix $\kappa \in \N_e$, $\kappa \geq 4$, and let $n \in \{1,2,3,4\}$. Let $(X_q)_{q=0}^{n+1}$, $\mathcal{E}_n (\kappa)$ be a solution such that $\mathcal{E}_n (\kappa) \in \mathfrak{T}_{n,\kappa} ^+$. Then for all $j \in \N^*$,
\begin{equation}
\label{problem1_tosolve_again_101_88___}
g_{j \kappa} ^{\mathcal{E}_n}  ( X_{\frac{n+1}{2}} ) = \sum _{q =0} ^{\frac{n-1}{2}} \omega_q \cdot g_{j \kappa} ^{\mathcal{E}_n} ( X_{q}), \quad \text{for} \ n=1,3 \quad \text{and} \quad g_{j \kappa} ^{\mathcal{E}_n}  ( X_{\frac{n}{2}} ) = \sum _{q =0} ^{\frac{n}{2}-1} \omega_q \cdot g_{j \kappa} ^{\mathcal{E}_n} ( X_{q}), \quad \text{for} \ n=2,4.
\end{equation}
where the $\omega_q$ are given below, see \eqref{first_omega}--\eqref{last_omega}.
\end{proposition}

\begin{remark}
I suspect the validity of Proposition \ref{prop1_intro_go} extends to $n \in \N^*$ but I wasn't able to guess a general formula for the $\omega_q$.
\end{remark}

The formulas for the $\omega_q$ under the alignment conditions \eqref{o3bisounours} and \eqref{o33bisounours} are:  

\noindent $\bullet$ For $(n,p)=(1,0)$, $T_{\kappa}(X_{n+1}) = T_{\kappa}(X_p)$:

\begin{equation}
\label{first_omega}
\omega_0 = 2\cdot  \frac{\frac{m(X_1) U_{\kappa-1}(X_1)}{X_1}}{\frac{m(X_0) U_{\kappa-1}(X_0)}{X_0}  + \frac{m(X_2) U_{\kappa-1}(X_2)}{X_2}}.
\end{equation}

\noindent $\bullet$ For $(n,p)=(3,1)$, $T_{\kappa}(X_{n+1}) = T_{\kappa}(X_p)$:

\begin{equation}
\omega_0 = -2 \cdot \frac{\frac{m(X_2) m(X_3) U_{\kappa-1}(X_2) U_{\kappa-1}(X_3)}{X_2 X_3}}{\frac{m(X_0) m(X_1) U_{\kappa-1}(X_0) U_{\kappa-1}(X_1) }{X_0 X_1}  -  \frac{m(X_3) m(X_4) U_{\kappa-1}(X_3) U_{\kappa-1}(X_4)}{X_3 X_4}}
\end{equation}

\begin{equation}
\omega_1 = 2 \cdot \frac{\frac{m(X_0) m(X_2) U_{\kappa-1}(X_0) U_{\kappa-1}(X_2)}{X_0 X_2} }{\frac{m(X_0) m(X_1) U_{\kappa-1}(X_0) U_{\kappa-1}(X_1) }{X_0 X_1}  -  \frac{m(X_3) m(X_4) U_{\kappa-1}(X_3) U_{\kappa-1}(X_4)}{X_3 X_4}}.
\end{equation}


\noindent $\bullet$ For $(n,p)=(3,0)$, $T_{\kappa}(X_{n+1}) = T_{\kappa}(X_p)$:

\begin{equation}
\omega_0 = - 2 \cdot \frac{\frac{m(X_2) m(X_3) U_{\kappa-1}(X_2) U_{\kappa-1}(X_3) }{X_2 X_3}}{\frac{m(X_0) m(X_1) U_{\kappa-1}(X_0) U_{\kappa-1}(X_1) }{X_0 X_1} + \frac{m(X_1) m(X_4) U_{\kappa-1}(X_1) U_{\kappa-1}(X_4) }{X_1 X_4}}
\end{equation}

\begin{equation}
\omega_1 = 2 \cdot \frac{\frac{m(X_0) m(X_2) U_{\kappa-1}(X_0) U_{\kappa-1}(X_2) }{X_0 X_2} + \frac{m(X_2) m(X_4) U_{\kappa-1}(X_2) U_{\kappa-1}(X_4) }{X_2 X_4}}{\frac{m(X_0) m(X_1) U_{\kappa-1}(X_0) U_{\kappa-1}(X_1) }{X_0 X_1} + \frac{m(X_1) m(X_4) U_{\kappa-1}(X_1) U_{\kappa-1}(X_4) }{X_1 X_4}}.
\end{equation}

\noindent $\bullet$ For $(n,p)=(2,0)$, $T_{\kappa}(X_{n+1}) = T_{\kappa}(X_p)$:
\begin{equation}
\omega_0 =  \frac{\frac{m(X_2) U_{\kappa-1}(X_2)}{X_2}}{ \frac{m(X_0) U_{\kappa-1}(X_0)}{X_0}  + \frac{m(X_3) U_{\kappa-1}(X_3)}{X_3}}
\end{equation}


\noindent $\bullet$ For $(n,p)=(4,0)$, $T_{\kappa}(X_{n+1}) = T_{\kappa}(X_p)$:

\begin{equation}
\omega_0 =  - \frac{\frac{m(X_3) m(X_4) U_{\kappa-1}(X_3) U_{\kappa-1}(X_4) }{X_3 X_4}}{  \frac{m(X_0) m(X_1) U_{\kappa-1}(X_0) U_{\kappa-1}(X_1)}{X_0 X_1} + \frac{m(X_1) m(X_5) U_{\kappa-1}(X_1) U_{\kappa-1}(X_5)}{X_1 X_5}}
\end{equation}

\begin{equation}
\omega_1 =  \frac{\frac{m(X_0) m(X_3) U_{\kappa-1}(X_0) U_{\kappa-1}(X_3) }{X_0 X_3} + \frac{m(X_3) m(X_5) U_{\kappa-1}(X_3) U_{\kappa-1}(X_5) }{X_3 X_5}}{ \frac{m(X_0) m(X_1) U_{\kappa-1}(X_0) U_{\kappa-1}(X_1)}{X_0 X_1} + \frac{m(X_1) m(X_5) U_{\kappa-1}(X_1) U_{\kappa-1}(X_5)}{X_1 X_5}}.
\end{equation}

\noindent $\bullet$ For $(n,p)=(4,1)$, $T_{\kappa}(X_{n+1}) = T_{\kappa}(X_p)$:

\begin{equation}
\omega_0 =  - \frac{\frac{m(X_3) m(X_4) U_{\kappa-1}(X_3) U_{\kappa-1}(X_4) }{X_3 X_4}}{\frac{m(X_0) m(X_1) U_{\kappa-1}(X_0) U_{\kappa-1}(X_1)}{X_0 X_1} - \frac{m(X_4) m(X_5) U_{\kappa-1}(X_4) U_{\kappa-1}(X_5)}{X_4 X_5} }
\end{equation}

\begin{equation}
\label{last_omega}
\omega_1 = \frac{\frac{m(X_0) m(X_3) U_{\kappa-1}(X_0) U_{\kappa-1}(X_3) }{X_0 X_3}}{ \frac{m(X_0) m(X_1) U_{\kappa-1}(X_0) U_{\kappa-1}(X_1)}{X_0 X_1} - \frac{m(X_4) m(X_5) U_{\kappa-1}(X_4) U_{\kappa-1}(X_5)}{X_4 X_5} }.
\end{equation}

Let us explain how the formulas for the $\omega_q$'s were determined. For general $m \geq 1$, $E \in \boldsymbol{\Theta}_{m,\kappa} (D)$ in dimension 2 iff $\exists (X_q)_{q=0}^m \subset [-1,-|E|] \cup [|E|, 1]$, and $(\omega_q)_{q=0}^{m-1} \subset \R$, $\omega_q \leq 0$, such that 
\begin{equation*}
\label{special relationship}
g_{j \kappa} ^E (X_m) = \sum_{q=0} ^{m-1} \omega_q \cdot g_{j \kappa} ^E ( X_q), \quad \forall j \in \N^* \quad (\omega_q \ \text{independent of} \ j).
\end{equation*} 
If this linear relationship holds, it must be that for any choice of disctinct $j_1$, $j_2$, ..., $j_{m} \in \N^*$,
\begin{equation}
\begin{aligned}
\label{matrix_system_tosolve00}
\begin{pmatrix}
\omega_{m-1} \\
\omega_{m-2} \\
... \\
\omega_{0}
\end{pmatrix}
&= \begin{pmatrix} g_{\kappa} ^E (X_{m-1}) & g_{\kappa} ^E (X_{m-2}) & ... & g_{\kappa} ^E (X_0) \\ g_{2\kappa} ^E (X_{m-1}) & g_{2 \kappa} ^E (X_{m-2}) & ... &  g_{2 \kappa} ^E (X_0) \\
... & ... & ... & ... \\
g_{m \kappa} ^E (X_{m-1}) & g_{m\kappa} ^E (X_{m-2}) & ... &  g_{m \kappa} ^E (X_0) \\
\end{pmatrix} ^{-1}
   \begin{pmatrix} g_{ \kappa} ^E (X_{m}) \\ g_{2 \kappa} ^E (X_{m}) \\
   ... \\
  g_{m  \kappa} ^E (X_{m})  \end{pmatrix}  \\
  &=
\begin{pmatrix} g_{j_{1}  \kappa} ^E (X_{m-1)} & g_{j_{1}  \kappa} ^E (X_{m-2}) & ... & g_{j_{1}  \kappa} ^E (X_0) \\ g_{j_{2}  \kappa} ^E (X_{m-1}) & g_{j_{2}  \kappa} ^E (X_{m-2}) & ... &  g_{j_{2}  \kappa} ^E (X_0) \\
... & ... & ... & ... \\
g_{j_{m} \kappa} ^E (X_{m-1}) & g_{j_{m} \kappa} ^E (X_{m-2}) & ... &  g_{j_{m} \kappa} ^E (X_0) \\
\end{pmatrix} ^{-1}
   \begin{pmatrix} g_{j_{1} \kappa} ^E (X_m) \\ g_{j_{2} \kappa} ^E (X_m) \\
   ... \\
  g_{j_{m}  \kappa} ^E (X_m)  \end{pmatrix}.
\end{aligned}
\end{equation}
We therefore performed the above matrix multiplication (after computing the inverse matrix) for $m=0,1$ used $j_1,j_2 = 1,2$, the definition of $g_{j \kappa} ^E(x)$ given in \eqref{def:gE22} and applied the assumptions \eqref{conjecture_system_conj_intro}, \eqref{o1}, \eqref{o2} and \eqref{o3bisounours} (respectively, \eqref{conjecture_system_conj_intro_even}, \eqref{o11}, \eqref{o22}, \eqref{o33bisounours} and \eqref{o44}). A key identity that was used in those calculations and that is worth highlighting is:
\begin{equation}
\label{extra_id}
[ U_{2\kappa-1} (X_q), U_{\kappa-1} (X_{q'}) ] = U_{\kappa-1} (X_q) \cdot U_{\kappa-1} (X_{q'}) \cdot (T_{\kappa} (X_q) - T_{\kappa} (X_{q'})).
\end{equation}
For $m \geq 2$ it may be necessary to find a formula for $[ U_{j_1\kappa-1} (X_q), U_{j_2\kappa-1} (X_{q'}) ]$.

Unfortunately I was not able to formulate general assumptions like \eqref{ao3} or \eqref{ae3} adapted to this section. Looking at the formulas for $\omega_0$ and $\omega_1$ one can come up with various assumptions that could ensure that the $\omega_q$'s are negative. Since $m(X_q) \geq 0$, it is likely that the most convenient assumption is a statement about the signs of the ratios $U_{\kappa-1}(X_q) / X_q$, or equivalently $T'_{\kappa}(X_q) / X_q$. The following example illustrates the idea.
\begin{example} Let $n=1$, $\kappa \geq 4$, $\kappa \in \N_e$. Suppose that $\E_n \in \mathfrak{T}_{n,\kappa} ^+$. Then $T_{\kappa}(X_0) = T_{\kappa}(X_1) = T_{\kappa}(X_2)$. If $T'_{\kappa}(X_0) / X_0$, $T'_{\kappa}(X_2) / X_2 \geq 0$ and $T'_{\kappa}(X_1) / X_1 \leq 0$, then $\omega_0 =$ \eqref{first_omega} is negative and so $\E_n \in \boldsymbol{\Theta}_{1,\kappa} (D[d=2])$.

\end{example}
An inspection of the ratios $T'_{\kappa-1}(X_q) / X_q$ in the first 3 graphs in Figure \ref{fig:example_2nd_gap2} and all the graphs in Figure \ref{fig:example_2nd_gap22} reveal that the $\omega_q$'s given by \eqref{first_omega}--\eqref{last_omega} are strictly negative. Therefore the corresponding energies $\E_n (\kappa)$ are thresholds. In the last 3 graphs in Figure \ref{fig:example_2nd_gap2} we computed the $\omega_q$'s numerically using \eqref{matrix_system_tosolve00} for a handful of indices $j_1,j_2,j_3$. These appeared to be independent of $j_1,j_2,j_3$ and :

$$(\omega_0, \omega_1, \omega_2) \simeq (-0.31, -1.02, -1.71), \quad 4^{th} \ \text{graph in Figure} \ \ref{fig:example_2nd_gap2},$$
$$(\omega_0, \omega_1, \omega_2) \simeq (-0.31, -1.02, -1.72), \quad 5^{th} \ \text{graph in Figure} \ \ref{fig:example_2nd_gap2},$$
$$(\omega_0, \omega_1, \omega_2) \simeq (-0.30, -1.04, -1.72), \quad 6^{th} \ \text{graph in Figure} \ \ref{fig:example_2nd_gap2}.$$

This indicates that the corresponding energies $\E_n$ are also thresholds.

Finally, to construct these threshold energy solutions graphically, we note that the dynamical algorithms described in subsection \ref{subby} hold provided that the last step be changed to : the energy $E$ must be calibrated in such a way that the last point constructed, $X_{q_{n+1}}$, satisfies the alignment condition $T(X_{q_{n+1}}) = T(X_{r})$ where $X_r$ is one of the previously constructed points. 



\section{Description of the polynomial interpolation in dimension 2}
\label{desc_scheme}

This entire section is in dimension 2. In this section we adapt the linear system \eqref{interpol_intro} to the interval $J_2(\kappa) := (\cos^2(\pi / \kappa), \cos(\pi / \kappa))$. This will setup our framework behind Conjecture \ref{conjecture22}. In sections \ref{appli_scheme_k2} and \ref{appli_scheme_k3} we numerically implement the equations of this section.

Fix $\kappa \geq 4$, $\kappa \in \N_e$. First, let $\E_n, X_{0,n}, ..., X_{n+1,n}$ be the solutions of Propositions \ref{prop1} and \ref{prop2} (or equivalently Theorem \ref{thm_decreasing energy general}). Our aim is to find the coefficients $\rho_{j_q \kappa}$ of $\mathbb{A}(n) = \sum_{q=1} ^{N(n)} \rho_{j_{q} \kappa} (n) A_{j_q \kappa}$ so that a strict Mourre estimate holds on the interval $( \E_n , \E_{n-1} )$ -- which we refer to as the $n^{th}$ \textit{band}. 

For $n$ \underline{odd}, the linear system \eqref{interpol_intro} becomes (using notation \eqref{def:gE} and \eqref{def:GGE} instead) :

\begin{equation}
\label{system_interpol}
\begin{cases}
G_{\kappa} ^{\E_n} (X_{q,n}) = 0 & q = 0,...,(n-1)/2 \\
\frac{d}{dx} G_{\kappa} ^{\E_n} (X_{q,n}) = 0 & q = 1,...,(n-1)/2 \\
G_{\kappa} ^{\mathcal{E}_{n-1}} (X_{q,n-1}) = 0 & q = 0,...,(n-3)/2 \\
\frac{d}{dx} G_{\kappa} ^{\mathcal{E}_{n-1}} (X_{q,n-1}) = 0 & q = 1,...,(n-3)/2+1. \\
\end{cases}
\end{equation}
This system of $2n-1$ equations has at most rank $2n-1$, but part of our conjecture is that it always has rank $2n-1$.

For $n$ \underline{even}, the linear system \eqref{interpol_intro} becomes (using notation \eqref{def:gE} and \eqref{def:GGE} instead) :

\begin{equation}
\label{system_interpol_even}
\begin{cases}
G_{\kappa} ^{\E_n} (X_{q,n}) = 0 & q = 0,...,n/2-1 \\
\frac{d}{dx} G_{\kappa} ^{\E_n} (X_{q,n}) = 0 & q = 1,...,n/2 \\
G_{\kappa} ^{\mathcal{E}_{n-1}} (X_{q,n-1}) = 0 & q = 0,...,n/2-1 \\
\frac{d}{dx} G_{\kappa} ^{\mathcal{E}_{n-1}} (X_{q,n-1}) = 0 & q = 1,...,n/2-1. \\
\end{cases}
\end{equation}
Again, this system of $2n-1$ equations has at most rank $2n-1$, but part of our conjecture is that it always has rank $2n-1$.

For the coefficients $\rho_{j_q \kappa}$ we will assume $\Sigma = \{ \rho_{j_1 \kappa}, \rho_{j_2 \kappa}, ..., \rho_{j_{2n} \kappa} \}$ and further always take the convention that $j_1 = 1$ and $\rho_{j_1 \kappa} = \rho_{\kappa} = 1$. Thus we have a system of $2n-1$ unknowns and $2n-1$ equations.

\begin{remark}
Note that by writing out systems \eqref{system_interpol} and \eqref{system_interpol_even} in this way with the $X_{q,n}$'s coming from Propositions \ref{prop1} and \ref{prop2}, we are implicitly conjecturing that the sequence $\{ \E_n \}$ of Propositions \ref{prop1} and \ref{prop2} are the only thresholds in $J_2(\kappa)$, i.e.\ $\{ \E_n \} = J_2(\kappa) \cap \boldsymbol{\Theta}_{n, \kappa}(D)$.
\end{remark}

Let us justify the range of the index $q$ in the first two lines of \eqref{system_interpol}. Fix $n$ odd. By Lemma \ref{obvious_lem}, $G_{\kappa} ^{\mathcal{E}_n} (X_{1+q,n}) = G_{\kappa} ^{\mathcal{E}_n} (X_{n-q,n})$ for any choice of coefficients $\rho_{j \kappa}$ and $q=-1,0, ...,(n-1)/2$. Additionally, thanks to \eqref{problem1_tosolve_again_101_88}, $G_{\kappa} ^{\mathcal{E}_n} (X_{(n+1)/2,n}) = \sum_{q=0} ^{(n-1)/2} \omega_q \cdot G_{\kappa} ^{\mathcal{E}_n} (X_q,n)$. So to avoid obvious linear dependencies, we require the first line of system \eqref{system_interpol} only for $q=0,...,(n-1)/2$. As for the second line of system \eqref{system_interpol}, Lemma \ref{Lemma_not_trivial} entails (since the $X_{i,n}$'s are non-zero) that $\frac{d}{dx} G_{\kappa} ^{\mathcal{E}_n}(X_{1+q,n}) = 0 \Leftrightarrow \frac{d}{dx} G_{\kappa} ^{\mathcal{E}_n}(X_{n-q,n}) = 0$ for $q=-1,0,...,(n-1)/2$ and any choice of coefficients $\rho_{j \kappa}$. Additionally, thanks to Lemma \ref{symmetry_minima_E/2_mv}, $\frac{d}{dx} G_{\kappa} ^{\mathcal{E}_n} (X_{(n+1)/2,n}) =0$ always holds. So to avoid obvious linear dependencies, we only require the second line of system \eqref{system_interpol} only for $q=0,1,...,(n-1)/2$. Furthermore, we don't include $q=0$ but that is for a separate reason based only on numerical and graphical evidence -- and probably related to the fact that $X_{0,n} = \E_n$ does not belong to the interior of $S_{\E_n}$.

\begin{remark}
The justification for the range of the index $q$ in the third and fourth lines of \eqref{system_interpol} follows from the above discussion, by replacing $n$ with $n-1$. A similar reasoning as above establishes the range of the index $q$ for the case of $n$ even, i.e.\  system \eqref{system_interpol_even}.
\end{remark}

\section{Application of polynomial interpolation to the case $\kappa = 4$ in dimension 2}
\label{appli_scheme_k2}

In this section the results of the polynomial interpolation (described in section \ref{desc_scheme}) are displayed for $\kappa =4$, and $1 \leq n \leq 5$, in dimension 2.

Table \ref{sol_kappa_2_____} below gives inputs we need to feed linear systems \eqref{system_interpol} and \eqref{system_interpol_even} into the computer (we used numbers with much higher precision in the solver and to draw the graphs).

\begin{table}[H]
\footnotesize
  \begin{center}
    \begin{tabular}{c|c|c|c} 
    $n$ & Left endpoint & Right endpoint & $\Sigma = $     \\ [0.2em]
      \hline
1 & $\E_1 \simeq 0.6180 $, $X_1 \simeq 0.7861$  & $\E_0 = \cos(\pi/4) \simeq 0.7071$ & \{4,8\} \\[0.25em]
2 & $\E_2 \simeq 0.5773$, $X_1 \simeq 0.7071$, $X_2 \simeq 0.8164$ & $=$ left endpoint for $n=1$ & $\{4,8,12,28\}$  \\[0.25em]
3 & $\E_3 \simeq 0.5549$, $X_1 \simeq 0.6671$,  & $=$ left endpoint for $n=2$ & $\{4,8,12,16,20,24\}$ \\[0.15em]
& $X_2 \simeq 0.7449$, $X_3 \simeq 0.83187$ & & \\[0.15em]
4 & $\E_4 \simeq 0.5411$, $X_1 \simeq 0.6435$, $X_2 \simeq 0.7071$,  & $=$ left endpoint for $n=3$ & $\{4,8,12,16,20,24,28,64\}$ \\[0.15em]
& $X_3 \simeq 0.7653$, $X_4 \simeq 0.8408$ & & \\[0.25em]
5 & $\E_5 \simeq 0.5411$, $X_1 \simeq 0.6284$, $X_2 \simeq 0.6840$,  & $=$ left endpoint for $n=4$ & $\{4,8,12,16,20,24,28,32,36,40\}$ \\[0.15em]
& $X_3 \simeq 0.7294$, $X_4 \simeq 0.7778$, $X_5 \simeq 0.8466$ & & 
    \end{tabular}
  \end{center}
    \caption{Data to setup polynomial interpolation on $(\E_n, \E_{n-1})$. $\kappa=4$. $d=2$.}
        \label{sol_kappa_2_____}
\end{table}
\normalsize

Table \ref{sol_kappa_2____2_} below reports the solutions to the polynomial interpolation. Higher precision is obviously desirable and available, but we have chosen to report only a limited number of decimal points.
\begin{table}[H]
\footnotesize
  \begin{center}
    \begin{tabular}{c|c} 
    $n$ & Coefficients $\rho_{j \kappa}$     \\ [0.2em]
      \hline
1 & $[\rho_4, \rho_{8}] ^{T} = [1, \frac{17 + 8 \sqrt{5}}{62}  ]^{T} \simeq [1,0.56271]^{T}$ \\[0.25em]
2 & $[\rho_4, \rho_8, \rho_{12},\rho_{28}]^{T} \simeq [1, 0.79123, 0.19359, 0.02771]^T$ \\[0.25em]
3 & $[\rho_4, \rho_8, \rho_{12}, \rho_{16}, \rho_{20}, \rho_{24} ]^{T} \simeq [1, 1.328058, 0.98129, 0.45526, 0.12626, 0.01635]^{T}$ \\[0.25em]
4 & $[\rho_4, \rho_8, \rho_{12}, \rho_{16}, \rho_{20}, \rho_{24}, \rho_{28}, \rho_{64} ]^{T} \simeq$ \\[0.1em]
& $ [1, 1.44284, 1.23820, 0.72957, 0.29600,  0.07629, 0.00959, -0.0000079]^{T}$  \\[0.25em]
5 & $[\rho_4, \rho_8, \rho_{12}, \rho_{16}, \rho_{20}, \rho_{24}, \rho_{28}, \rho_{32}, \rho_{36}, \rho_{40}]^{T} \simeq$ \\[0.1em]
& $[1, 1.55110, 1.51418, 1.09196, 0.60426, 0.25666, 0.08171, 0.01849, 0.00266, 0.00018]^{T}$
\end{tabular}
  \end{center}
    \caption{Results of the linear interpolation for $\kappa =4$}
        \label{sol_kappa_2____2_}
\end{table}
\normalsize
We noticed that only for very small values of $\kappa$ and $n$ is Python able to produce exact solutions. For example, for $\kappa=4$, $n=2$ we have a rare ocurrence of an exact solution (exact values of $X_0$, $X_1$, $X_2$ and $X_3$ were used) :
\begin{equation*}
\label{solution_3rd_band_D,2d simpler comb}
\begin{cases} 
\rho_4 = 1 \\
\rho_{8} = \frac{122211686455285242171392}{145241499888481405699961} - \frac{3260892257389038020608 \sqrt{5}}{145241499888481405699961} \simeq 0.79123 & \\
\rho_{12} = \frac{3077650500266123893602}{13203772717134673245451} - \frac{233195044949152756416\sqrt{5}}{13203772717134673245451} \simeq  0.19359 & \\
\rho_{28} = \frac{2012748316225734218931}{145241499888481405699961} + \frac{900129194439352172352\sqrt{5}}{145241499888481405699961} \simeq 0.02771. & \\
\end{cases}
\end{equation*}

Figures \ref{fig:test_k4n1} -- \ref{fig:test_k4n555} below depict $G_{\kappa} ^E(x)$, the functional representation of the commutator $[D, \i \mathbb{A}]_{\circ}$ localized at energy $E$. When looking at the Figures we want to see that $G_{\kappa} ^E(x)$ is $\geq 0$ when $E$ is a threshold energy, i.e.\ $E \in \{ \E_n\}$, but strictly positive for $x\in [E,1]$ and $E \in (\E_n, \E_{n-1})$. We choose only one value $E \in (\E_n, \E_{n-1})$ for illustrative purposes. Furthermore, thanks to Lemma \ref{even_pol} it is enough to observe $G_{\kappa}^E(x)$ for $x$ positive.

\begin{figure}[H]
  \centering
 \includegraphics[scale=0.322]{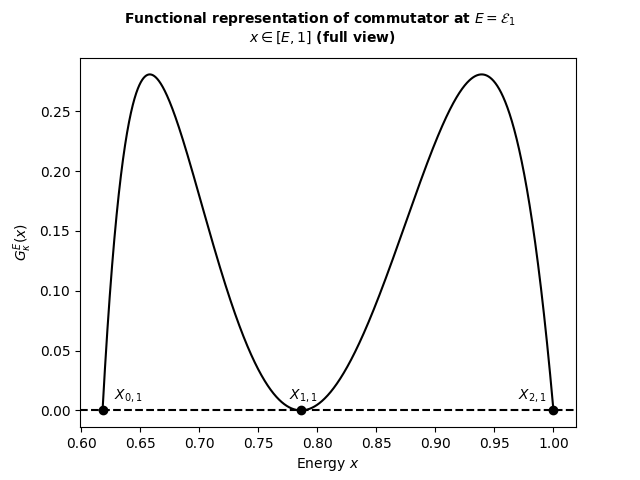}
  \includegraphics[scale=0.322]{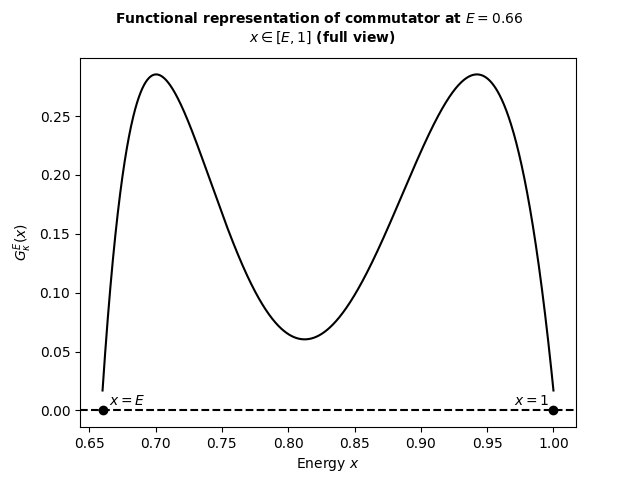}
  \includegraphics[scale=0.322]{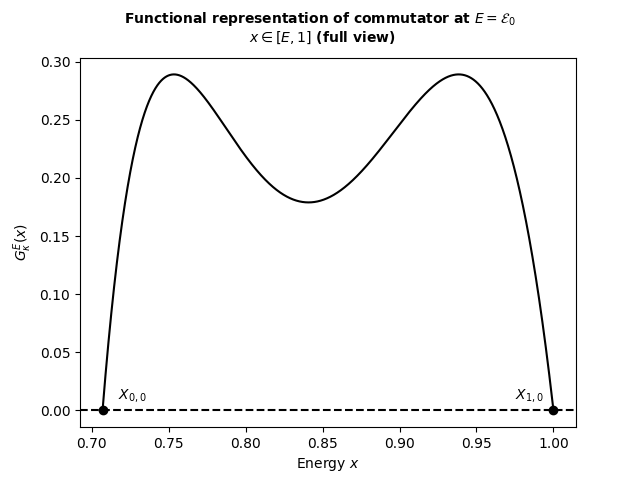}
\caption{$G_{\kappa=4} ^E(x)$, $x \in [E,1]$. Graphs at: $E=\E_1$, $E=0.66$, $E=\E_0$. $1^{st}$ band}
\label{fig:test_k4n1}
\end{figure}

\begin{figure}[H]
  \centering
 \includegraphics[scale=0.322]{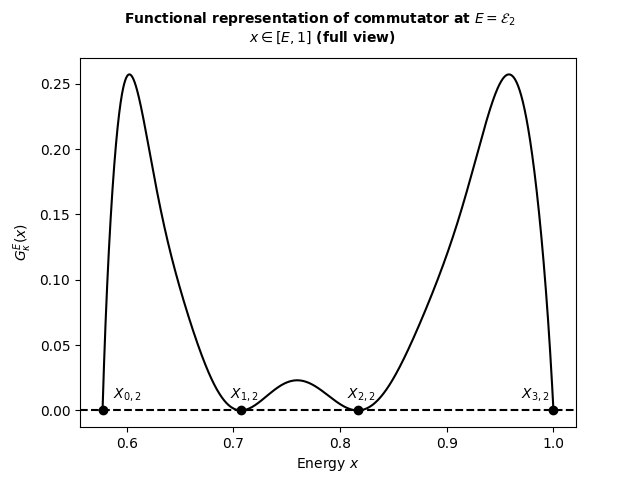}
  \includegraphics[scale=0.322]{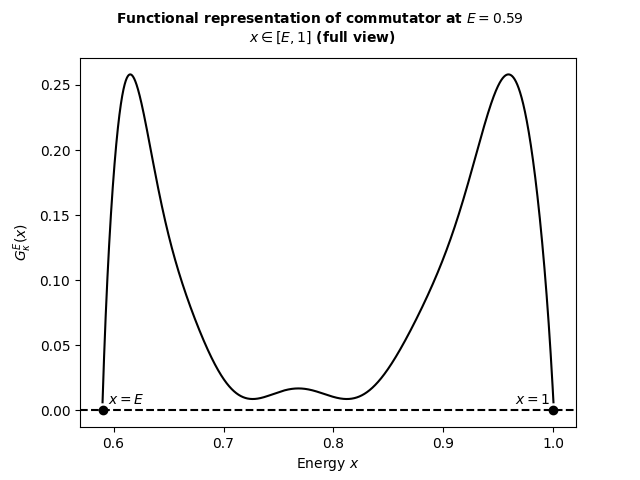}
  \includegraphics[scale=0.322]{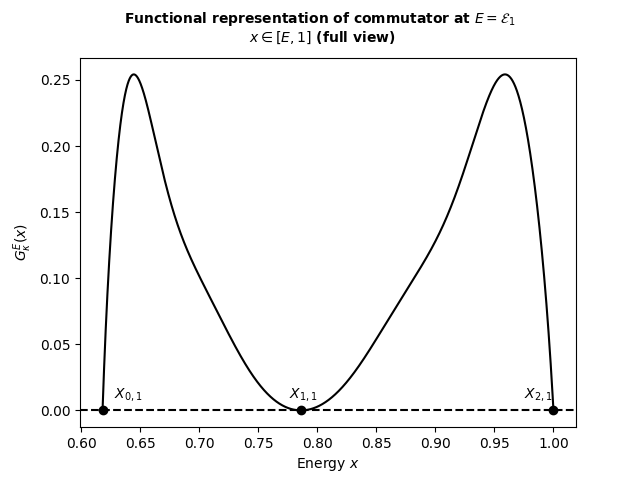}
\caption{$G_{\kappa=4} ^E(x)$, $x \in [E,1]$. Graphs at: $E=\E_2$, $E=0.59$, $E=\E_1$. $2^{nd}$ band}
\label{fig:test_k4n2}
\end{figure}

\begin{figure}[H]
  \centering
 \includegraphics[scale=0.24]{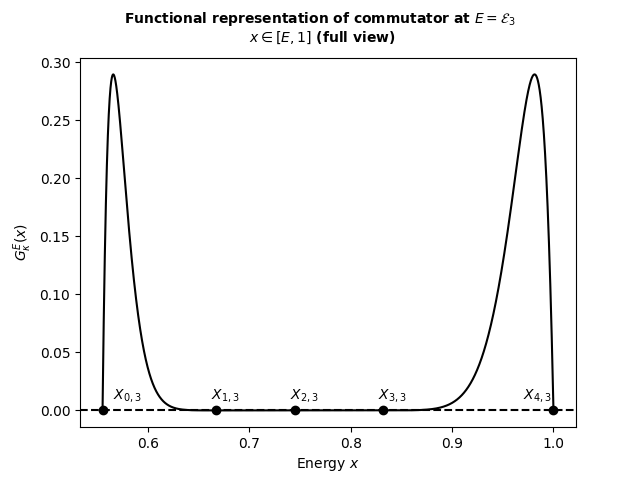}
  \includegraphics[scale=0.24]{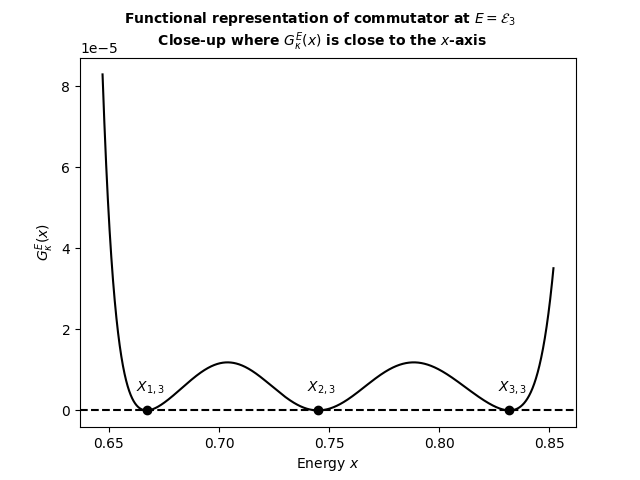}
  \includegraphics[scale=0.24]{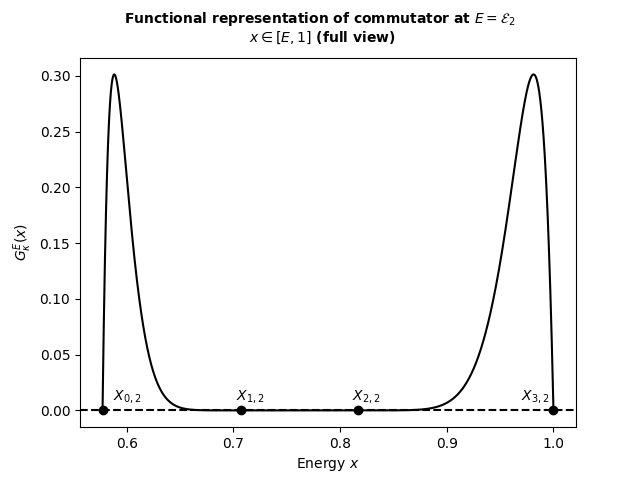}
   \includegraphics[scale=0.24]{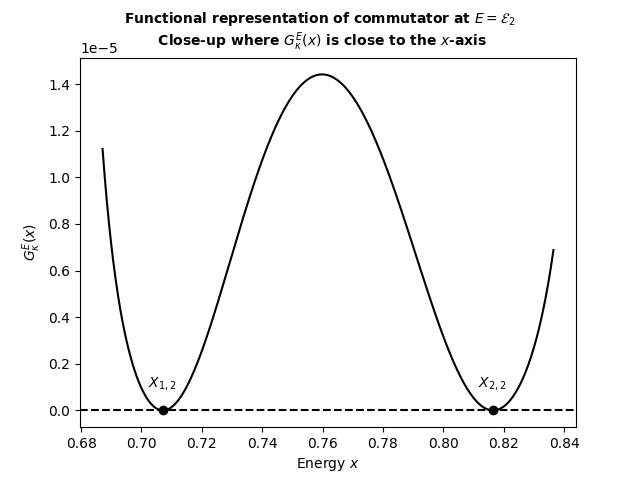}
     \includegraphics[scale=0.322]{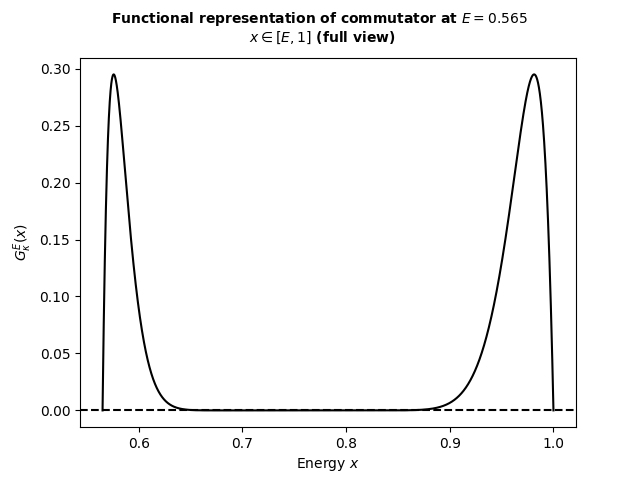}
  \includegraphics[scale=0.322]{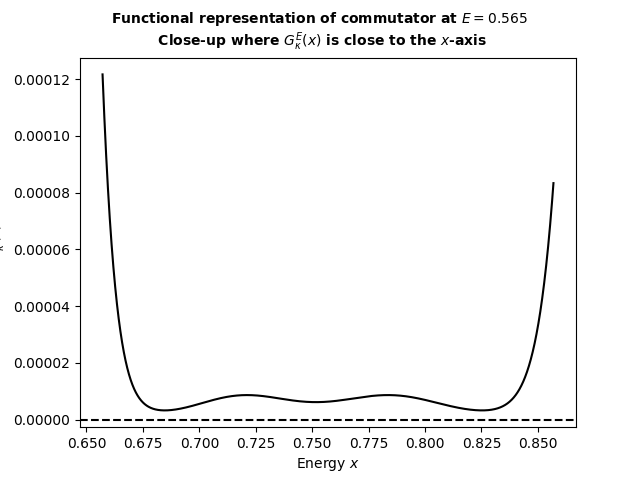}
   \includegraphics[scale=0.27]{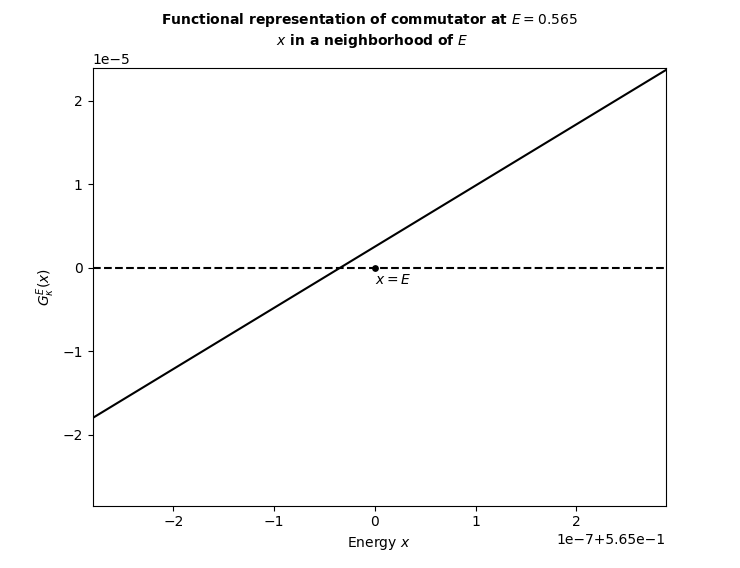}
\caption{$G_{\kappa=4} ^E(x)$, $x \in [E,1]$. Graphs at: $E=\E_3$, $E=0.565$, $E=\E_2$. $3^{rd}$ band}
\label{fig:test_k4n3}
\end{figure}

\begin{figure}[H]
  \centering
  \includegraphics[scale=0.24]{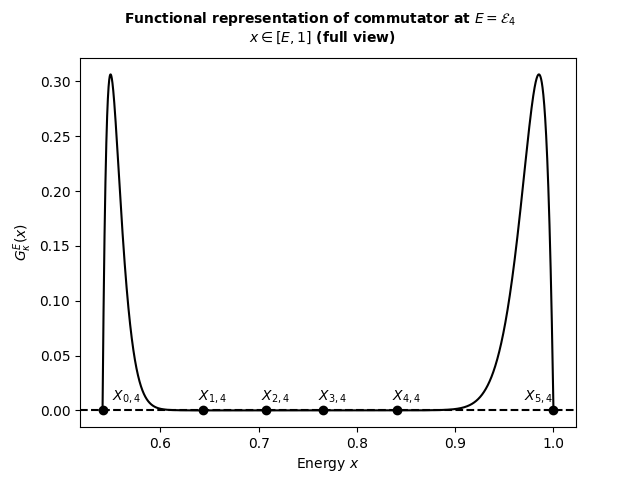}
   \includegraphics[scale=0.24]{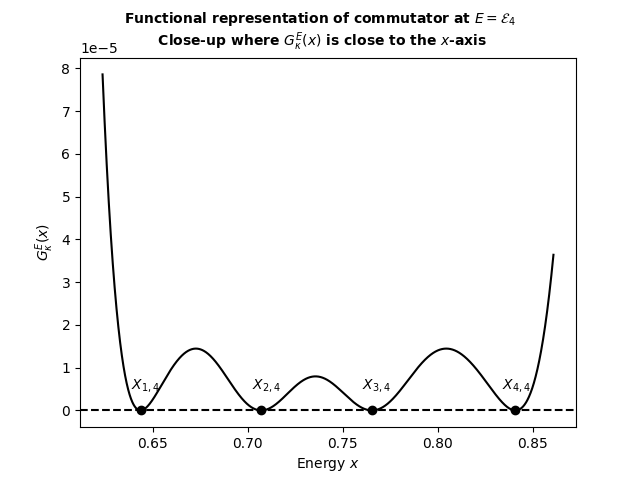}
    \includegraphics[scale=0.24]{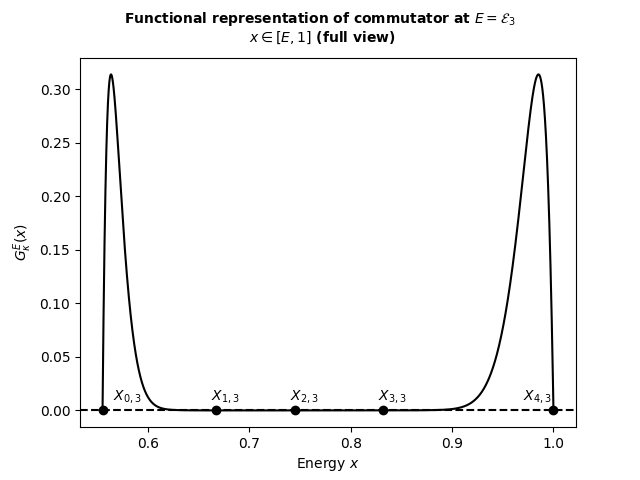}
  \includegraphics[scale=0.24]{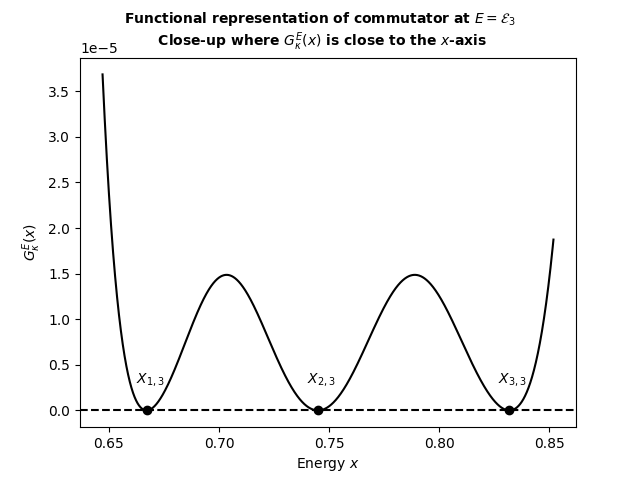}
     \includegraphics[scale=0.318]{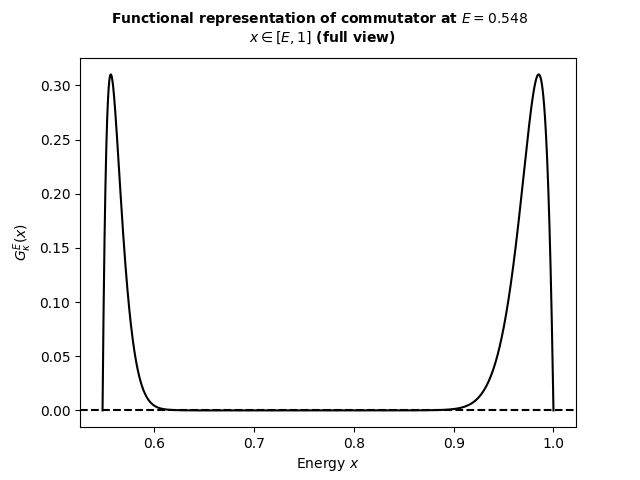}
  \includegraphics[scale=0.318]{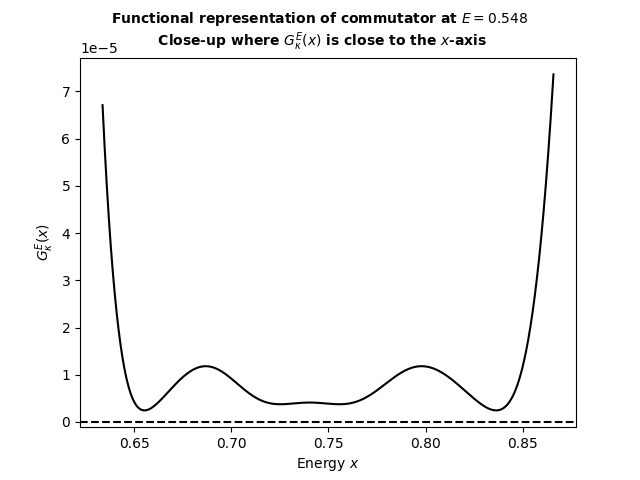}
   \includegraphics[scale=0.262]{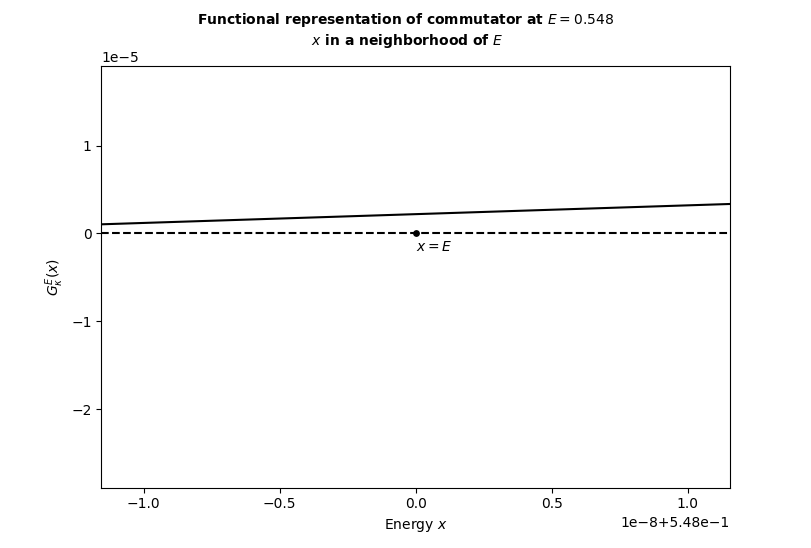}
\caption{$G_{\kappa=4} ^E(x)$, $x \in [E,1]$. Graphs at: $E=\E_4$, $E=0.548$, $E=\E_3$. $4^{th}$ band}
\label{fig:test_k4n4}
\end{figure}

\begin{figure}[H]
  \centering
  \includegraphics[scale=0.24]{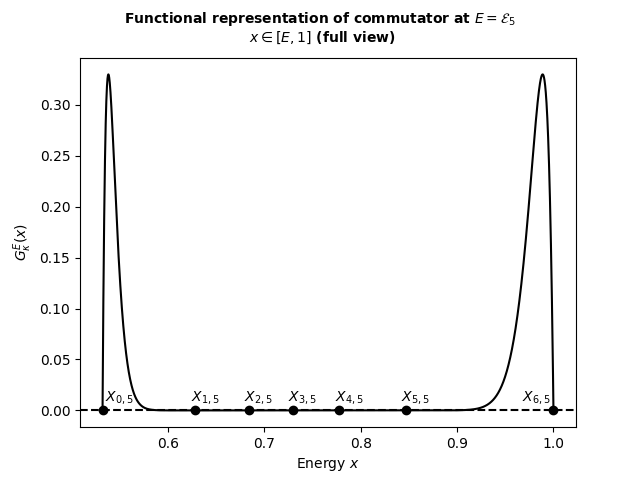}
   \includegraphics[scale=0.24]{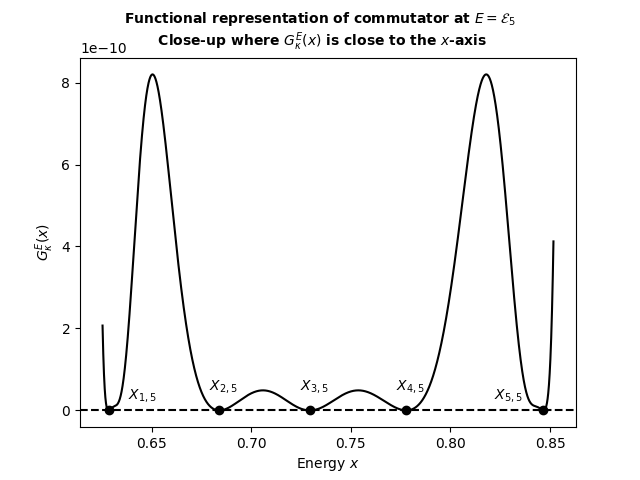}
    \includegraphics[scale=0.24]{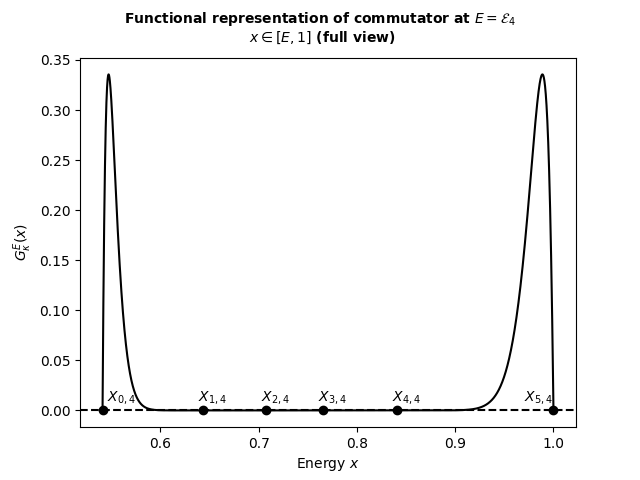}
  \includegraphics[scale=0.24]{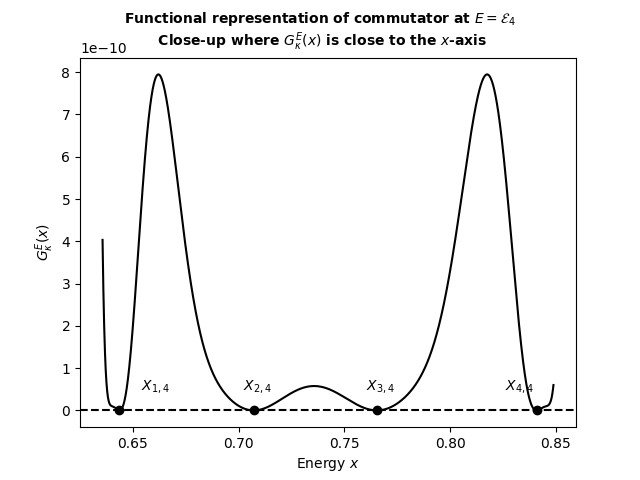}
     \includegraphics[scale=0.322]{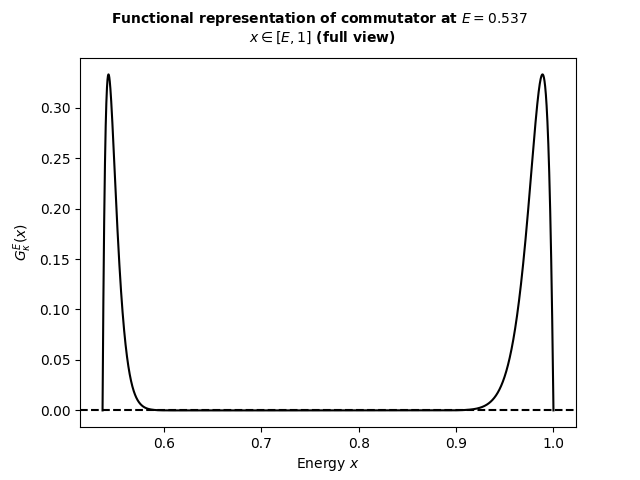}
  \includegraphics[scale=0.322]{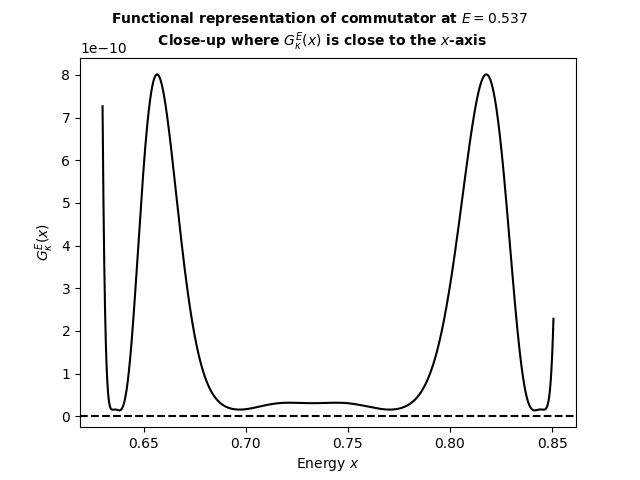}
   \includegraphics[scale=0.322]{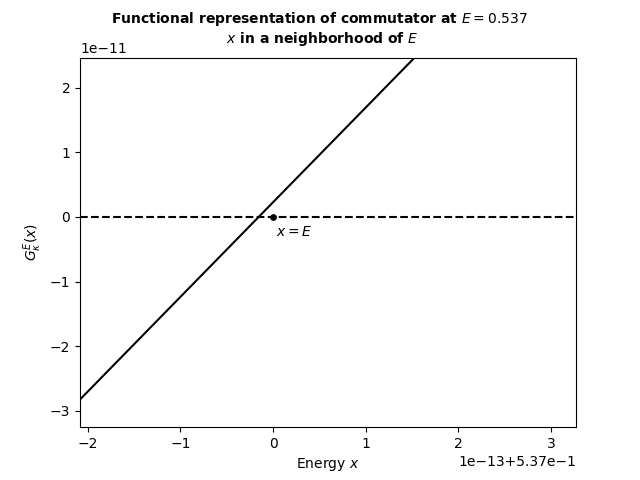}
\caption{$G_{\kappa=4} ^E(x)$, $x \in [E,1]$. Graphs at: $E=\E_5$, $E=0.537$, $E=\E_4$. $5^{th}$ band}
\label{fig:test_k4n555}
\end{figure}

\section{Application of polynomial interpolation to the case $\kappa = 6$ in dimension 2}
\label{appli_scheme_k3}

In this section the results of the polynomial interpolation (described in section \ref{desc_scheme}) are displayed for $\kappa =6$, and $1 \leq n \leq 5$, in dimension 2.

Table \ref{sol_kappa_3_____} gives inputs we need to feed linear systems \eqref{system_interpol} and \eqref{system_interpol_even} into the computer (we used numbers with much higher precision to draw the graphs).

\begin{table}[H]
\footnotesize
  \begin{center}
    \begin{tabular}{c|c|c|c} 
    $n$ & Left endpoint & Right endpoint & $\Sigma = $     \\ [0.2em]
      \hline
1 & $\E_1 \simeq 0.8201 $, $X_1 \simeq 0.9056$  & $\E_0 = \cos(\pi/6) \simeq 0.8660$ & \{6,12\} \\[0.25em]
2 & $\E_2 \simeq 0.7977$, $X_1 \simeq 0.8660$, $X_2 \simeq 0.9211$ & $=$ left endpoint for $n=1$ & $\{6,12,18,42\}$  \\[0.25em]
3 & $\E_3 \simeq 0.7848$, $X_1 \simeq 0.8447$,  & $=$ left endpoint for $n=2$ & $\{6, 12, 18, 24, 30, 36\}$ \\[0.15em]
& $X_2 \simeq 0.8859$, $X_3 \simeq 0.9291$ & & \\[0.15em]

4 & $\E_4 \simeq 0.7766$, $X_1 \simeq 0.8316$, $X_2 \simeq 0.8660$,  & $=$ left endpoint for $n=3$ & $\{6, 12, 18, 24, 30, 36, 42, 90\}$ \\[0.15em]
& $X_3 \simeq 0.8968$, $X_4 \simeq 0.9339$ & & \\[0.25em]

5 & $\E_5 \simeq 0.7711$, $X_1 \simeq 0.8229$, $X_2 \simeq 0.8534$,  & $=$ left endpoint for $n=4$ & $\{6, 12, 18, 24, 30, 36, 42, 48, 54, 60\}$ \\[0.15em]
& $X_3 \simeq 0.8781$, $X_4 \simeq 0.9035$, $X_5 \simeq 0.9370$ & & 
    \end{tabular}
  \end{center}
    \caption{Data to setup polynomial interpolation on $(\E_n, \E_{n-1})$. $\kappa=6$. $d=2$.}
        \label{sol_kappa_3_____}
\end{table}
\normalsize

Table \ref{sol_kappa_3____3_} below reports the solutions to the polynomial interpolation. Higher precision is obviously desirable and available, but we have chosen to report only a limited number of decimal points.

\begin{table}[H]
\footnotesize
  \begin{center}
    \begin{tabular}{c|c} 
    $n$ & Coefficients $\rho_{j \kappa}$     \\ [0.2em]
      \hline
1 & $[\rho_6, \rho_{12}] ^{T} \simeq [1, 0.57405]^{T}$ \\[0.25em]
2 & $[\rho_6, \rho_{12}, \rho_{18},\rho_{42}]^{T} \simeq [1, 0.87489, 0.26445, 0.01568]^T$ \\[0.25em]
3 & $[\rho_6, \rho_{12}, \rho_{18}, \rho_{24}, \rho_{30}, \rho_{36} ]^{T} \simeq [1, 1.34434, 1.01470, 0.48606, 0.14113, 0.01952]^{T}$ \\[0.25em]
4 & $[\rho_6, \rho_{12}, \rho_{18}, \rho_{24}, \rho_{30}, \rho_{36}, \rho_{42}, \rho_{90} ]^{T} \simeq$ \\[0.1em]
& $ [1, 1.43702, 1.22364, 0.71194, 0.28365, 0.07146, 0.00878, 4.05e-6]^{T}$  \\[0.25em]
5 & $[\rho_6, \rho_{12}, \rho_{18}, \rho_{24}, \rho_{30}, \rho_{36}, \rho_{42}, \rho_{48}, \rho_{54}, \rho_{60}]^{T} \simeq$ \\[0.1em]
& $[1, 1.56415, 1.54850, 1.139066, 0.64678, 0.28365, 0.09384, 0.02222, 0.00337, 0.00024]^{T}$
\end{tabular}
  \end{center}
    \caption{Results of the linear interpolation  for $\kappa =6$}
        \label{sol_kappa_3____3_}
\end{table}

Figures \ref{fig:test_k6n1} -- \ref{fig:test_k6n5} below depict $G_{\kappa} ^E(x)$, the functional representation of the commutator $[D, \i \mathbb{A}]_{\circ}$ localized at energy $E$. When looking at the Figures we want to see that $G_{\kappa} ^E(x)$ is $\geq 0$ when $E$ is a threshold energy, i.e.\ $E \in \{ \E_n\}$, but strictly positive for $x\in [E,1]$ and $E \in (\E_n, \E_{n-1})$. We choose only one value $E \in (\E_n, \E_{n-1})$ for illustrative purposes. Furthermore, thanks to Lemma \ref{even_pol} it is enough to observe $G_{\kappa}^E(x)$ for $x$ positive.

\begin{figure}[H]
  \centering
 \includegraphics[scale=0.322]{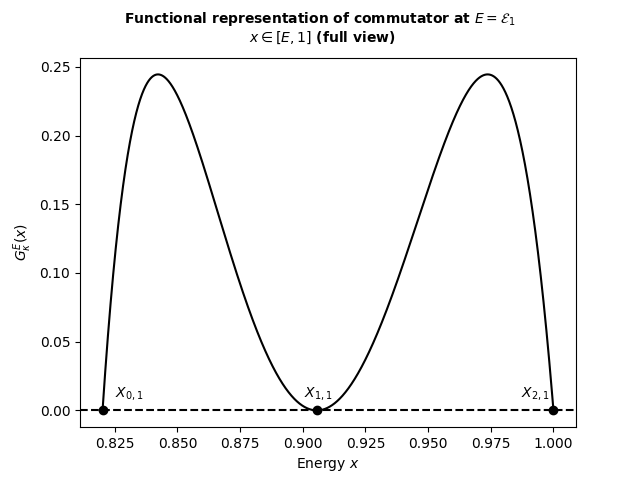}
  \includegraphics[scale=0.322]{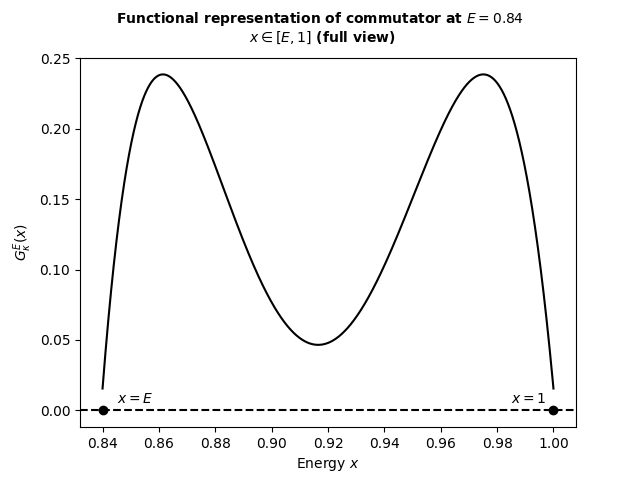}
  \includegraphics[scale=0.322]{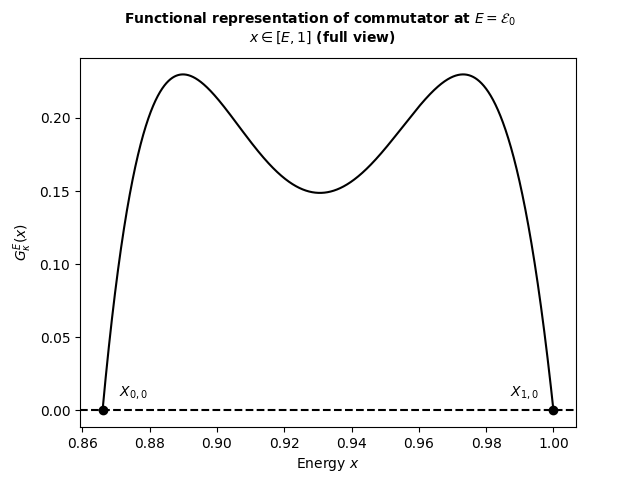}
\caption{$G_{\kappa=6} ^E(x)$, $x \in [E,1]$. Graphs at: $E=\E_1$, $E=0.84$, $E=\E_0$. $1^{st}$ band}
\label{fig:test_k6n1}
\end{figure}

\begin{figure}[H]
  \centering
 \includegraphics[scale=0.322]{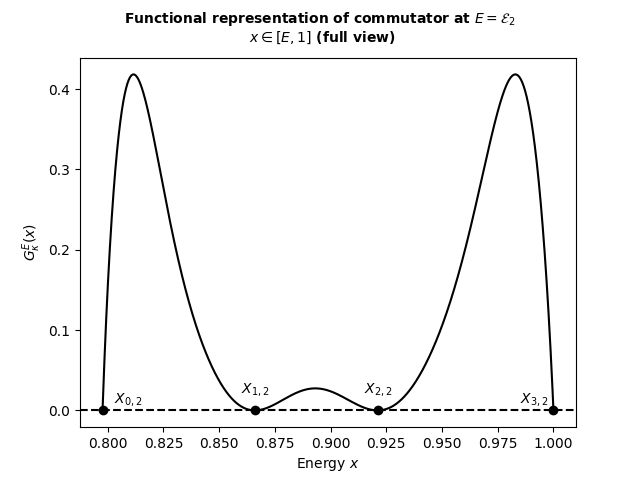}
  \includegraphics[scale=0.322]{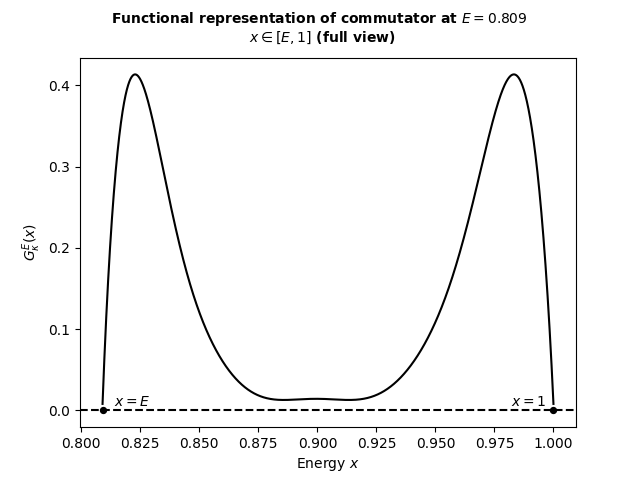}
  \includegraphics[scale=0.322]{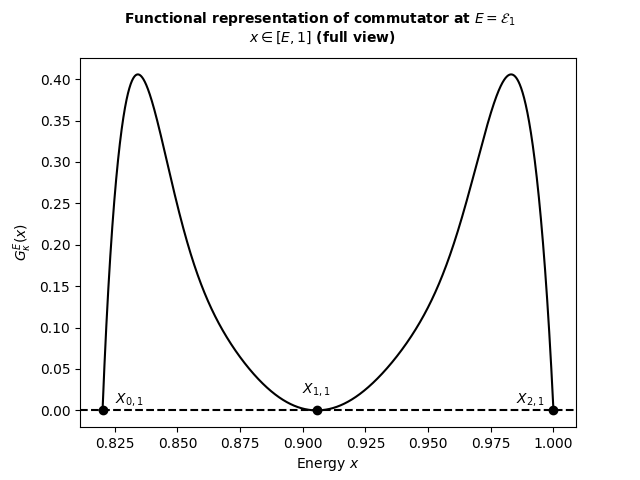}
\caption{$G_{\kappa=6} ^E(x)$, $x \in [E,1]$. Graphs at: $E=\E_2$, $E=0.809$, $E=\E_1$. $2^{nd}$ band}
\label{fig:test_k6n2}
\end{figure}

\begin{figure}[H]
  \centering
 \includegraphics[scale=0.24]{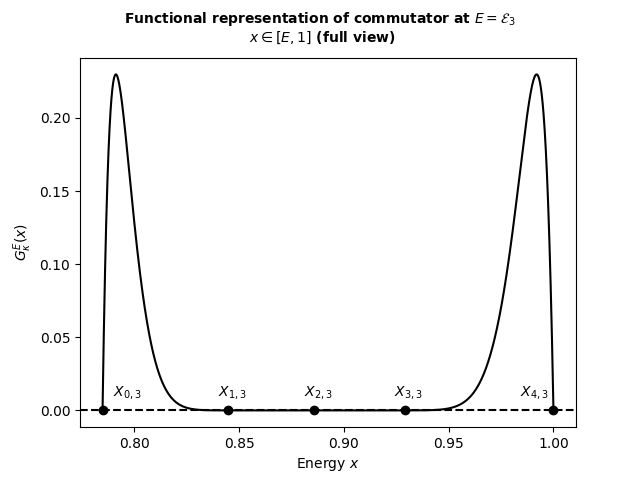}
  \includegraphics[scale=0.24]{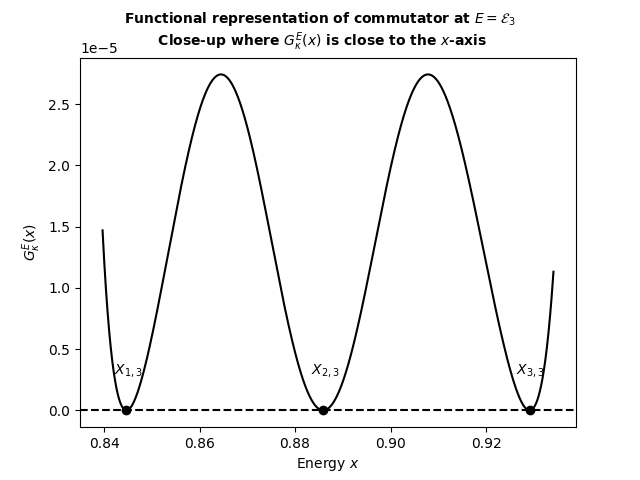}
  \includegraphics[scale=0.24]{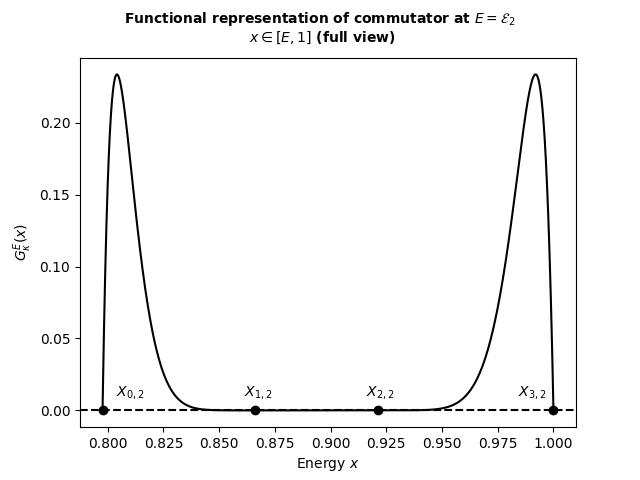}
   \includegraphics[scale=0.24]{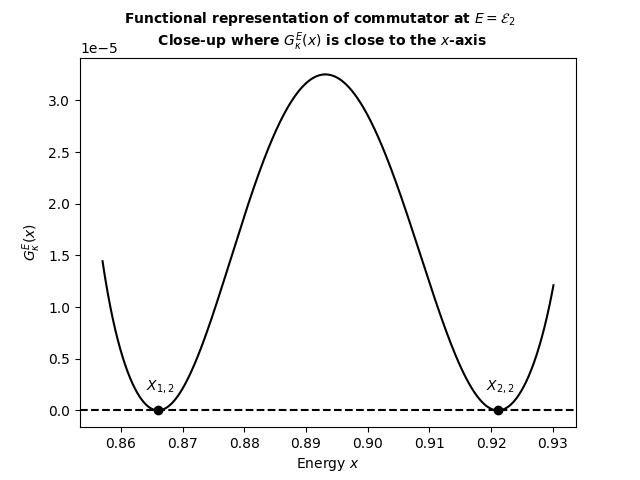}
     \includegraphics[scale=0.31]{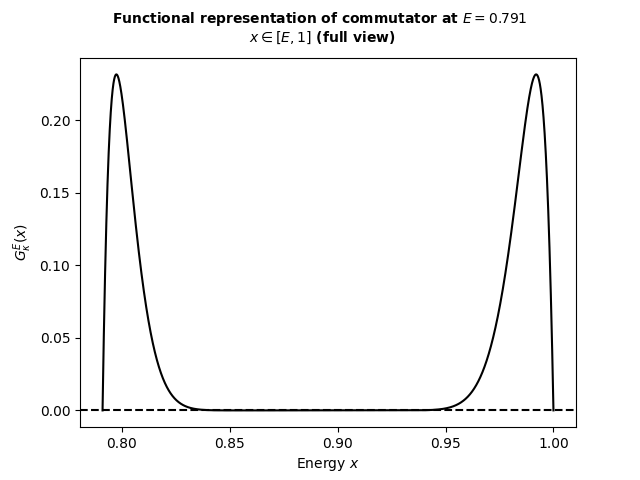}
  \includegraphics[scale=0.31]{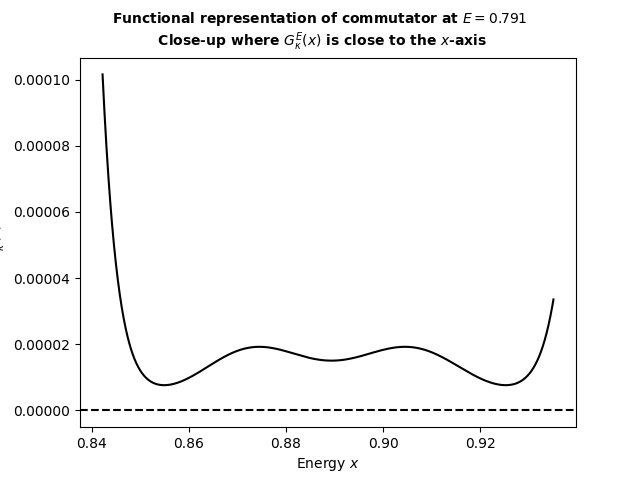}
   \includegraphics[scale=0.266]{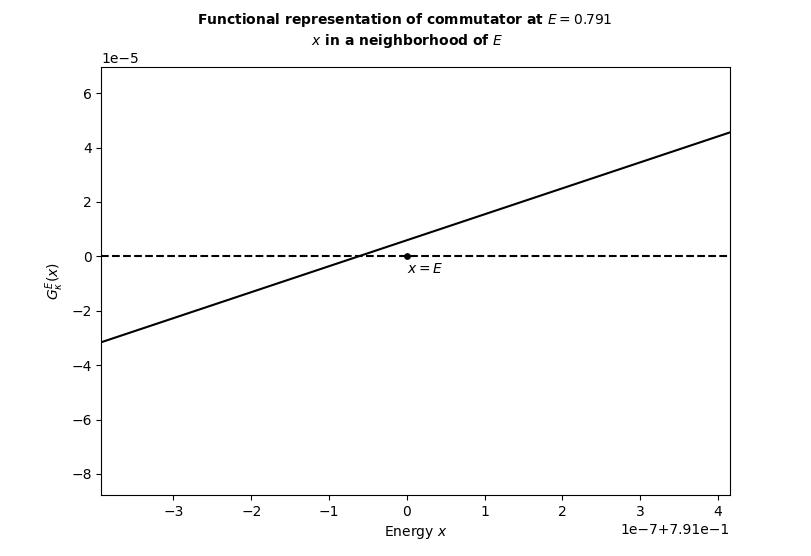}
\caption{$G_{\kappa=6} ^E(x)$, $x \in [E,1]$. Graphs at: $E=\E_3$, $E=0.791$, $E=\E_2$. $3^{rd}$ band}
\label{fig:test_k6n3}
\end{figure}

\begin{figure}[H]
  \centering
 \includegraphics[scale=0.24]{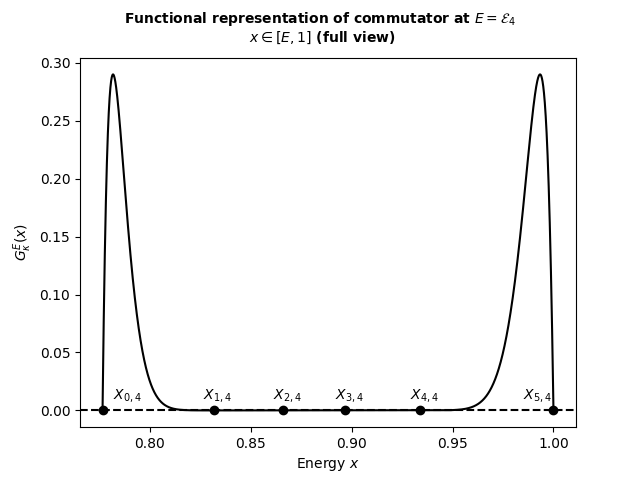}
  \includegraphics[scale=0.24]{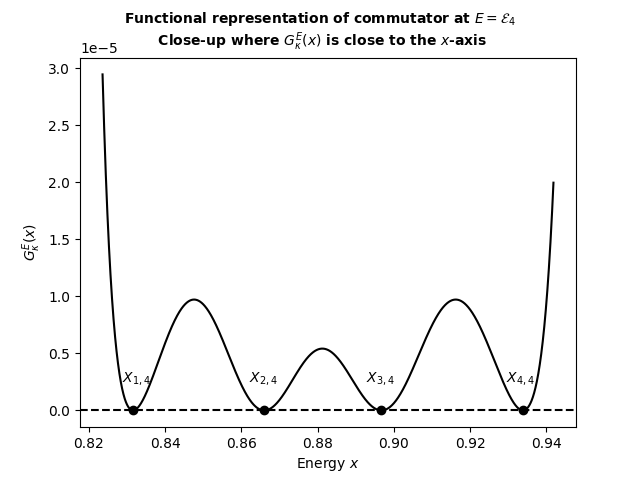}
  \includegraphics[scale=0.24]{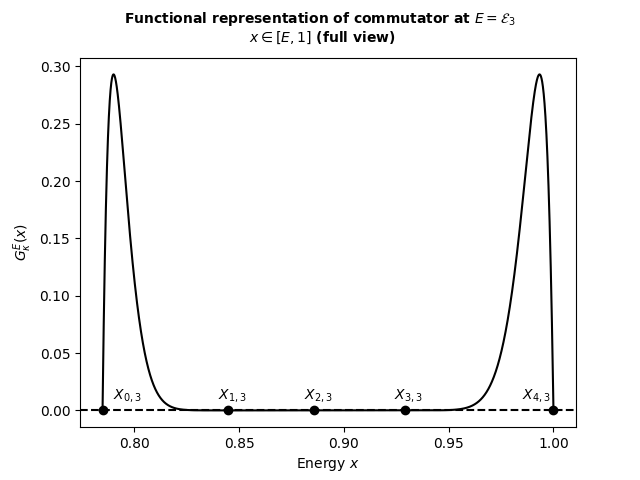}
   \includegraphics[scale=0.24]{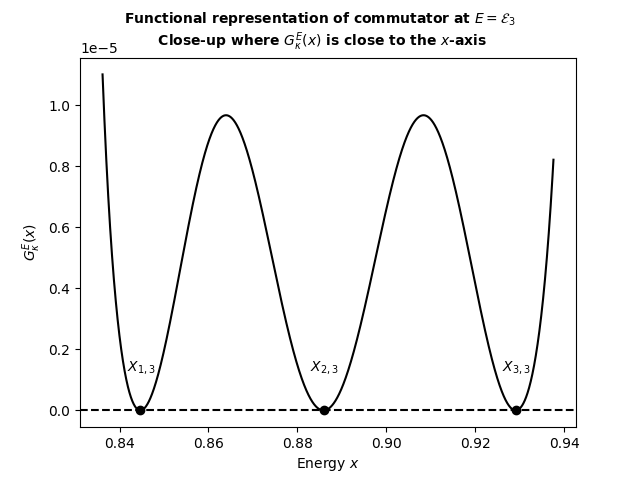}
     \includegraphics[scale=0.308]{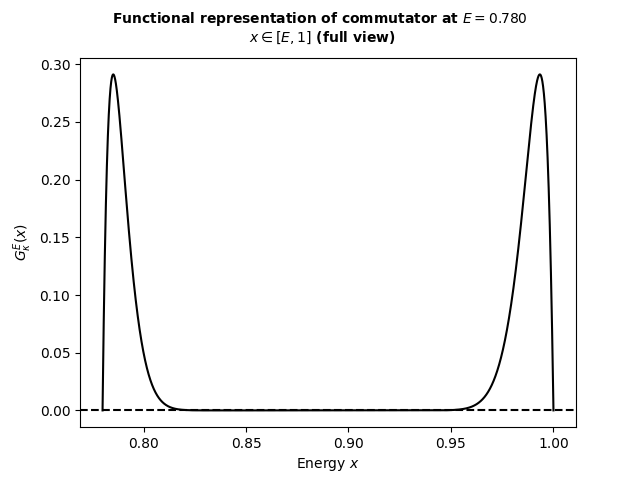}
  \includegraphics[scale=0.308]{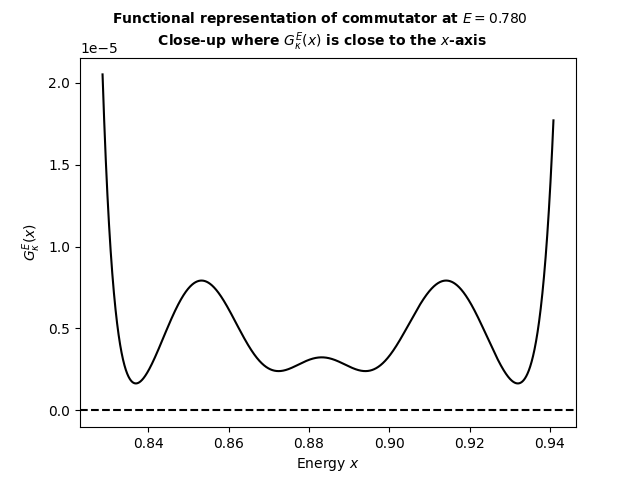}
   \includegraphics[scale=0.272]{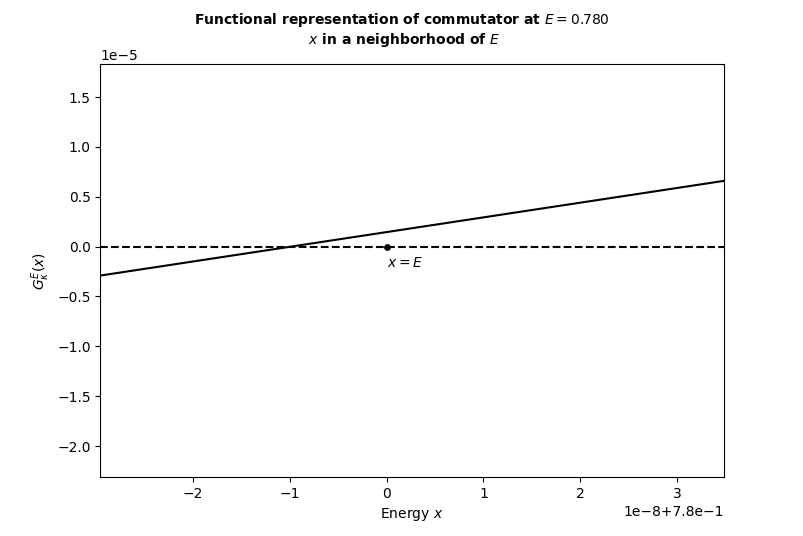}
\caption{$G_{\kappa=6} ^E(x)$, $x \in [E,1]$. Graphs at: $E=\E_4$, $E=0.780$, $E=\E_3$. $4^{th}$ band}
\label{fig:test_k6n4}
\end{figure}

\begin{figure}[H]
  \centering
 \includegraphics[scale=0.24]{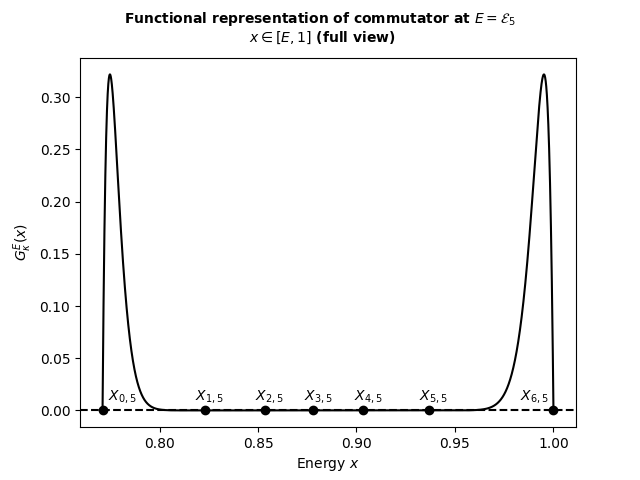}
  \includegraphics[scale=0.24]{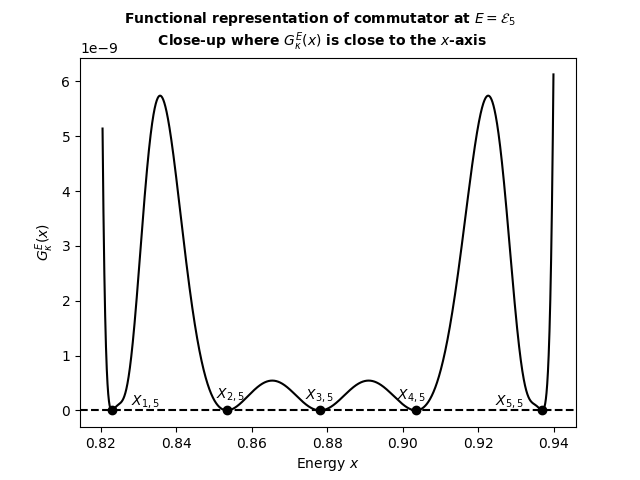}
  \includegraphics[scale=0.24]{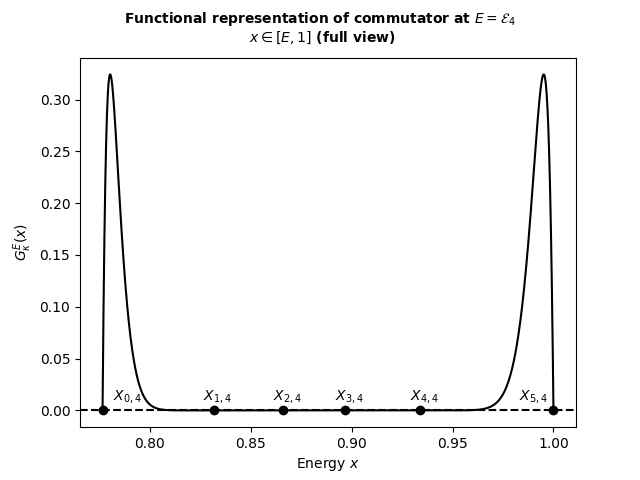}
   \includegraphics[scale=0.24]{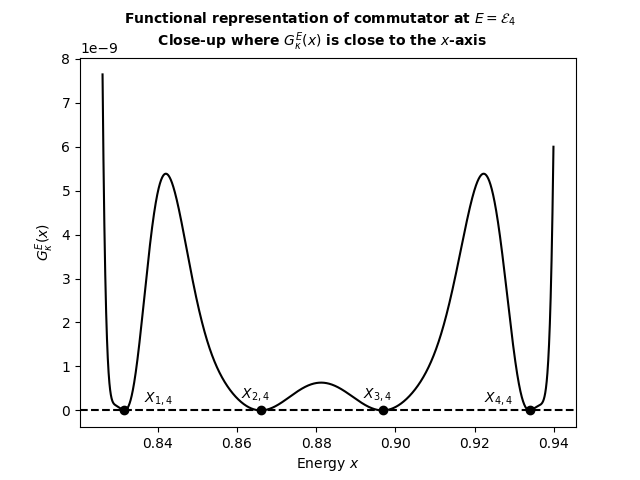}
     \includegraphics[scale=0.32]{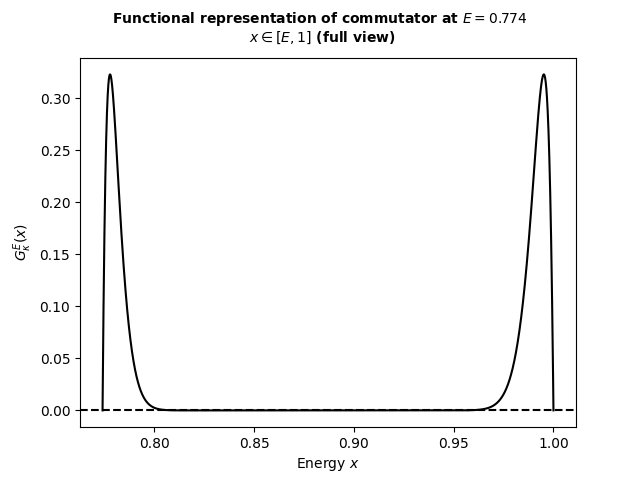}
  \includegraphics[scale=0.32]{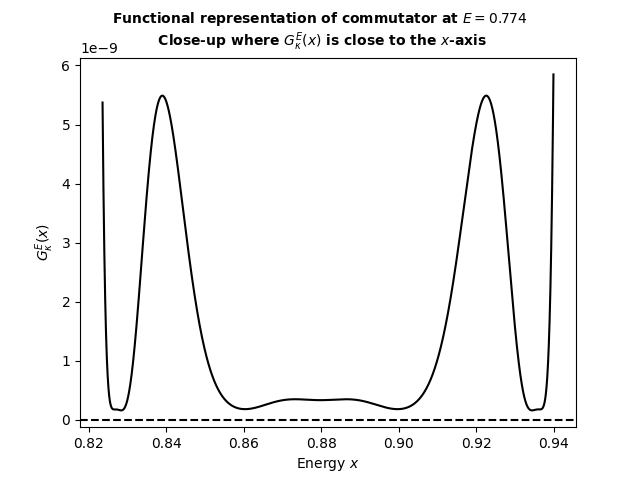}
   \includegraphics[scale=0.32]{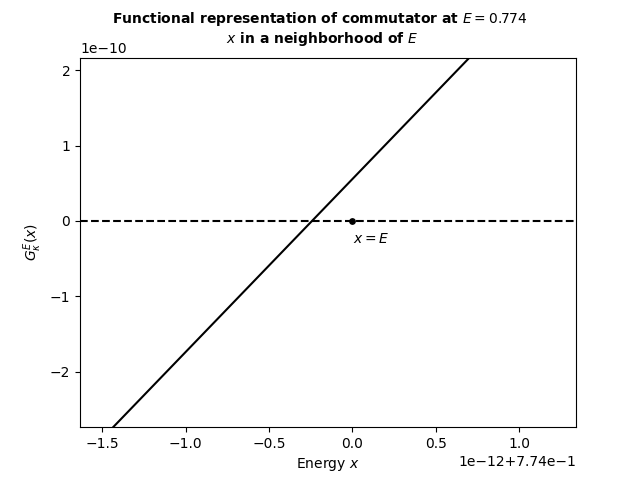}
\caption{$G_{\kappa=6} ^E(x)$, $x \in [E,1]$. Graphs at: $E=\E_5$, $E=0.774$, $E=\E_4$. $5^{th}$ band}
\label{fig:test_k6n5}
\end{figure}

\section{How to choose the correct indices $\Sigma =  \{ j_1 \kappa , j_2 \kappa , ..., j_{2n}\kappa \}$ ?}
\label{sectionCHOICE}
This section is in dimension 2. We discuss the possibility of using other index sets $\Sigma$ when performing the linear interpolation.

For $(\kappa,n) = (4,2)$, we checked that the indices $\Sigma = [4,8,12,4l]$ are equally valid for $l=7,8,9$ but not valid for $l=4,5,6$. $\Sigma = [4,8,16,4l]$ are valid for $l = 7,8,9$ but not valid for $l=5,6$. The point is that there is generally not only one valid index set $\Sigma$. 

\begin{figure}[H]
  \centering
 \includegraphics[scale=0.322]{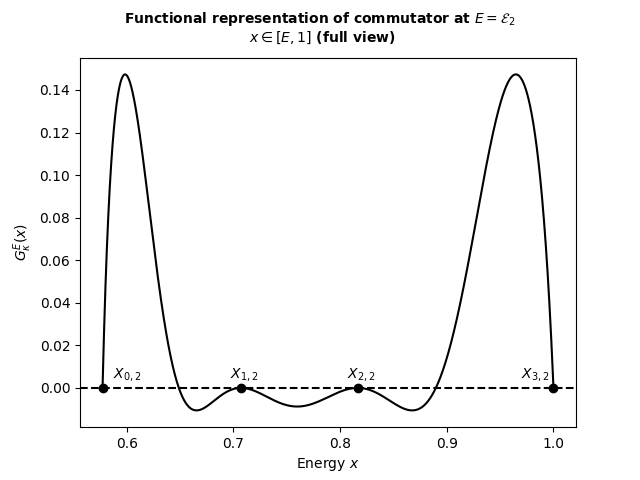}
  \includegraphics[scale=0.322]{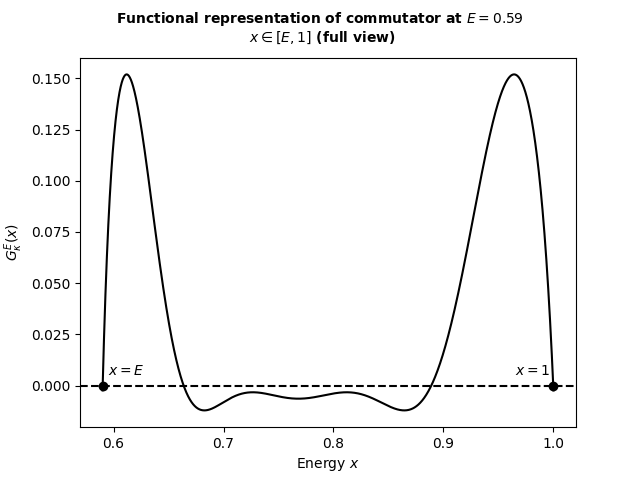}
  \includegraphics[scale=0.322]{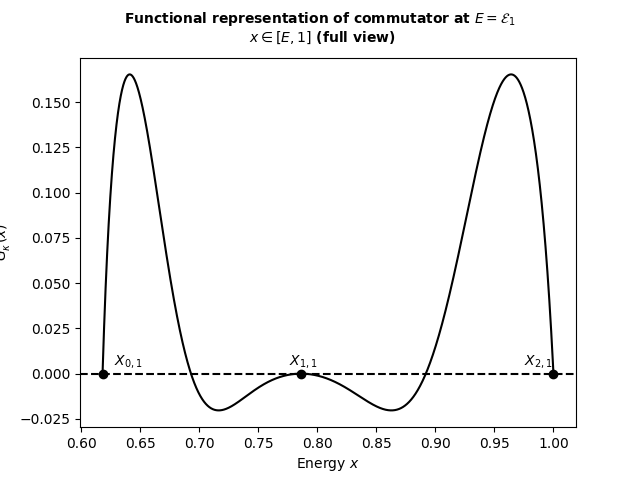}
\caption{$G_{\kappa=4} ^E(x)$, $x \in [E,1]$. Graphs at: $E \in \{\E_2, 0.59, \E_1 \}$, $\Sigma = \{4, 8, 12, 24\}$}
\label{fig:test_k4n2counter}
\end{figure}

Figure \ref{fig:test_k4n2counter} gives an idea of what $G_{\kappa} ^E$ looks like for a non-valid index set $\Sigma$. While $G_{\kappa}^E$ satisfies the interpolation constraints, it is not strictly positive.

\section{Conjecture for the interval $J_3(\kappa) := (\cos(2\pi / \kappa), \cos^2(\pi / \kappa))$}
\label{sec_Conjecture J3}

In this section we give some evidence for Conjecture \ref{conjecture_k344_2d}. We only do $\kappa=6,8$ in dimension 2. Note that this cumbersome phenomenon also happened for the standard Laplacian, see \cite[section 16]{GM3}.

For $\kappa=6$, the index set $\Sigma = \{6,12\}$, for which the coefficients $\rho_{6}$ and $\rho_{12}$ are given in section \ref{appli_scheme_k3}, also gives strict positivity on $(0.5024,0.672)$. This is an improvement compared to what was known previously (see Table \ref{tab:table1012sup}).

For $\kappa=8$, the index set $\Sigma = \{8,16\}$ gives strict positivity (using the coefficients of the linear interpolation) on $(0.70897,0.804)$. This is an improvement compared to what was known previously (see Table \ref{tab:table1012sup}). 

\section{The case of $\kappa=4$ in dimension 3}

\begin{figure}[H]
  \centering
  \includegraphics[scale=0.24]{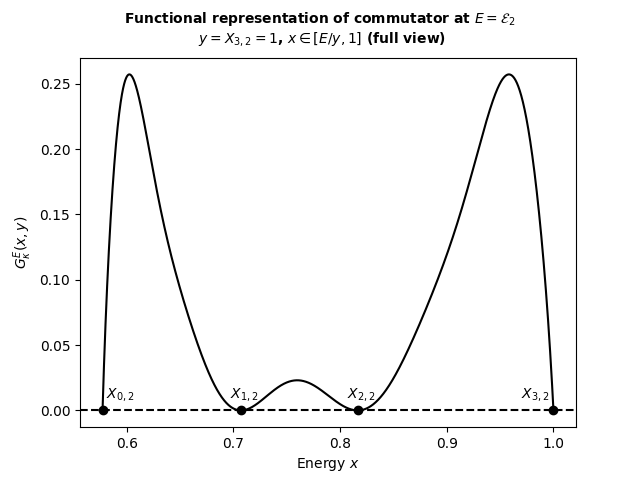}
   \includegraphics[scale=0.24]{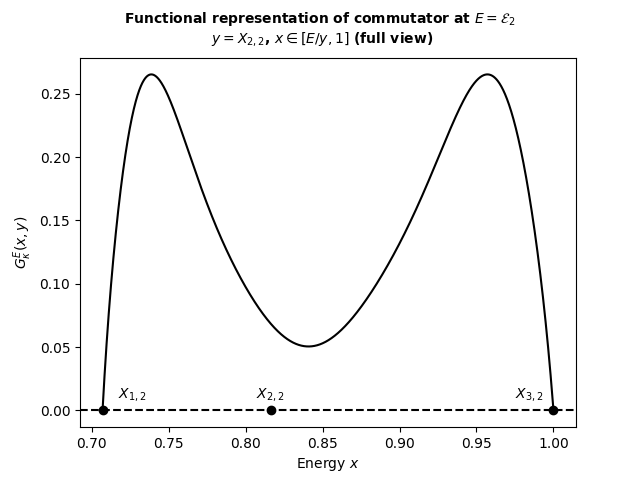}
    \includegraphics[scale=0.24]{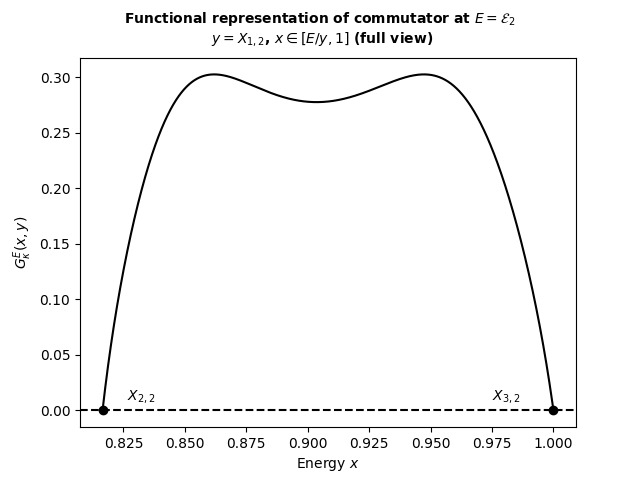}
 \includegraphics[scale=0.24]{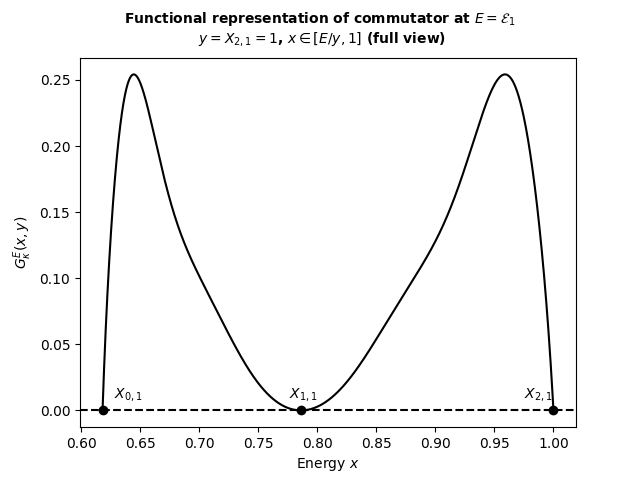}
  \includegraphics[scale=0.24]{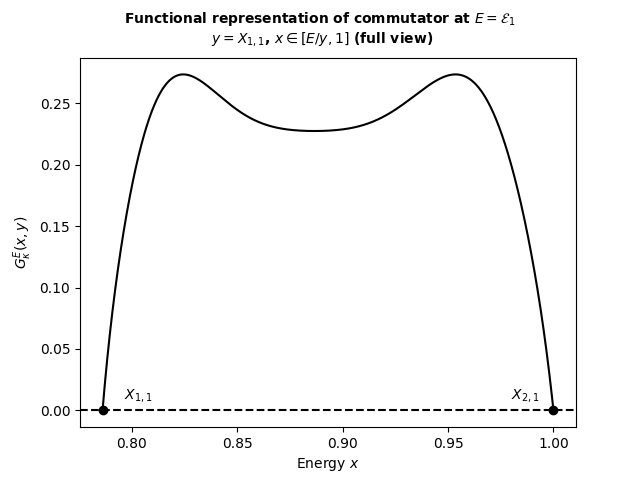}
  \includegraphics[scale=0.24]{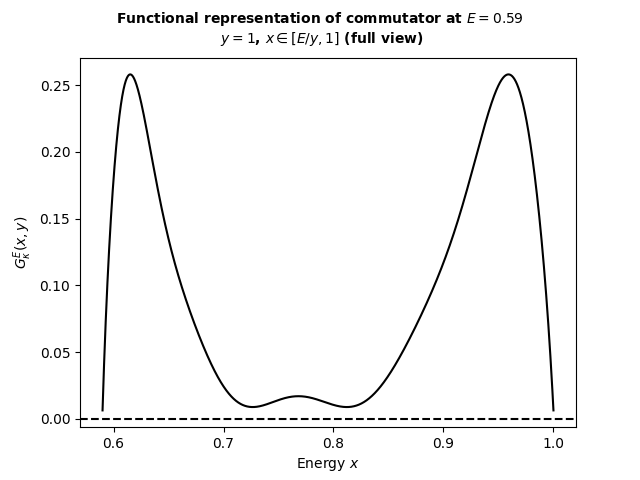}
    \includegraphics[scale=0.24]{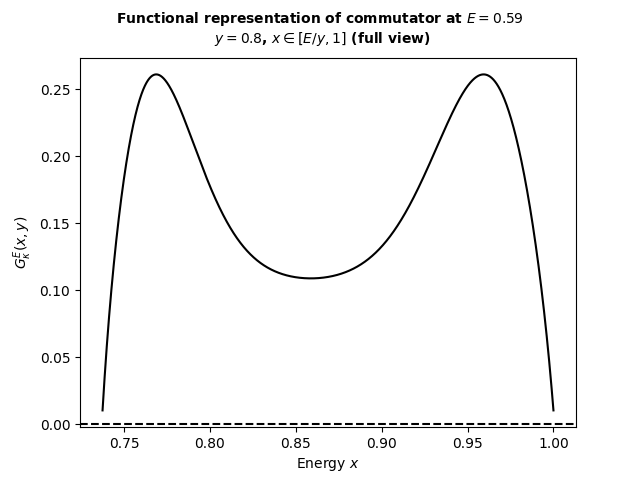}
     \includegraphics[scale=0.24]{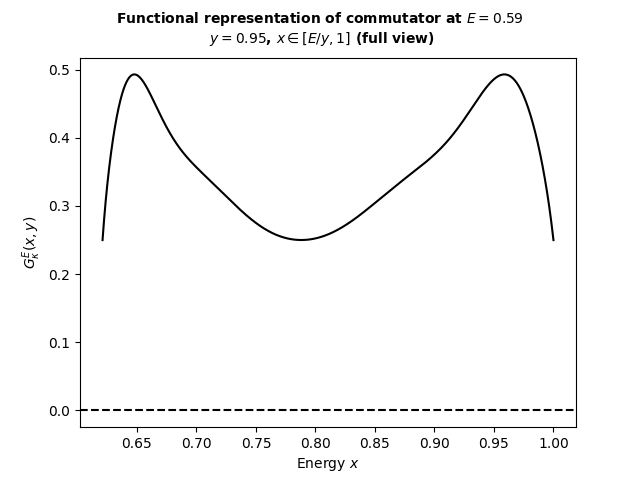}
\caption{$G_{\kappa=4} ^E(x)$, $x \in [E/y,1]$. Graphs at: $E=\E_2$, $E=0.59$, $E=\E_1$. $2^{nd}$ band}
\label{fig:test_k4n50}
\end{figure}

This section is in dimension 3. We illustrate the situation for $\kappa=4$, and the $2^{nd}$ band, namely $(\E_2, \E_1)$. In other words, $n=2$. We use the linear combination $\sum_j \rho_{j\kappa} A_{j\kappa}$ where the coefficients $\rho_{j\kappa}$ are the same as in dimension 2, i.e.\ the ones found in section \ref{appli_scheme_k2}. Figure \ref{fig:test_k4n50} shows the function $G_{\kappa=4} ^E (x,y)$ at $E \in \{ \E_2, 0.59, \E_1 \}$ and certain values of $y$. Note that $G_{\kappa}^E$ is $\geq 0$ at the band endpoints, but strictly positive at $E=0.59$. Of course, we have to check all the values of $y$ in the range $[-1,-|E|] \cup [|E|,1]$ to test the strict positivity of $G_{\kappa}^E$.


\section{Prior results for the Molchanov-Vainberg Laplacians} 
\label{appendix_mv}

Table \ref{tab:table1012sup} recalls the bands identified (numerically for the most part) in \cite{GM2} for the Molchanov-Vainberg Laplacian. These were obtained using the linear combination \eqref{LINEAR_combinationA} with $\rho_{j \kappa} = 1$ if $j=1$ and $0$ if $j>1$. Results for more values of $\kappa$ are listed in \cite[Tables XV and XVII]{GM2}.

\begin{table}[H]
\footnotesize
  \begin{center}
    \begin{tabular}{c|c|c} 
      $\kappa$ &  Intervals $\subset \boldsymbol{\mu} _{\kappa} (D)$. $d=2$.  &  Intervals $\subset \boldsymbol{\mu} _{\kappa} (D)$. $d=3$. \\ [0.5em]
      \hline
   \footnotesize   $2$  & $(0,1)$  &  $(0,1)$  \\ [0.5em]
    \footnotesize    $4$ &  $(0,0.5) \cup ( 0.7071 , 1 )$ & $(0,0.3535) \cup (0.7071,1)$ \\[0.5em]
     \footnotesize   $6$ &  $(0,0.25) \cup (0.5064,0.75) \cup (0.8660,1)$ & $(0,0.125) \cup (0.5148,0.6495) \cup (0.8660,1)$ \\[0.5em]
    \footnotesize    $8$ &  $(0,0.1464) \cup (0.3826,0.5) \cup (0.7121 ,0.8535) \cup (0.9238,1)$ &  $(0, 0.0560) \cup (0.7187,0.7885) \cup (0.9238,1)$  \\[0.5em]
    \end{tabular}
  \end{center}
    \caption{Sets $\subset \boldsymbol{\mu}_{\kappa} (D) \cap [0,1]$ found in \cite{GM2} using the trivial linear combination.}
        \label{tab:table1012sup}
\end{table}
We had conjectured exact expressions for the band endpoints in \cite{GM2}. Namely in dimension 2, we had conjectured the intervals in Table         \ref{tab:table1012sup} are in fact :
\begin{equation*}
\begin{aligned}
& \left(0,\sin^2(\pi / 4) \right) \cup \left( \sin( \pi /4) , 1 \right) & \kappa =4, \\
& \left(0,\sin^2(\pi / 6) \right) \cup \left( \ \_\_ \ , \sin^2(2\pi / 6) \right) \cup \left( \sin(2\pi / 6),1 \right) & \kappa =6, \\
& \left(0,\sin^2(\pi / 8) \right) \cup \left( \sin(\pi / 8),\sin^2(2\pi / 8) \right) \cup \left( \ \_\_ \  ,\sin^2(3\pi / 8) \right) \cup \left(\sin(3\pi / 8),1 \right) & \kappa =8.
\end{aligned}
\end{equation*}
Thanks to the identity 
$$\sin \left( \kappa/2 - j) \pi /\kappa \right) = \cos(j \pi / \kappa),$$
these conjectures can be reformulated in terms of cosines, which are more adapted to this paper:
 \begin{equation*}
\begin{aligned}
& \left(0,\cos^2(\pi / 4) \right) \cup \left( \cos( \pi /4) , 1 \right) & \kappa =4, \\
& \left(0,\cos^2(2\pi / 6) \right) \cup \left( \ \_\_ \ , \cos^2(\pi / 6) \right) \cup \left( \cos(\pi / 6),1 \right) & \kappa =6, \\
& \left(0,\cos^2(3\pi / 8) \right) \cup \left( \cos(3\pi / 8),\cos^2(2\pi / 8) \right) \cup \left( \ \_\_ \  ,\cos^2(\pi / 8) \right) \cup \left(\cos(\pi / 8),1 \right) & \kappa =8.
\end{aligned}
\end{equation*}
In dimension 3, we had conjectured :
\begin{equation*}
\begin{aligned}
& \left(0,\sin^3 \left( \pi /4 \right) \right) \cup \left(\sin \left(\pi /4 \right),1 \right) & \kappa =4, \\
& \left(0, \sin^3 \left( \pi / \kappa \right) \right) \cup \left( \ \_\_ \  , \sin^3 \left((\kappa / 2-1)\pi / \kappa \right)  \right) \cup \left(\sin \left( (\kappa /2-1)\pi / \kappa  \right) ,1\right) & \kappa =6,8.
\end{aligned}
\end{equation*}

\section{Appendix : Numerical / graphical algorithm to analyze the positivity of $G_{\kappa} ^E$}

In dimension $2$, we used the simple algorithm :
\begin{itemize}
\item For all $E \in [-1,1]$ : 
\begin{itemize}
\item let $y = E/x$
\item check if the function $x \mapsto G_{\kappa} ^E (x)$ has same sign on the interval $x \in \left[ -1, -|E|] \cup [|E|,1 \right]$.
\end{itemize}
\end{itemize}

In dimension $3$, we used the simple algorithm :
\begin{itemize}
\item For all $E \in [-1,1]$ : 
\item For all $y \in [-1,-|E|] \cup [|E|,1]$ :
\begin{itemize}
\item let $z = E/(xy)$
\item check if the function $x \mapsto G_{\kappa} ^E (x,y)$ has same sign on the interval $x \in \left[ -1, -|E/y|] \cup [ |E/y|,1 \right]$.
\end{itemize}
\end{itemize}

\end{document}